%% file: UC-waves-CIRM-course.tex
\documentclass[10pt]{article}
\usepackage[english,activeacute]{babel}
\usepackage{epsfig}
\usepackage[latin1]{inputenc}
\usepackage[T1]{fontenc}
\usepackage{stmaryrd}
\usepackage{amsmath,amsfonts,amssymb,mathrsfs,amsthm}
\usepackage{xcolor}
\usepackage{dsfont}
\usepackage{graphicx}
\usepackage[mathscr]{eucal}
\usepackage{makeidx}
\usepackage{verbatim}
\usepackage{graphics,graphicx}
\usepackage{textcomp}
\usepackage{float}
\usepackage[colorlinks=true]{hyperref}

\input mymacros.tex

\newcommand\bna{\begin{eqnarray*}}
\newcommand\ena{\end{eqnarray*}}

\newcommand\bnan{\begin{eqnarray}}
\newcommand\enan{\end{eqnarray}}

\newcommand\bnp{\begin{proof}}
\newcommand\enp{\end{proof}}

\newcommand\bneq{\begin{eqnarray*}\left\lbrace \begin{array}{rcl}}
\newcommand\eneq{\end{array} \right.\end{eqnarray*}}
\newcommand\bneqn{\begin{eqnarray}\left\lbrace \begin{array}{rcl}}
\newcommand\eneqn{\end{array} \right.\end{eqnarray}}

\newcommand\bni{\begin{itemize}}
\newcommand\eni{\end{itemize}}

\newcommand\nor[2]{\left\|#1\right\|_{#2}}
\newcommand{\pc}{ \usefont{T1}{cmtl}{m}{n} \selectfont}
\def\e{{\varepsilon}}

\def\x{{\bf x}}

\DeclareMathOperator{\differential}{Diff}
\newcommand{\dif}[1]{\differential_{\tau}^{#1}}

\newcommand{\difclas}[1]{\differential^{#1}}

\newcommand{\Ht}[1]{H^{#1}_{\tau}}

\title{Lectures on unique continuation for waves}
\author{Camille \textsc{Laurent}\footnote{CNRS, UMR 7598, Laboratoire Jacques-Louis Lions, F-75005, Paris, France} \footnote{UPMC Univ Paris 06, UMR 7598, Laboratoire Jacques-Louis Lions, F-75005, Paris, France, email: { \pc laurent@ann.jussieu.fr }}
and
Matthieu \textsc{L\'eautaud}\footnote{Laboratoire de Math\'ematiques d'Orsay, Universit\'e Paris-Saclay, CNRS, B\^atiment 307, 91405 Orsay Cedex France, email: { \pc  matthieu.leautaud@universite-paris-saclay.fr}}
}

\begin{document}
\maketitle

\begin{abstract}
These notes are intended as an introduction to the question of unique continuation for the wave operator, and some of its applications. The general question is whether a solution to a wave equation in a domain, vanishing on a subdomain has to vanish everywhere. We state and prove  two of the main  results in the field.
We first give a proof of the classical local H\"ormander theorem in this context which holds under a pseudoconvexity condition.  
We then specialize to the case of wave operators with time-independent coefficients and prove the Tataru theorem: local unique continuation holds across any non-characteristic hypersurface. This local result implies a global unique continuation statement which can be interpreted as a converse to finite propagation speed.
We finally give an application to approximate controllability, and present without proofs the associated quantitative estimates.
\end{abstract}

\tableofcontents
\section{Introduction and generalities}
\label{ch:intro-general}

These notes propose an introduction to the question of unique continuation for waves. We present the H\"ormander theorem and the Tataru theorem in this context, which are two of the main local results in the field. 
Before entering to heart of the subject, we discuss motivation for studying unique continuation for wave operators. Related references and further readings are presented in Section~\ref{notes}.

\subsection{Motivation and applications}
\label{s:motivation}
We start with presenting different applications to motivate the more technical parts of these notes. All these applications are discussed in Section~\ref{s:GLUC}.

\subsubsection{Penetration of waves into the shadow region}
\label{s:shadow-intro}
In this section, we consider the wave equation outside a convex obstacle in $\R^d$. Namely, let $\mathcal{O} \subset \R^d$ be a bounded smooth open subset, and consider $\M = \R^d \setminus \mathcal{O}$. We consider the Laplace operator $\Delta$ and $u(t,x)$ the solution to the wave equation
\begin{equation}
\label{e:onde}
\left\{
\begin{array}{rl}
\d_t^2u - \Delta u  = 0,  & \text{ on } \R \times \Int(\M) ,\\
                     u  = 0,  & \text{ on } \R \times \d \M , \\
 (u, \d_t u )|_{t=0} = (u_0, u_1),  & \text{ on } \M . \\
\end{array}
\right.
\end{equation}
Now, we consider a compact set $K \subset \M$, and assume that the initial data $(u_0, u_1)$ are supported in $K$. If the set $K$ is not too large, there is a whole region of $\M$ which does not intersect any ray of geometric optics in $\M$  (i.e. straight line in $\Int(\M)$, which reflects according to Snell-Descartes laws at the boundary $\d\M$) passing through $K$. Taking an open set $\omega$ in this shadow region, the question under consideration is the following: 
$$
\text{Can one recover } (u_0, u_1) \text{ from the observation of } u \text{ on the set } (-T,T)\times \omega ?
$$
And if so, what is the time $T$ required? By linearity of~\eqref{e:onde}, this can be reformulated under the following unique continuation question:
\begin{align}
\label{e:UCshadow}
\big( u \text{ solution to } \eqref{e:onde} , \quad \supp(u_0,u_1) \subset K,   \quad u|_{(-T,T)\times \omega} =0  \big)\implies (u_0,u_1)=0 ?
\end{align}
(and hence $u \equiv 0$ on account to the well-posedness of the Cauchy problem, see e.g.~\cite{Taylor:vol1} or~\cite{Evans:98}).
We shall see that Property~\eqref{e:UCshadow} is false if $T$ is too small, but holds true if $T$ is large enough. The limit time will be expressed as a natural geometric quantity.

\subsubsection{Approximate controllability for the wave equation}
\label{s:approx-wave-intro}
In this section, we consider a wave equation in a compact $d$-dimensional manifold $\M$ (or the closure of a bounded open set $\M \subset \R^d$), controlled from a subdomain. Namely, given $\chi_\omega \in C^\infty(\M)$, the equation
\begin{equation}
\label{e:controlled-onde}
\left\{
\begin{array}{rl}
\d_t^2u - \Delta_g u  = \chi_{\omega} f,  & \text{ on } (0,T) \times \Int(\M) ,\\
                     u  = 0,  & \text{ on } (0,T) \times \d \M , \\
 (u, \d_t u )|_{t=0} = (u_0, u_1),  & \text{ on } \M . \\
\end{array}
\right.
\end{equation}
The term $f$ in this equation plays the role of a forcing term, acting only on $\omega :=\{\chi_\omega \neq 0\}$. Controllability problems concern the ability of driving the solution $u$ to~\eqref{e:controlled-onde} from the initial state $(u_0, u_1)$ to a final target state $(v_0, v_1)$ at time $T$, using only the action of $f$ on $\omega$. This property depends a priori on the data/target, and is very complicated. More tractable questions, arising from applications in engineering are the following.

\begin{definition}
We say that~\eqref{e:controlled-onde} is {\em exactly controllable} from $(\chi_\omega, T)$ if for all data $(u_0, u_1) \in L^2(\M) \times H^{-1}(\M)$ and all target state $(v_0, v_1) \in L^2(\M) \times H^{-1}(\M)$, there is a function $f \in L^2((0,T) ; H^{-1}(\M))$ such that the solution to~\eqref{e:controlled-onde} satisfies $(u, \d_t u)|_{t=T} = (v_0, v_1)$. 

We say that~\eqref{e:controlled-onde} is {\em approximately controllable} from $(\chi_\omega, T)$ if for all data $(u_0, u_1) \in L^2(\M) \times H^{-1}(\M)$ all target state $(v_0, v_1) \in L^2(\M) \times H^{-1}(\M)$, and any precision $\eps>0$, there is a function $f=f_\eps \in L^2((0,T) ; H^{-1}(\M))$ such that the solution to~\eqref{e:controlled-onde} satisfies $\nor{(u, \d_t u)|_{t=T} - (v_0, v_1)}{L^2(\M) \times H^{-1}(\M)} \leq \eps$. 
\end{definition}
Notice that multiplication by $\chi_\omega$ maps $H^1_0(\M)$ into itself continuously and thus $H^{-1}(\M)$ into itself continuously as well; the Cauchy problem is thus well-defined in these spaces.

Due to finite speed of propagation for waves, if $\ovl{\omega}\neq \M$, a minimal time will be required for controllability to hold. 
Here, we will mostly be interested in the (weaker) approximate controllability question. Linearity of the equation shows it is enough to consider zero initial conditions $(u_0, u_1)=(0,0)$. Introducing the ``final value'' linear map 
$$
\begin{array}{rcl}
F_T :  L^2((0,T) ; H^{-1}(\M)) & \to & L^2(\M) \times H^{-1}(\M) \\
f & \mapsto & (u, \d_t u)|_{t=T} ,
\end{array}
$$
 where $u$ denotes the solution of~\eqref{e:controlled-onde} associated to $(u_0, u_1)=(0,0)$, approximate controllability is equivalent to $\range(F_T)$ being dense in $L^2(\M) \times H^{-1}(\M)$. This can be reformulated as $\ker(\transp{F}_T) = \{0\}$, where $\transp{F}_T$ is an appropriate transpose of $F_T$. Multiplying Equation~\eqref{e:controlled-onde} by $w$ solution to 
 \begin{equation}
\label{e:ondes-obs}
\left\{
\begin{array}{rl}
\d_t^2w - \Delta_g w  = 0,  & \text{ on } (0,T) \times \Int(\M) ,\\
                     w  = 0,  & \text{ on } (0,T) \times \d \M , \\
 (w, \d_t w )|_{t=T} = (w_0, w_1),  & \text{ on } \M , 
\end{array}
\right.
\end{equation}
 integrating on $(0,T)\times \M$, and integrating by parts in time and space, we obtain
 \begin{align*}
 \langle \d_t u(T) , w_0 \rangle_{H^{-1},H^1_0} -  ( u(T) , w_1 )_{L^2} = \langle f, \chi_\omega w \rangle_{L^2(0,T; H^{-1}), L^2(0,T;H^1_0)} .
  \end{align*}
As a consequence,  with an adequate choice of duality, one can identify $\transp{F}_T$ to the map
 $$
\begin{array}{rcl}
\transp{F}_T :   H^1_0 (\M) \times L^2(\M) & \to & L^2((0,T) ; H^1_0(\M))\\
(w_0 , w_1)& \mapsto &\chi_\omega w,
\end{array}
$$
where $w$ is the unique solution to~\eqref{e:ondes-obs}.
Again, $\ker(\transp{F}_T) = \{0\}$ is the unique-continuation property
\begin{equation}
\label{e:global-UC}
\big( w \text{ solution to } \eqref{e:ondes-obs} , \quad (w_0, w_1) \in H^1_0(\M) \times L^2(\M) \quad w|_{(0,T)\times \omega} =0  \big)\implies (w_0,w_1)=0 ,
\end{equation}
which now appears to characterize the approximate controllability of~\eqref{e:controlled-onde}.
 
\subsubsection{Inverse problems and the boundary control method}
\label{s:BC-method}
In this section, we still consider a wave equation in the compact manifold $\M$ with $\d\M\neq \emptyset$, with a time independent potential $q \in C^\infty(\M)$ and a forcing term $f \in C^\infty_c((0,\infty) \times \d\M)$ at the boundary and vanishing initial data:
\begin{equation}
\label{e:inverse-onde}
\left\{
\begin{array}{rl}
\d_t^2u - \Delta_g u + q u =0,  & \text{ on } (0,\infty) \times \Int(\M) ,\\
                     u  =f,  & \text{ on } (0,\infty) \times \d \M , \\
 (u, \d_t u )|_{t=0} = (0,0),  & \text{ on } \M . \\
\end{array}
\right.
\end{equation}
Existence and uniqueness of a solution $u \in C^\infty((0,\infty) \times \M)$ follows for instance from lifting the (smooth) boundary data to $\M$ and using well-posedness of~\eqref{e:controlled-onde}.
We define the dynamical Dirichlet-to-Neumann map 
$$
\Lambda_q : C^\infty_c((0,\infty) \times \d\M) \ni f \mapsto \d_nu|_{\d\M} \in C^\infty((0,\infty) \times \d\M) ,
$$
where $u$ is the solution to~\eqref{e:inverse-onde} and $\d_n$ denotes the unit outward normal derivative to $\d\M$.
A general question in inverse problem is whether the knowledge of the Dirichlet-to-Neumann map $\Lambda_q$ determines the potential $q$ uniquely?
That is to say, probing the wave in the domain $\M$ by means of a boundary source $f$ and knowing the response for all possible inputs $f$, is it possible to determine the potential $q$. In mathematical terms, do we have 
$\Lambda_{q_1}=\Lambda_{q_2} \implies q_1=q_2$?
We refer to~\cite{NO:22} for a presentation of the boundary control method of~\cite{Belishev:87} to solve this question. This methods relies in a key fashion on a unique continuation result presented in Section~\ref{s:UC-problem} below.

\subsection{Generalities about unique continuation}
\subsubsection{The unique continuation problem}
\label{s:intro-UC-problem}
All above described problems amount to a {\em unique continuation} property for the wave equation.
The general problem of  {\em unique continuation} can be set into the following form: given a differential operator $P=\sum_{|\alpha|\leq m} a_{\alpha}(x)\d_x^{\alpha}$ on an open set $\Omega \subset \R^n$, and given a small subset $U$ of $\Omega$, do we have (for $u$ regular enough):
\bnan
\label{e:UCP-general}
\left\{\begin{array}{rcl}
Pu &=& 0  \textnormal{ in  } \Omega,\\ u &=&0 \textnormal{ in  } U \end{array}\right.\Longrightarrow u= 0 \textnormal{ on } \Omega .
\enan
A more tractable problem than~\eqref{e:UCP-general} is the so called {\em local unique continuation across an hypersurface} problem: given an oriented local hypersurface $S=\left\{\Psi=0\right\}$ at a point $x_0$ (that is $\Psi(x_0) = 0$ and $d\Psi(x_0)\neq 0$), do we have the following implication:

There is a neighborhood $\Omega$ of $x_0$ so that 
\bnan
\label{e:UCP-general-local}
\left\{\begin{array}{rcl}
Pu &=& 0  \textnormal{ in  } \Omega,\\ u&=&0\textnormal{ in  }\Omega \cap S^+ \end{array}\right.\Longrightarrow u= 0 \textnormal{ in a neighborhood of } x_0 .
\enan 
Here $S^+=\left\{\Psi> 0\right\}$ is one side of $S$.
It turns out that proving~\eqref{e:UCP-general-local} for a suitable class of hypersurfaces $S$ (with regards to the operator $P$) is in general a key step in the proof of properties of the type~\eqref{e:UCP-general}.  

\bigskip
Let us discuss briefly  the local unique continuation property~\eqref{e:UCP-general-local} in the simple case where $P$ is a real non-degenerate vector fields. More precisely, consider $P$ a general real vector field (or, equivalently, first order ($m=1$) homogeneous differential operator) near $0$, that is $P = \sum_{k=0}^n a_k(x) \d_{x_k}$, with $a_k$ smooth and real-valued. Assume further that it is nondegenerate at $0$, that is $a(0) = (a_1(0) , \cdots , a_n(0)) \neq 0$. 
\begin{enumerate}
\item \label{i:lin-vec-field-1} (non-characteristic hypersurface) Take 
$S = \{\Psi = 0\}$ where $\Psi(0) = 0$ and $d\Psi(0)\neq 0$. Then, a sufficient condition for having local unique continuation~\eqref{e:UCP-general-local} is that $\langle d\Psi(0) , a(0)\rangle \neq 0$, i.e. the vector field $P$ is transversal to $S$ at $0$. This condition is a ``non-characteristicity assumption'', see Definition~\ref{def:char-non-char} below.
Note that the condition  $\langle d\Psi(0) , a(0)\rangle \neq 0$ is not necessary for unique continuation to hold, see the discussion in Example~\ref{i:lin-vec-field-2} below.
\item \label{i:lin-vec-field-2} (Constant vector fields and curved hypersurface) Here (as opposed to Example~\ref{i:lin-vec-field-1}), we shall see that the orientation of the hypersurface may play a role.
Consider for simplicity the operator $P=\frac{\d}{\d x_1}$  in $\R^2$ in a neighborhood of $0$, but the curved hypersurface $S = \{x=(x_1,x_2) \in \R^2 , \Psi(x)= 0\}$, where $\Psi(x_1, x_2) =x_2 - x_1^2$. Notice first that $S$ is tangent to $P$ at $0$ since $\langle d\Psi(0) , P \rangle = 0$.
We shall see that unique continuation holds from $S^- = \{\Psi <0\}$ (outside the parabola) to $S^+ = \{\Psi >0\}$ (inside the parabola), but not from $S^+$ to $S^-$.

Indeed solutions $u$ to $Pu=0$ write $u(x_1, x_2) = u^0( x_2)$ for all $x_1 \in \R$. 
The first statement then follows from the fact that any line $x_2 = cst>0$ intersects $S^+$ in a neighborhood of zero, thus showing that if $u^0( x_2)=0$ for all $x_1$ in a neighborhood of zero, then $u^0=0$.
Choosing $u^0 \in C^\infty_c(\R)$ such that $u^0(x_2) \neq 0$ on $0>x_2>-1$ and $u^0(x_2)=0$ on $x_2 \geq 0$ yields the second statement.
\end{enumerate}

The above examples~\ref{i:lin-vec-field-1}-\ref{i:lin-vec-field-2}  concerning first order partial differential operators (namely, vector fields) show that geometrical conditions linking the operator $P$ and the hypersurface $S$ are often needed for unique continuation to hold.
Let us now discuss related properties for the wave operator, for which the situation is far more difficult.

\subsubsection{Remarks on (non-)unique continuation for the flat/Minkowski wave operator}
\label{s:discussion-waves}

In this section, we collect known facts for the wave equation
\begin{equation}
\label{e:wave-flat}
(\partial_t^2-\Delta) u = 0 \quad \text{ on } \R\times \R^d ,
\end{equation}
 in the flat space $\R^{1+d}$, that are related to unique continuation questions.  

We start with the simpler case $d=1$ and consider the wave operator $P=\partial_t^2-\partial_x^2$ on $\R_t\times \R_x$. Then $P$ factorizes as $P = (\d_t + \d_x)(\d_t -\d_x)$ and all solutions to $Pu =0$ write 
 $u(t,x)=f(x+t)+g(x-t) + C_0 t + C_1 x + C_2$, where $f,g$ are functions and $C_j$ constants. Take for instance $g=0$, $C_j=0$ and $f \in C^\infty(\R)$ with $\supp(f) = [0,1]$. Then $u(t,x) = f(x+t)$ and the hypersurface $S=\left\{x+t=0\right\}$ thus does not satisfy the unique continuation property (at any point). More precisely, up to linear changes of variables, this problem reduces to that of examples~\ref{i:lin-vec-field-1}-\ref{i:lin-vec-field-2} discussed above, and one sees that the only hyperplanes $S$ not satisfying the unique continuation property (at any point) are $S_\pm^\alpha =\left\{x \pm t=\alpha \right\}$, for $\alpha \in \R$.

Let us now discuss the situation  in higher dimensions $d \geq 2$, which is radically different. This is linked with the fact that the polynomial $\xi_t^2 - \sum_{j=1}^d \xi_{x_j}^2$ does not factorize in a product of polynomials of degree $1$ and translates the fact that the values of solutions to~\eqref{e:wave-flat} are not ``transported''. To see this, we can actually solve the wave equation~\eqref{e:wave-flat}: for instance, in $\R^{1+3}$, the Kirchhoff formula (see e.g.~\cite{Evans:98})
\bnan
\label{e:kirchhoff}
u(t, x) = \frac{1}{4\pi t } \int_{|y-x|=t}u_1(y) dS_t (y) = \frac{t}{4\pi} \int_{\S^2}u_1^h(x - t \sigma) dS_1(\sigma) , \quad u(-t) = -u(t) ,\quad  t>0 
\enan
gives the unique solution to~\eqref{e:wave-flat} with $(u,\d_tu)|_{t=0} = (0,u_1)$, $u_1\in C^0(\R^3)$. In the first formula, the integration set is the ($2$ dimensional) sphere centered at $x$ and of radius $t$; in the second it is the unit sphere. The integration measure $dS$ is the surface measure on the sphere of radius $t$ (induced by the Euclidean measure $dx$ on $\R^3$). 
As a consequence of this explicit solution, we see that if we choose $u_1(x) = \chi(x)$ with $\chi\in C^\infty_c(\R^3)$, $\chi\geq 0$ and $\chi >0$ on $B(0,r)$, $r>0$ the associated solution $u$ is smooth and satisfies $u\geq0$ on $\R^{1+3}$. 
Moreover, notice that $u_1(x - t \sigma)= 0$ iff $x - t \sigma \notin B(0,r)$, and hence $u(t,x) = 0$ as soon as $t \S^2 \cap B(x, r) = \emptyset$. As a consequence, we have
\bnan
\label{e:suppsubset}
\supp ( u ) \cap \{t\geq 0\} = \{(t,x) \in \R^+ \times \R^3 , t-r \leq |x| \leq t + r \}   .
\enan
Several remarks are in order. 
The fact that the solution $u$ at time $t$ vanishes in the ball $|x| \leq t- r$ corresponds to the strong Huygens principle; this is strongly related to the fact that the dimension $3$ of $\R^3$ is odd, the metric is flat, and the wave operator has no lower order term.  
A contrario, the fact that the support of the solution at time $t$ is contained in the ball  $|x| \leq t +r$ translates the finite speed of propagation, discussed in more details in Theorem~\ref{t:wave-finite-speed} below.
Finally,~\eqref{e:suppsubset} also tells us that any point in the annulus $t-r \leq |x| \leq t + r $ is actually in the support of $u(t,\cdot)$.
This new piece of information is very important for what follows. It implies in particular that unique continuation cannot hold across an hyperplane tangent to the cone $|x|=t+r$.
We recall in this flat geometric setting the finite speed of propagation for waves.
\begin{theorem}[Finite speed of propagation for the wave equation]
\label{t:wave-finite-speed}
Let $u$ be a $C^2(\R^{1+d})$ (real-valued) solution of~\eqref{e:wave-flat}. 
 If $u|_{t=0}(x)=\d_t u|_{t=0}(x)=0$ for $|x|\leq r_0$, then $u=0$ in the cone 
\bna
C_{r_0}=\left\{(t,x)\in \R^{1+d}\textnormal{ s.t. }t\in [0,r_0] \text{ and }|x|\leq r_0-t\right\}.
\ena
\end{theorem}

We can infer an interesting consequence of Theorem~\ref{t:wave-finite-speed} concerning the unique continuation property for the wave operator:
unique continuation holds across the hypersurface $\{t=0\}$ and actually, we have some nice local linear quantification of the unique continuation. This situation is actually a particular case of a more general situation in which the differential operator $P$ (here $\d_t^2 - \Delta$) is said to be hyperbolic with respect to the hypersurface $S$ (here e.g. $\{t=0\}$). We refer e.g. to~\cite{Taylor:vol1} or~\cite{Evans:98} for more precisions and usual proofs of finite propagation speed (which is also a consequence of Theorem~\ref{thmreal2} below, see Section~\ref{s:examples-hor}). 

\subsubsection{Differential operators}
For later purposes, we give a definition of differential operators. 
Recall first that a function $f$ on $\R^n$ is said homogeneous of degree $m>0$ if 
$$
f(\lambda \xi) = \lambda^m f(\xi) , \quad \text{ for all }\lambda>0 \text{ and }\xi \in \R^n.
$$

\begin{definition}[Classical differential operators]
\label{defdiffclassic}
Let $\Omega \subset \R^n$ be an open set and $m \in \N$.
  We say that $P$ is a (linear) differential operator of order $m$ on $\Omega$ if there are coefficients $a_\alpha \in C^\infty(\Omega)$ having all derivatives bounded uniformly on $\Omega$, such that 
$P = \sum_{|\alpha|\leq m} a_\alpha(x) D^\alpha$ with $m = \max\{|\alpha| , a_\alpha \neq 0\}$. 
  We denote $\difclas{m}(\Omega)$ the set of differential operators of order $m$ on $\Omega$ (in the class $\difclas{m}(\Omega)$).
  We say that the function $p_m(x,\xi) = \sum_{|\alpha|= m} a_\alpha(x) \xi^\alpha$ is the principal symbol of $P$.  It is a {\em homogeneous} polynomial of degree $m$ in the variable $\xi$.
 \end{definition}
 \begin{example}
If $\Omega \subset \R^n$ and $a,b \in C^\infty(\Omega;\C)$ have all derivatives bounded uniformly on $\Omega$, then $a(x)D_j + b(x) \in \difclas{1}(\Omega)$ with principal symbol $p_1(x,\xi) =a(x)\xi_j$, and $-\Delta + a(x)D_j + b(x) \in \difclas{2}(\Omega)$ with principal symbol $p_2(x,\xi)=|\xi|^2$.
\end{example}
The augmented set $\Omega \times \R^n$, in which the principal symbol  $ p_m$ lives, may be seen as a ``phase space'' containing both the position variable $x$ and the Fourier/frequency/momentum variable $\xi\in \R^n$. The latter is to be understood as a cotangent variable $\xi \in T^*_x\Omega$, as we shall see below. 

\subsubsection{A general local unique continuation result in the analytic Category}

The first general unique continuation result of the form~\eqref{e:UCP-general-local} is the Holmgren-John Theorem, stating that, for operators with analytic coefficients, unique continuation holds across any noncharacteristic hypersurface $S$. 

\begin{definition}
\label{def:char-non-char}
Let $P$ be a differential operator of order $m$ on $\Omega$, $x_0 \in \Omega$ and $S$ a local hypersurface passing through $x_0$, that is $S= \{\Psi=0\}$, $\Psi(x_0)=0$ and $d\Psi(x_0)\neq 0$ with $\Psi\in C^1(\Omega)$. We say that $S$ is characteristic (resp. non-characteristic) for $P$ at $x_0$ if $p_m(x_0, d\Psi(x_0)) = 0$ (resp. $p_m(x_0, d\Psi(x_0)) \neq 0$).
\end{definition}

Also, given a local hypersurface $S= \{\Psi=0\}$, it has locally two sides which we write 
$$
S^\pm  = \left\{x\in\Omega; \pm \Psi(x)>  0\right\}.
$$
\begin{theorem}[Holmgren-John Theorem]
\label{thmholmgren}
Let $P$ be a differential operator of order $m$ on $\Omega$, having all coefficients {\em real analytic} in a neighborhood of $x_0 \in \Omega$ and $S \ni x_0$ being a local hypersurface.
Assume that $S$ is non characteristic for $P$ at $x_0$.
Then, there exists a neighborhood $V$ of $x_0$ so that every $u\in \mathcal{D}'(\Omega)$ satisfying $Pu=0$ on $\Omega$ and $u=0$ in the set $S^+$ vanishes identically in $V$.
\end{theorem}
Another (slightly weaker) way of writing the conclusion is to say that $x_0 \notin \supp(u)$. We refer e.g. to~\cite[Theorem~5.3.1]{Hoermander:63} for a proof of Theorem~\ref{thmholmgren}.
Note that this unique continuation property does not take into account the orientation of the hypersurface $S$, i.e. it holds from $S^+$ to $S^-$ as well as from $S^-$ to $S^+$.

The non-characteristicity condition is very weak, and in some sense optimal. Indeed, we saw in Example~\ref{i:lin-vec-field-1} in Section~\ref{s:intro-UC-problem} for linear vector-fields that unique continuation holds for non-characteristic hypersurfaces, and does not hold for some characteristic hypersurfaces.  
We also saw in Section~\ref{s:discussion-waves} for the wave operator that local uniqueness does not hold across some hypersurfaces that are tangent to the cone $|x|=t+r$. These are precisely characteristic hypersurfaces: the principal symbol of the wave operator $\d_t^2-\Delta$ is given by $p_2(t,x,\xi_t,\xi_x) = - \xi_t^2 + |\xi_x|^2$, and a hypersurface $\{\Psi(t,x) = 0\}$ tangent to $\{|x|=t+r\}$ at the point $(t_0,x_0)$ has $|\d_t\Psi(t_0,x_0)| = |d_x\Psi(t_0,x_0)|$.
Remark however that the non-characteristicity condition is a ``first order condition'': it only cares about the tangent space of the hypersurface. We saw in Example~\ref{i:lin-vec-field-2} in Section~\ref{s:intro-UC-problem} in the case of first order differential operators a more subtle ``second order condition'' (curvature condition) on the hypersurface that  may yield unique continuation across a characteristic hypersurface. This is linked to the so-called pseudoconvexity condition (see e.g. Definition~\ref{defconvexsurface} below).

We recall that a function $f : \Omega \subset \R^n \to \C$  is {\em real analytic} if for every $y\in \Omega$, there is a convergence radius $R>0$ and coefficients $a_\alpha \in \C^n$, $\alpha\in \N^n$ such that 
$$
f(x) = \sum_{\alpha \in \N^n} a_\alpha (x-y)^\alpha = \sum_{\alpha_1,\cdots,\alpha_n \in \N} a_\alpha (x_1-y_1)^{\alpha_1}\cdots (x_n-y_n)^{\alpha_n} , \quad \text{ for all } x\in B(y,R)\subset \Omega ,
$$
where the series is absolutely convergent. For every compact set $K\Subset \Omega\subset \R^n$, such a function $f$ can be extended to a complex neighborhood of $K$ in $\C^n$ as a complex analytic function. Analyticity is a very demanding regularity assumption. In Theorem~\ref{thmholmgren}, we stress that {\em all} the coefficients of $P$ should have this regularity. In most situations, however, this requirement is much too strong.
As an example, even for the wave equation on a flat (and hence analytic) metric, this theorem does not allow for the addition of a $C^\infty$ time independent potential $V(x)$. This is a very strong drawback to the result.
Therefore, we would like to avoid the analyticity assumption on the coefficients. This will lead to consider stronger geometric assumption, of convexity type (see e.g. Definition~\ref{defconvexsurface} below) and will be the object of Chapter~\ref{chapterclassical}. Then Chapter~\ref{chapterwave} will deal with an intermediate case where the analyticity is with respect to only one variable (we will actually treat the simpler case where it is independent on one variable).

\subsubsection{The general strategy of Carleman}
\label{subsectstratCarle}
We consider here $\Omega$ a bounded open subset of $\R^n$, $P$ a differential operator on $\Omega$, $x_0 \in \Omega$ a point, and a hypersurface $S=\left\{\Psi=0\right\}$ containing $x_0$.
We aim at proving local unique continuation for an operator $P$ across the hypersurface $S=\left\{\Psi=0\right\}$ (say, a statement like~\eqref{e:UCP-general-local}).
In particular, we want to prevent the situation in which a smooth function $w$ both solves $Pw=0$ and vanishes (possibly ``flately'', in the sense that all its derivatives vanish) on $S$.
We thus need to ``emphasize'' the local behavior of functions close to the hypersurface $S$. 

The general idea of Carleman to do so, and thus prove unique continuation, is to consider weighted estimates of the form
\bnan
\label{Carlintro}
 \nor{e^{\tau\Phi} w}{L^2(\Omega)} \leq C \nor{e^{\tau\Phi}P w}{L^2(\Omega)} , 
\enan
which hold: 
\begin{itemize}
\item for some well-chosen weight function $\Phi : \ovl{\Omega}\to \R$ (related to $\Psi$ as discussed below);
\item for all $w \in C^{\infty}_c(\Omega)$ (related to $u$ as discussed below);
\item and {\em uniformly} for $\tau$ sufficiently large, i.e. $\tau \geq \tau_0$.
\end{itemize}

To prove the relevance/efficiency of this approach, two different things need to be explained:
\begin{enumerate}
\item \label{Carl-imp-PU} what is the link between Carleman estimates like~\eqref{Carlintro} and unique continuation properties like~\eqref{e:UCP-general}?
\item \label{prove-carl} how to prove such Carleman estimates?
\end{enumerate}

Let us first discuss point~\ref{Carl-imp-PU}.
Note first that~\eqref{Carlintro} says directly that if $w \in C^{\infty}_c(\Omega)$ is solution of $Pw =0$ on $\left\{\Phi\geq 0\right\}$, then the right hand side will tend to zero as $\tau$ tends to infinity. Therefore, the left hand side will converge to zero, which implies that $w$ is supported in $\left\{\Phi\leq 0\right\}$.

However, statements like~\eqref{e:UCP-general-local} that are useful in applications are not concerned with functions $w$ having compact support. Moreover, in general, as we shall see, usual differential operators $P$ do not admit solutions $w$ to $Pw=0$ having compact support!

The heart of the Carleman method to pass from the estimate~\eqref{Carlintro} to the unique continuation statement~\eqref{e:UCP-general-local} resides in applying~\eqref{Carlintro} to $w= \chi u$, where $u$ is the function for which unique continuation has to be proved (hence solving $Pu=0$ in $\Omega$ and $u=0$ on $\Psi\geq 0$), and $\chi \in C^\infty_c(\Omega)$ is a cut-off function (to be chosen) allowing to apply~\eqref{Carlintro}. 

Using that $P \chi u = \chi P  u + [P ,\chi] u =[P ,\chi] u$ (where $[P ,\chi]$ denotes the commutator of $P$ and the multiplication operator by $\chi$), this then yields
 $$
 \nor{e^{\tau\Phi} \chi u}{L^2(\Omega)} \leq C \nor{e^{\tau\Phi}[P ,\chi] u}{L^2(\Omega)} .
$$
We then notice that $\supp [P ,\chi] \subset \supp \nabla \chi$. If we now {\em assume} (this can be achieved if $\Phi$ is a slight convexification of $\Psi$),
that the functions $\Psi, \Phi, \chi$ are chosen such that $\supp(\nabla \chi)\cap \{\Psi \leq 0\} \subset \{\Phi \leq -\eta\}$, for some $\eta>0$ (small!), then the support property of $u$ (namely $u=0$ on $\Psi\geq 0$) implies that $\supp ( [P ,\chi] u ) \subset \{\Phi \leq -\eta\}$, and we thus obtain
 $$
 \nor{e^{\tau\Phi} \chi u}{L^2(\Omega)} \leq C_u e^{-\eta \tau} , \quad \text{ for all } \tau \geq \tau_0 .
$$
The following lemma then implies that $\chi u$ vanishes identically in $\{\Phi \geq-\eta\}$ which contains a neighborhood of the point $x_0$. 
\begin{lemma}
\label{l:uniqueness-from-carleman}
Assume $w \in L^2(\Omega)$ satisfies $\nor{e^{\tau\Phi}w}{L^2(\Omega)} \leq C e^{-\eta \tau}$ for all $\tau \geq \tau_0$. Then we have $w=0$ a.e. on $\{\Phi \geq-\eta\}$.
\end{lemma}
The proof of the lemma reduces first to the case $\eta=0$ by changing $\Phi$ in $\Phi+ \eta$. Then, it suffices to notice that if $w$ does not vanish a.e. on $\{\Phi > 0\}$, there are $\eps>0$ and a compact set $E \subset \{\Phi > 0\}$ of positive measure such that $|w| \geq \eps>0$ a.e. on $E$. This yields 
$$
C^2 \geq \nor{e^{\tau\Phi}w}{L^2(\Omega)}^2 \geq \int_E e^{2\tau \Phi} |w|^2 \geq \eps^2 \int_E e^{2\tau \min_E \Phi}  = \eps^2 |E| e^{2\tau \min_E \Phi}   \to_{\tau \to + \infty} + \infty ,
$$
and hence a contradiction.

To conclude, this brief discussion of point~\ref{Carl-imp-PU} suggests that unique continuation~\eqref{e:UCP-general-local} will hold (across $\{\Psi=0\}$) provided the Carleman estimate~\eqref{Carlintro} is true for some weight function $\Phi$ satisfying an appropriate {\em geometric convexity condition}. 

\bigskip
As stated in point~\ref{prove-carl}, the other issue is how to prove Carleman estimates, and, in particular, understand the conditions on $\Phi$ for which~\ref{prove-carl} can hold. As far as this analysis is concerned, the exponential weight is not convenient to work with. One might thus want to eliminate it by setting $v=e^{\tau \Phi} w$. Then~\eqref{Carlintro} is equivalent to $\nor{v}{L^2(\Omega)}\leq C \nor{P_{\Phi} v}{L^2(\Omega)}$, with $P_{\Phi}=e^{\tau\Phi}Pe^{-\tau\Phi}$ is the so-called conjugated operator. 
Note that again here, we slightly abuse notation and make the confusion between the function $e^{\tau\Phi}$ and the operator of multiplication by $e^{\tau\Phi}$.
We are thus left to prove a lower bound for the operator $P_\Phi$.

Writing $\d_j(e^{-\tau\Phi} u) = e^{-\tau\Phi}(\d_j u -\tau u \d_j\Phi)$ implies that 
\bnan
\label{e:conj-Dj}
e^{\tau\Phi}D_je^{-\tau\Phi}=D_j+i\tau \partial_j \Phi.
\enan 
The first effect of conjugation is that there is no exponential factor in the right-handside, which is much more convenient.
Second, the conjugation changes $D_j$ into an operator having one derivative and one exponent of $\tau$. 
We thus expect (and we will check) that for general differential operators $P=\sum_{\alpha}a_{\alpha}(x)D^{\alpha}$, the associated conjugated operator $P_{\Phi}$ will have as many derivatives as exponents of $\tau$. Since we want to obtain estimates that are uniform for large $\tau$, we have to think of $\tau$ as having the same weight as a derivative. We describe this calculus in the next section.

\subsection{Operators depending on a large parameter $\tau$}
\label{s:operators-dif-tau}
In this section, we describe the setting in which Carleman estimates like~\eqref{Carlintro} shall be proved (see Chapters~\ref{chapterclassical} and~\ref{chapterwave} below).
\subsubsection{Differential operators depending on a large parameter $\tau$}
We discuss  the calculus for differential operators depending on a large parameter $\tau$.
One may think  to $\tau$ as having the same weight as a derivative, i.e. as the Fourier variable $\xi$.
Since $\tau$ is aimed at being large, we will always assume $\tau\geq 1$ when dealing with estimates uniform in $\tau$.
We first define the $H^s_\tau$ norm of a function $u \in \calS(\R^n)$ as 
\bna
\nor{u}{\Ht{s}}=\nor{(|D|^2+\tau^2)^{\frac{s}{2}}u}{L^2(\R^n)}=(2\pi)^{-n/2}\nor{(|\xi|^2+\tau^2)^{\frac{s}{2}}\widehat{u}}{L^2(\R^n)}.
\ena
 Note that for fixed $\tau$, this norm is equivalent to the usual $H^s$ norm, since, for $\tau \geq 1$, we have
 $$
 |\xi|^2 + 1 \leq |\xi|^2+\tau^2 \leq \tau^2( |\xi|^2 + 1) .
 $$
 That is to say that $ \nor{u}{H^s(\R^n)} \leq \nor{u}{\Ht{s}} \leq \tau^s \nor{u}{H^s(\R^n)}$ for all $\tau \geq 1$. This is however not uniform as $\tau \to +\infty$.
Note also that, as for usual Sobolev spaces, the definition of the $\Ht{s}$ norm has a uniformly equivalent definition in case $s=k \in \N$.
\begin{lemma}
\label{l:equiv-norms-tau}
Let $k\in \N$. Then, there is $C>1$ such that for all $\tau\geq 1$ and all $u \in H^k(\R^n)$, we have
\bna
C^{-1} \nor{u}{H^k_\tau(\R^n)}^2 \leq  \sum_{|\alpha|+ \beta \leq k} \tau^{2\beta} \nor{ \d^\alpha u}{L^2(\R^n)}^2 \leq C \nor{u}{H^k_\tau(\R^n)}^2 .
\ena
\end{lemma}
In particular, we often use the case $k=1$ for which 
\bna
\nor{u}{\Ht{1}}\approx \nor{u}{H^1}+\tau \nor{u}{L^2}, \quad \text{ uniformly for } \tau \geq 1 .
\ena
\begin{definition}[Differential operators depending on $\tau$]
\label{defdiffclassic-tau}
Let $m\in \N$ and $\Omega \subset \R^n$ an open set. We denote $\dif{m}(\Omega)$ the set of differential operators of the form $P=\sum_{|\alpha|+\beta\leq m}p_{\alpha,\beta}(x)\tau^{\beta}D^{\alpha}$ with $p_{\alpha,\beta} \in C^\infty(\Omega)$ such that all derivatives of $p_{\alpha,\beta}$ are bounded uniformly on $\Omega$. 
We say that $p_m(x,\xi,\tau)=\sum_{|\alpha|+\beta= m}p_{\alpha,\beta}(x)\tau^{\beta}\xi^{\alpha}$ is its principal symbol (in the class $\dif{m}(\Omega)$). It is homogeneous of degree $m$ in $(\xi,\tau)$, in the sense that
$$
p_m(x, \lambda \xi , \lambda \tau) = \lambda^m p_m(x, \xi , \tau) , \quad \text{ for all } x\in \Omega ,\xi \in \R^n , \tau \geq 0 \text{ and } \lambda >0 .
 $$
\end{definition}
Recall that, the order $m$ being fixed, smooth homogeneous functions of degree $m$ in this sense identify (through the restriction map) to smooth functions on the half-sphere bundle over $\Omega$, namely
$$
\{(x, \xi , \tau)  \in  \Omega \times \R^n \times \R^+ ,  |\xi|^2 + \tau^2 = 1\} .
$$
Remark that  Definition~\ref{defdiffclassic-tau} is almost the same definition as Definition~\ref{defdiffclassic}, except for the dependence on  $\tau$ which changes the definition of the principal symbol.

\bigskip
Now, we  describe the calculus of differential operators with a large parameter. This consists in explaining the properties of such operators with respect to usual operations (composition, commutators, taking the adjoint), their mapping properties (in $\tau$ dependent Sobolev spaces) and positivity properties. Moreover, we link such properties with those of the symbol of the operators. The philosophy is that we recover certain properties of the operators only from their principal symbols (which are simpler objects to manipulate, namely functions on the augmented space $\Omega_x \times \R^n_\xi \times \R^+_\tau$).
The general Heuristic is that a these differential operators act as if they were multiplication by $p_m(x,\xi,\tau)$, modulo lower order terms.
A rough summary of the calculus properties of operators of $\dif{m}$ is as follows. 
\begin{proposition}[Symbolic calculus for differential operators]
\label{p:calcul}
Let $\Omega \subset \R^n$ be an open set, $m, m_1, m_2 \in \N$ and  $P\in \dif{m}(\Omega)$, $A \in \dif{m_1}(\Omega)$, $B\in \dif{m_2}(\Omega)$ having respective principal symbol $p$, $a$, $b$, then
\begin{enumerate}
\item \label{i:composition} (composition) $AB= A\circ B \in \dif{m_1+m_2}(\Omega)$ with principal symbol $ab$;
\item \label{i:commutator} (commutators) $[A,B]=AB-BA\in \dif{m_1+m_2-1}(\Omega)$ is of order $m_1+m_2-1$ with principal symbol $\frac{1}{i}\left\{a,b\right\}$, where
\begin{equation}
\label{d:def-poisson-bracket}
\left\{a,b\right\} :=  \d_\xi a \cdot \d_x b - \d_x a \cdot \d_\xi b = \sum_{j=1}^n \left(\d_{\xi_j} a \d_{x_j} b - \d_{x_j} a \d_{\xi_j} b \right) ,
 \end{equation}
is the Poisson bracket.
\item \label{i:propadjoint} (adjoint) $P^*$, the formal adjoint in $L^2$ (tested with functions in $C^\infty_c(\Omega)$) belongs to $\dif{m}(\Omega)$ with principal symbol $\overline{p}$;
\item \label{i:propSob} (action on Sobolev spaces) 
if $\Omega =\R^n$, $P$ maps continuously $H^{s}_\tau$ into $H^{s-m}_\tau$ for all $s\in \R$, uniformly for $\tau \geq 1$. 
\end{enumerate}
\end{proposition}
Proofs are elementary as operators in $\dif{m}$ are linear combinations of mononomials $\tau^\beta a(x)D^\alpha$, $\beta \in \N, \alpha \in \N^n$ for which the properties can be checked by hand.
We refer e.g. to~\cite{LL:book} for direct elementary proofs, and to~\cite{Hoermander:V3} or~\cite{LRLR:book1} for the whole machinery of pseudodifferential calculus (with a large parameter).

\subsubsection{The conjugated operator}
As described in Section~\ref{subsectstratCarle}, the introduction of the calculus with the large parameter $\tau$ is motivated by the conjugated operator $P_{\Phi}:=e^{\tau\Phi}Pe^{-\tau\Phi}$.
We prove here that it belongs to the class $\dif{m}$, and compute its principal symbol.

\begin{lemma}[The conjugated operator]
\label{lmsympphi}
Let $P=\sum_{|\alpha|\leq m}p_{\alpha}(x)D^{\alpha}\in \difclas{m}(\Omega)$ be a (classical) differential operator with principal symbol $p_m$ and let  $\Phi\in C^{\infty}(\Omega)$ be real-valued and bounded as well as all its derivatives.
Then, the operator $P_{\Phi} $ defined by $P_{\Phi}v =e^{\tau \Phi}P(e^{-\tau \Phi }v)$ satisfies $P_{\Phi} \in \dif{m}(\Omega)$, and its principal symbol, denoted by $p_{\Phi} = p_{\Phi,m}$ (with a slight abuse of notation), is given by
\bna
p_{\Phi}(x,\xi, \tau)=p_m(x,\xi+i\tau d \Phi(x))=\sum_{|\alpha| = m} p_{\alpha}(x)(\xi+i\tau d\Phi(x))^{\alpha} .
\ena
\end{lemma}
Roughly speaking, Lemma~\ref{lmsympphi} says that $p_{\Phi}$ is obtained by replacing $\xi$ by $\xi+i\tau d\Phi(x)$ in $p_m$.
Note that it implies in particular that $p_{\Phi}$ has a complex-valued symbol if $p_m$ is real-valued: the conjugation turns selfadjoint operators into non-selfadjoint ones.

\bnp
As already checked in~\eqref{e:conj-Dj}, we have $e^{\tau\Phi}D_j(e^{-\tau\Phi}u)=D_j u+i\tau (\partial_j \Phi) u$. In particular, the conjugated operator $e^{\tau\Phi}D_je^{-\tau\Phi}$ lies in the class $\dif{1}$ with principal symbol $\xi_j+i\tau \partial_j \Phi$.
We now write  
\begin{align*}
e^{\tau\Phi}D_j^{\alpha_j}e^{-\tau\Phi}&= (e^{\tau\Phi}D_je^{-\tau\Phi})(e^{\tau\Phi}D_je^{-\tau\Phi})\cdots (e^{\tau\Phi}D_j e^{-\tau\Phi})  \qquad (\alpha_j\text{ times}).
\end{align*}
Therefore, using Proposition~\ref{p:calcul} $\alpha_j-1$ times, we obtain that this is a differential operator depending on $\tau$ of order $\alpha_j$ with principal symbol $(\xi_j+i\tau \partial_j \Phi)^{\alpha_j}$.
Since $D^{\alpha}= D_1^{\alpha_1}\cdots D_j^{\alpha_j}\cdots D_n^{\alpha_n}$, we obtain similarly that $e^{\tau\Phi}D^{\alpha}e^{-\tau\Phi} \in \dif{|\alpha|}$, with principal symbol
\bna
\prod_{j=1}^n (\xi_j+i\tau \partial_j \Phi)^{\alpha_j}= (\xi+i\tau d \Phi)^{\alpha}  .
\ena
Since $p_{\alpha}$ commutes with $e^{\tau\Phi}$ and $P=\sum_{\alpha}p_{\alpha}(x)D^{\alpha}$, this proves the sought result.
\enp
\begin{example}[second order operators with real-valued principal symbol]
\label{e:general-operator-2}
In these notes, we are particularly interested in {\em second order differential operators with real-valued principal symbol} (and in particular wave operators), namely $P \in \difclas{2}(\Omega)$ with $p_2$ real-valued. The principal symbol of such operators write $p_2(x,\xi)= \sum_{i,j=1}^n a^{ij}(x) \xi_i \xi_j$ with real coefficients $a^{ij}$.
This encompasses of course the case of the Laplace operator and the wave operator.
Notice first that $\xi \mapsto p_2(x,\xi)= \sum_{i,j=1}^n a^{ij}(x) \xi_i \xi_j$ is a real quadratic form for all $x\in \Omega$. In particular, we have the canonical polar form: 
$$
p_2(x,\xi)= \sum_{i,j=1}^n a^{ij}(x) \xi_i \xi_j  = \sum_{i,j=1}^n \frac12 \left( a^{ij}(x) + a^{ji}(x) \right) \xi_i \xi_j ,
$$
and we can thus assume that 
\bnan
\label{e:symmetric-aij}
\text{the matrix } (a^{ij}(x))_{i,j} \text{ is symmetric, i.e. } a^{ij}(x) = a^{ji}(x) \text{ for all } 1\leq i,j\leq n .
\enan
Concerning the operator $P$, we then have 
$$
P - p_2(x,D) \in \difclas{1}(\Omega) , \quad p_2(x,D) = \sum_{i,j=1}^n a^{ij}(x) D_i D_j = \sum_{i,j=1}^n D_i a^{ij}(x) D_j  + R_1 ,
$$
where $R_1 = -  \sum_{i,j=1}^n D_i (a^{ij}) D_j \in  \difclas{1}(\Omega)$.
The operator $\sum_{i,j=1}^n D_i a^{ij}(x) D_j$ is formally (i.e. tested with functions in $C^\infty_c(\Omega)$) selfadjoint in $L^2(\Omega,dx)$ (equivalently, one says that it is of divergence form with respect to the measure $dx$).  
This last form thus states in a clearer way that the operator is formally self-adjoint modulo $\difclas{1}(\Omega)$.
Also, $\xi \mapsto p_2(x,\xi)= \sum_{i,j=1}^n a^{ij}(x) \xi_i \xi_j$ is a real quadratic form with~\eqref{e:symmetric-aij}, and thus Lemma~\ref{lmsympphi} states that the principal symbol of the associated conjugated operator $P_\Phi$ is given by
\begin{equation*}
p_{\Phi}(x, \xi,\tau) = p_2(x, \xi + i \tau d\Phi(x)) = p_2(x, \xi ) - \tau^2 p_2(x, d\Phi(x)) +2 i \tau \widetilde{p_2}(x, \xi, d\Phi(x)) ,
\end{equation*}
where $\widetilde{p_2}(x, \xi, \eta) = \sum_{i,j=1}^n a^{ij}(x) \xi_i \eta_j$ is the polar bilinear form of the quadratic form $p_2(x,\xi)$.
\end{example}

\subsubsection{A  G\aa rding inequality for a class of operators with a large parameter}
In this section, we prove that operators having a real positive principal symbol are positive (referred to as a G\aa rding inequality).
However, for the need of Carleman estimates, the class of differential operators is not quite sufficient. We need to consider a slightly larger class, that also includes the operator 
$$(-\Delta+\tau^2)^{-1} = (|D|^2+\tau^2)^{-1},  \quad \tau\geq 1,
$$ 
defined as a Fourier multiplier:
$$
\F((-\Delta+\tau^2)^{-1}u) (\xi)= (|\xi|^2+\tau^2)^{-1}\hat{u}(\xi) , \quad u \in \calS(\R^n) .
$$
Note that, as opposed to differential operators, the operator $(-\Delta+\tau^2)^{-1}$ is non-local (in the sense that it does not satisfy $\supp(Pu) \subset \supp(u)$ for all $u \in C^\infty_c(\R^n)$). 
We write in this section a weak form of G\aa rding estimates for (almost-)differential operators of order $2$, which is at the core of the Carleman method.
\begin{proposition}[A local G\aa rding inequality for particular operators]
\label{lmweakGard}
Assume $\Omega$ is an open set with $0 \in \Omega$ and let $P$ be an operator of the form
\bnan
\label{form}
P=A+\sum_{i=1}^kB_i\circ (-\Delta+\tau^2)^{-1}\circ B_i
\enan
with $A,B_i\in \dif{2}(\Omega)$ with {\em real} principal symbols $a_2(x,\xi,\tau)$ and $b_{2,i}(x,\xi,\tau)$. Define 
\bnan
\label{e:def-p2}
p_2(x,\xi,\tau)=a_2(x,\xi,\tau)+\sum_{i=1}^k\frac{b_{2,i}^2(x,\xi,\tau)}{|\xi|^2+\tau^2}, 
\enan
 and assume that there is $C_0>0$ such that
\bnan
\label{pos}
p_2(0,\xi,\tau) \geq C_0(|\xi|^2+\tau^2) , \quad \text{ for all } \xi\in \R^n ,\tau\geq 0.
\enan
Then, there exist $r,C ,\tau_0>0$ such that 
\bna
\Re \left(P u,u\right)_{L^2}\geq  C \nor{u}{H^1_\tau}^2 , \quad \text{ for all } u\in C^{\infty}_c(B(0,r)) ,\tau\geq \tau_0.
\ena
\end{proposition}
Note that formally, such operators $P$ are ``of order $2$''. The ``principal symbol'', defined in~\eqref{e:def-p2} is indeed a homogeneous function of degree $2$. Inequality~\eqref{pos} is thus a homogeneous inequality, and it is sufficient to assume it on the half-sphere $\S^n_+ :=\{(\xi,\tau) \in \R^n \times \R_+ , |\xi|^2 + \tau^2 = 1\}$. We refer e.g. to~\cite{LL:book} for a proof relying on freezing coefficients, together with the fact that the result is elementary for the Fourier multiplier $p_2(0,D,\tau)$.

\section{Unique continuation under a convexity condition}
\label{chapterclassical}

This section is devoted to a classical unique continuation result under a pseudoconvexity condition, due to H\"ormander~\cite{Hoermander:63}. 

\subsection{Statement and examples}
\subsubsection{Statement of H\"ormander theorem}
Given $a,b \in C^\infty(\Omega \times \R^n; \C)$, the Poisson bracket $\left\{a,b\right\}$ is defined in~\eqref{d:def-poisson-bracket} and its geometric content is recalled in Section~\ref{rkconvexflow} below.
The geometric definition entering into the game is the following.
\begin{definition}[Strongly pseudoconvex hypersurface for operators of order two with real principal symbols]
\label{defconvexsurface}
Let $\Omega \ni x_0$ be an open set, $P\in \difclas{2}(\Omega)$ with real-valued principal symbol $p_2$ and $\Psi\in C^{\infty}(\Omega)$ real-valued. 
We say that the oriented hypersurface $S = \{\Psi = \Psi(x_0)\}$ is strongly pseudoconvex with respect to $P$ at $x_0$ if it satisfies
\bnan
\label{pseudoconvex2def}
p_2(x_0,\xi)=\{ p_2 ,\Psi\}(x_0, \xi)=0 \Longrightarrow \{p_2,\{p_2,\Psi\}\}(x_0, \xi)>0  \quad \text{for all } \xi \in  \R^n\setminus \{0\}  .
\enan
\end{definition}
Note that $\{ p , \Psi\}(x_0, \xi) =\d_\xi p(x_0,\xi)\cdot \d_x \Psi(x_0)$.
We can check that Definition \ref{defconvexsurface} is invariant if we change the defining function $\Psi$ (see e.g. Lemma~\ref{l:defining-function} and~\cite{LL:book}). That is why this is a geometric property of the oriented hypersurface solely. See Section~\ref{rkconvexflow} for an interpretation as convexity with respect to the bicharacteristic curves.
The geometric condition~\eqref{pseudoconvex2def} has to be compared with that discussed for vector fields in Examples~\ref{i:lin-vec-field-1}-\ref{i:lin-vec-field-2} in Section~\ref{s:intro-UC-problem}. The following result is a particular case of the H\"ormander theorem.
\begin{theorem}[Unique continuation for real operators of order $2$ under strong pseudoconvexity]
\label{thmreal2}
Let $\Omega \ni x_0$ be an open set of $\R^n$, and let $P\in \difclas{2}(\Omega)$ with {\em real} principal symbol $p_2$.
Assume that the oriented hypersurface $S = \{\Psi = \Psi(x_0)\}$ is strongly pseudoconvex with respect to $P$ at $x_0$. Then, there exists a neighborhood $V$ of $x_0$ so that for all $u\in H^1(\Omega)$, we have
\bnan
\label{e:UC-th-pseudo}
\left\{\begin{array}{l}
Pu = 0  \textnormal{ in  } \Omega,\\ u=0 \textnormal{ in } \Omega\cap \left\{\Psi> \Psi(x_0)\right\}\end{array}\right\}
\Longrightarrow u= 0 \textnormal{ in } V.
\enan
\end{theorem}
Another (slightly weaker) way to formulate the conclusion of the theorem is to say that $x_0 \notin \supp(u)$.
Here, we have assumed that all coefficients of $P$ are smooth for simplicity.  
\begin{remark}[Elliptic case]
\label{r:elliptic-rem}
Note now that in the particular case where the operator $P$ is elliptic at $x_0$, i.e. $p_2(x_0,\xi) \geq c|\xi|^2$, then the condition $p_2(x_0,\xi)=0$ is never fulfilled when $\xi \neq 0$ and~\eqref{pseudoconvex2def} is empty. Local unique continuation thus holds across {\em any} hypersurface. 
\end{remark}
Theorem~\ref{thmreal2} will be proved in Section \ref{sectUCP}. Before this, let us describe the underlying geometry of the condition~\eqref{pseudoconvex2def} together with typical examples of application of this result.

\subsubsection{Geometric content of the pseudoconvexity condition~\eqref{pseudoconvex2def}.}  
\label{rkconvexflow}

In this section, we explain the geometric content of the condition of Definition~\ref{defconvexsurface}. 
For this we need to introduce the Hamiltonian flow of the symbol $p_2$. Recall that $\{p_2 , \cdot\}$ (defined in~\eqref{d:def-poisson-bracket}) is a derivation on $C^\infty(\Omega \times \R^n)$  and can thus be identified with the vector field 
$$
H_{p_2}(x,\xi)= \d_\xi p_2(x,\xi) \cdot \d_x - \d_x p_2(x,\xi) \cdot \d_\xi  = \sum_{j=1}^n \d_{\xi_j} p_2(x,\xi) \frac{\d}{\d x_j}- \d_{x_j} p_2(x,\xi) \frac{\d}{\d \xi_j}
$$ 
on $\Omega \times \R^n$.  We denote by $\chi_s$ the associated flow, defined by 
\begin{equation}
\label{e:hamilton-1}
\left\{
\begin{array}{l}
\frac{d}{ds} \chi_s(x_0,\xi_0) = H_{p_2}\big(\chi_s(x_0,\xi_0) \big) , \\
\chi_0(x_0,\xi_0) =(x_0,\xi_0) ,
\end{array}
\right.
\end{equation}
and called the Hamiltonian flow of $p_2$.
Remark that $H_{p_2}(p_2) = \{p_2 , p_2\}=0$ so that $p_2$ is preserved along the flow: $p_2\circ\chi_s(x_0,\xi_0)= p_2(x_0,\xi_0)$.
Note also that the flow $\chi_s$ is (at least) locally defined for $(s,x,\xi)$ in a neighborhood of $(0,x_0,\xi_0)$ according to the Cauchy-Lipschitz theorem.

If we now denote by $(x_s, \xi_s)= \chi_s$, that is $\chi_s(x_0,\xi_0) = \big(x_s(x_0,\xi_0), \xi_s(x_0,\xi_0) \big)$ and recall the definition of the Poisson bracket $\{p_2 , \cdot\} = \d_\xi p_2 \cdot \d_x - \d_x p_2 \cdot \d_\xi$, \eqref{e:hamilton-1} now reads
\begin{equation}
\label{e:hamilton-2}
\left\{
\begin{array}{l}
\dsp \frac{d}{ds} x_s(x_0,\xi_0) = \d_\xi p_2 \big(\chi_s(x_0,\xi_0) \big) , \\
\dsp \frac{d}{ds} \xi_s(x_0,\xi_0) = - \d_x p_2 \big(\chi_s(x_0,\xi_0) \big) , \\
\dsp \big(x_s(x_0,\xi_0),\xi_s(x_0,\xi_0) \big) |_{s=0}=(x_0,\xi_0) .
\end{array}
\right.
\end{equation}

With these definitions in hand, we can now reformulate the strong pseudoconvexity condition of Definition~\ref{defconvexsurface}. 
Namely, note that we have
\begin{align*}
\{ p_2 ,\Psi\}(x_0, \xi) = H_{p_2} (\Psi)(x_0, \xi) = \frac{d}{ds} \Psi\circ x_s(x_0 ,\xi)|_{s=0}, \\
 \{ p_2, \{ p_2 ,\Psi\}\}(x_0, \xi) = H_{p_2} \big( H_{p_2} (\Psi)\big) (x_0, \xi) = \frac{d^2}{ds^2} \Psi\circ x_s(x_0 ,\xi)|_{s=0} .
 \end{align*}
Now, if for $\xi \in \R^n$ we define $c_\xi(s) =  \Psi\circ x_s(x_0 ,\xi)$, then \eqref{pseudoconvex2def} is equivalent to:
$$
\text{For all } \xi \in \R^n \setminus \{0\}, \text{ we have: } \quad p_2(x_0,\xi)=0 \text{ and } \dot{c}_\xi(0) = 0 \implies \ddot{c}_\xi(0) > 0 . 
$$
This means that for all $\xi \in \R^n \setminus \{0\}$, 
\begin{itemize}
\item if $\xi$ is noncharacteristic ($p_2(x_0,\xi)\neq 0$), we don't care;
\item if ($\xi$ is characteristic and) the (projected) Hamiltonian curve $x_s(x_0 ,\xi)$ is not tangent to $S=\{\Psi = \Psi(x_0)\}$ at $s=0$, we don't care;
\item if $\xi$ is characteristic and the curve $x_s(x_0 ,\xi)$ is tangent to $S=\{\Psi = \Psi(x_0)\}$, then it should have non-vanishing second derivative (tangency at order $2$) and the curve $(x_s(x_0 ,\xi))_{s\in (-\eps,\eps)}$ should stay in $\{\Psi \geq \Psi(x_0)\}$ (for $\eps>0$ small enough).
\end{itemize}
This excludes the situations where tangent characteristic curves stay in $\{\Psi \leq \Psi(x_0)\}$ or the case of contacts of higher order with the tangent space at $x_0$.
\begin{remark}
Note also that, as a consequence of the above discussion, the fact that $\{\Psi=0\}$ satisfies the strong pseudoconvexity condition at $x_0$ implies in particular that the set $\{\Psi>0\}$ ``controls geometrically'' a whole neighborhood of $x_0$, in the sense of~\cite{RT:74,BLR:92}: any null bicharacteristic curve of $P$ passing close enough to $x_0$ intersects $\{\Psi>0\}$. 
In particular if $P$ has smooth coefficients in $\Omega$ and $u \in \mathcal{D}'(\Omega)$ satisfies $Pu=0$ in $\mathcal{D}'(\Omega)$ and $u=0$ in $\{\Psi>0\}$, then propagation of singularities~\cite[Theorem~23.2.9]{Hoermander:V3} together with the fact that $\{\Psi>0\}$ satisfies this local form of geometric control condition imply that $u \in C^\infty$ in a neighborhood of $x_0$, and Theorem~\ref{thmreal2} applies. 
\end{remark}

\subsubsection{Examples}
\label{s:examples-hor}
\begin{remark}[Operators with constant coefficients]
\label{r:cte-coeff}
Consider here the simple case where $P = AD \cdot D$ where $A$ is a constant real symmetric matrix (think $A =\diag(-1, 1, \cdots,1)$ for the flat/Minkowski wave operator). This is $p_2(x,\xi) = A\xi\cdot \xi$.
We have $\{p_2, \Psi\} (x,\xi)= 2A \xi\cdot d\Psi(x)$ and $\{p_2, \{p_2, \Psi\} \} (x,\xi) = \{A\xi\cdot \xi , 2A \xi\cdot d\Psi(x)\} = 4\Hess\Psi(x) (A\xi , A\xi)$. Condition~\eqref{pseudoconvex2def} rewrites 
$$
A \xi\cdot \xi =0 \quad \text{and}\quad A \xi\cdot d\Psi(x_0)=0 \Longrightarrow \Hess\Psi(x_0) (A\xi , A\xi) >0 \quad \text{for all } \xi \in  \R^n\setminus \{0\}  .
$$
Notice that, in case $A$ is invertible, this cotangent formulation might be equivalently replaced by the following tangent one, setting $X=A\xi$:
\begin{equation}
\label{e:tangent-formulation-pseudoconv}
A^{-1} X\cdot X =0 \quad \text{and}\quad d\Psi(x_0)(X)=0 \Longrightarrow \Hess\Psi(x_0) (X , X) >0 \quad \text{for all } X \in  \R^n\setminus \{0\}  .
\end{equation}
\end{remark}

\begin{remark}[Pseudo-Riemannian metric and operator]
\label{r:waves-tangent}
As in Remark~\ref{r:cte-coeff}, in the general case of operators $P=D_ia^{ij}(x)D_j$ (plus lower order terms) with variable coefficients but nondegenerate cometric $a^{ij}(x_0)$, one may want to rephrase the cotangent formulation of strong pseudoconvexity~\eqref{pseudoconvex2def} on the tangent space at $x_0$. Note that this covers the wave operator in case the signature of $a^{ij}(x_0)$ is $(-1,1,\cdots,1)$ (i.e. the metric is Lorentzian).
Introducing the pre-dual metric $a(x) = a_{ij}(x) = (a^{ij}(x))^{-1}$ and setting $X=a(x)^{-1}\xi$, i.e. $X^i=a^{ij} \xi_j$ (where the Einstein summation convention is used), one can check that 
$\{p_2,\Psi\}(x_0,\xi) = 2 d\Psi(x_0)(X) = (\nabla_a \Psi ,X)_a$, where $\nabla_a$ is the usual pseudo-Riemannian gradient and $ (Y ,X)_a = a_{ij}(x)Y^iX^j$. Moreover, one can also check that 
$\{p_2,\{p_2,\Psi\}\}(x_0,\xi) = 4 \Hess_a\Psi(x_0)(X,X)$, where $\Hess_a\Psi(X,Y) = (D_Xd\Psi)(Y)$ is the usual pseudo-Riemannian Hessian of $\Psi$ (and $D$ denotes the Levi-Civita connection associated to the metric $a$, see e.g.~\cite[Chapter~II Section~B]{GallotHulinLaf}). Condition~\eqref{pseudoconvex2def} rewrites as 
$$
(X, X)_a =0 \quad \text{and}\quad d\Psi(x_0)(X)=0 \Longrightarrow \Hess_a \Psi(x_0) (X , X) >0 \quad \text{for all } X \in  T_{x_0}\Omega \setminus \{0\}  .
$$
\end{remark}

\begin{remark}[The flat/Minkowski wave operator]
\label{rkwave}
We discuss here the case of the wave operator with constant coefficients, which is a particular case of the above examples with $A = \diag(-1, 1, \cdots , 1)$, $P=\partial_t^2-\Delta$, $p_2=-\xi_t^2+|\xi_x|^2$.
We may now write the strong pseudoconvexity condition~\eqref{pseudoconvex2def} specialized in the point $(t,x)=(0,0)$ (the operator is translation invariant in $(t,x)$), in different situations. 
\begin{itemize}
\item Spacelike hypersurface: if $|\partial_t \Psi(0)| > |\nabla_x\Psi(0)|$: the hypersurface $\{\Psi = \Psi(0)\}$ is called spacelike (its normal vector $\nabla_{t,x} \Psi$ is timelike). The first two conditions imply  $|\xi_t \d_t\Psi(0)|=|\xi_x \cdot \d_x \Psi(0)| \leq |\xi_x| |\d_x\Psi(0)| =|\xi_t| |\d_x \Psi(0)| <|\xi_t|  |\partial_t \Psi(0)|$. This is a contradiction, and hence Condition~\eqref{pseudoconvex2def} is empty.

Any spacelike hypersurface satisfies the unique continuation property. This is natural since the Cauchy problem is hyperbolic and thus locally wellposed for any spacelike hypersurface (like for instance the wave equation posed with initial data at $t=0$, see Theorem~\ref{t:wave-finite-speed}). 
\item Time-invariant hypersurface: in the applications discussed in Section~\ref{s:motivation}, a typical unique continuation result needed is rather across hypersurfaces of the form $\Psi(t,x)=\varphi(x)$. The strong pseudoconvexity condition then writes
\bna
\xi_t^2=|\xi_x|^2 \textnormal{ and }\xi_x \cdot \nabla_x \varphi(0)=0 \Longrightarrow \Hess_x (\varphi)(0)(\xi_x,\xi_x)>0 \quad \forall \xi\in \R^n\setminus \{0\}.
\ena
Typically, if $\Psi(t,x)=|x|^2-1$, this condition holds (and Theorem~\ref{thmreal2} implies unique continuation) from the exterior of the cylinder $\{(t,x), \Psi(t,x)>0\}$ towards the interior but not in the other direction. 
Note also that for the $1D$ wave equation, the constraint $\xi_x \cdot \nabla_x \varphi(0)=0$ is much more demanding and implies $\xi_x=0$ and $\xi=0$ if $\xi_t^2=|\xi_x|^2$. This is natural since time and space variables play the same role. Hence, finite speed of propagation essentially implies the local unique continuation property across any non characteristic hypersurface.

\item More generally, if one considers the function $\Psi(t,x) = |x-x_1|^2-\gamma^2t^2$, then we have, for any $(t_0,x_0)\in \R^{1+d}$
$\Hess(\Psi)(t_0,x_0)(X,X)= 4( |X_x|^2-\gamma^2|X_t|^2)$. As a consequence, on the tangent cone $|X_x| = |X_t| , X \neq 0$, we have $\Hess(\Psi)(t_0,x_0)(X,X)= 4 |X_x|^2 (1-\gamma^2)$ with $X_x\neq 0$. That is to say that the strong pseudoconvexity condition~\eqref{e:tangent-formulation-pseudoconv} is satisfied at any point $(t_0,x_0)$ as soon as $\gamma \in [0,1)$.
Then, Theorem~\ref{thmreal2} implies that the local unique continuation property holds at any point $(t_0,x_0)$ of the one sheet hyperboloid of revolution $\{\Psi = \Psi(t_0,x_0)\}$ from the exterior of the hyperboloid $\{\Psi > \Psi(t_0,x_0)\}$ towards the interior $\{\Psi < \Psi(t_0,x_0)\}$.

Using a compactness argument and letting $\gamma$ close to one, one may deduce the following global result.
\begin{proposition}
\label{p:global-pseudoconvex}
Assume $\M$ is the closure of a bounded open set in $\R^d$, let $\omega$ be an open neighborhood of $\d \M$ in $\M$, and fix any $x_1 \in \R^d$. 
Then, assuming 
\begin{align*}
T>\sup\{|x-x_1|, x \in \M\setminus \omega\} ,
\end{align*}
and letting $q_0,q_1 \in L^\infty((-T,T)\times \M ; \C), q_2 \in L^\infty((-T,T)\times \M ; \C^d)$, any solution to
\begin{equation*}
\left\{\begin{array}{l}
(\d_t^2-\Delta)u +q_0 u + q_1\d_tu + q_2\cdot \nabla_x u= 0  \textnormal{ in  } (-T,T)\times \Int(\M),\\ 
u \in H^1((-T,T)\times \M) , \\
u =0 \textnormal{ in } (-T,T)\times \omega  ,
\end{array}\right.  
\end{equation*}
satisfies 
$$
 (u,\d_t u)|_{t=0}= (0,0) \textnormal{ in }  \Int(\M) .
$$
\end{proposition}
Notice that the unique continuation result first proves that $$u=0 \text{ a.e. in }(-T,T)\times \M \cap \{(t,x) \in \R^{1+d}, |x-x_1|^2-\gamma^2t^2>0\}$$ (for all $\gamma <1$, and thus for $\gamma =1$ as well), and in particular in a neighborhood of any point in $\{0\} \times( \M\setminus \{x_1\})$. Then, regularity of the wave equation implies that $u \in C^0((-T,T);H^1_{\loc}(\M)) \cap C^1((-T,T);L^2_{\loc}(\M))$, so that $(u,\d_t u)|_{t=0} \in H^1_{\loc}(\M)\times L^2_{\loc}(\M)$ have $\supp((u,\d_t u)|_{t=0} ) \subset \{x_1\}$ and thus vanish identically.
\end{itemize}
\end{remark}

\subsection{The Carleman estimate}
\label{s:Carleman}
Here, we recall that, $P\in \difclas{2}(\Omega)$ (with principal symbol $p_2(x,\xi)$) and $\Phi\in C^{\infty}(\Omega; \R)$ being given, the conjugated operator is $P_\Phi = e^{\tau \Phi }P e^{-\tau \Phi } \in \dif{2}(\Omega)$ and its principal symbol is $p_\Phi(x,\xi,\tau) = p_2(x, \xi+i d\Phi(x))$ (computed in Lemma \ref{lmsympphi} and Example~\ref{e:general-operator-2}). We write $p_\Phi = \Re(p_\Phi) + i \Im (p_\Phi)$.

\subsubsection{Carleman estimate under subellipticity condition}
As we have seen in Section \ref{subsectstratCarle}, the key intermediate step for proving Theorem~\ref{thmreal2} is to obtain estimates of the type of~\eqref{Carlintro}. 
In this section, we prove a Carleman estimate (Theorem \ref{thmCarlclass}) under a symbolic condition usually called ``H\"ormander subellipticity condition'' (namely~\eqref{ineginverse}). 
The next sections link this condition to the strong pseudoconvexity condition~\ref{pseudoconvex2def}.

\begin{theorem}[Local Carleman estimate]
\label{thmCarlclass}
Let $\Omega$ be an open subset of $\R^n$ and $x_0\in \Omega$. Let $P\in \difclas{2}(\Omega)$ be a (classical) differential operator with {\em real-valued principal symbol} $p_2$ and $\Phi\in C^{\infty}(\ovl{\Omega}; \R)$.
Assume there exist $C_1, C_2 >0$ such that for all $(\xi , \tau) \in \R^n \times \R_+^*$, 
 \bnan
\label{ineginverse}
 \frac{C_1}{|\xi|^2+\tau^2}\left[(\Re p_{\Phi})^2+(\Im p_{\Phi})^2\right]+\frac{1}{\tau} \left\{\Re p_{\Phi},\Im p_{\Phi}\right\}\geq C_2 \left(|\xi|^2+\tau^2\right),
 \enan
 where the symbols are taken at the point $(x_0,\xi,\tau)$.
Then, there  exist $C, r, \tau_0>0$  such that
\begin{align}
\label{Carlthm}
\tau \nor{v}{\Ht{1}}^2\leq C \nor{P_\Phi v}{L^2}^2,& \quad \text{for all }v \in C^{\infty}_c (B(x_0, r)) , \tau\geq \tau_0 ; \\
\label{Carlthm-bis}
 \tau^3 \nor{e^{\tau \Phi} u}{L^2}^2 +\tau \nor{e^{\tau \Phi }\nabla u}{L^2}^2 \leq C \nor{e^{\tau \Phi }P u}{L^2}^2, &\quad \text{for all }u \in C^{\infty}_c (B(x_0, r)) , \tau\geq \tau_0 .
\end{align}
\end{theorem}
Notice that $\frac{\Im p_{\Phi}}{\tau} = 2 \widetilde{p_2}(x, \xi, d\Phi(x))$ (see Example~\ref{e:general-operator-2}) is smooth, so this is not a problem to divide by $\tau$ in~\eqref{ineginverse}, even when $\tau\to 0^+$.
Before proceeding to the proof of this result, several comments are in order. 
First, the statement~\eqref{Carlthm-bis} is useful for applications unique continuation, see~\eqref{Carlintro} and the discussion in Section~\ref{subsectstratCarle}.
The statement~\eqref{Carlthm} is only a reformulation in terms of the conjugated operator, which belongs to $\dif{2}$, and is thus analyzable with the tools developed in Section~\ref{s:operators-dif-tau}. The statement~\eqref{ineginverse},   as opposed to~\eqref{Carlthm}--\eqref{Carlthm-bis}, is a ``symbolic estimate'', concerning only the principal symbol of the conjugated operator.
The interest of this result is that it reduces the problem of proving a Carleman estimate to a checkable property on the principal symbol of the conjugated operator.
The question of rephrasing the condition~\eqref{ineginverse} in geometric terms is addressed in Sections~\ref{s:pseudo-functions} and~\ref{s:strong-pseudo-surfaces}  below.
The useful information in this theorem is \eqref{ineginverse}$\implies$\eqref{Carlthm-bis}. The converse implication is also true, which indicates the limit of this classical Carleman approach.
This point is slightly more technical and we refer the reader to~\cite[Section~28.2]{Hoermander:V4} for a proof.

\bnp
The equivalence between~\eqref{Carlthm-bis} and~\eqref{Carlthm} comes from the change of unknown $v=e^{\tau\Phi}u$.  This yields $P_{\Phi}v = e^{\tau \Phi}Pe^{-\tau \Phi } v= e^{\tau \Phi }P u$. Moreover, we have 
$\nabla u = \nabla (e^{-\tau\Phi}v) = e^{-\tau\Phi} (\nabla v - \tau v \nabla \Phi)$ so that 
$$
 \tau^2 \nor{e^{\tau \Phi} u}{L^2}^2 + \nor{e^{\tau \Phi }\nabla u}{L^2}^2 
 \leq  \tau^2 \nor{v}{L^2}^2 + 2\nor{\nabla v}{L^2}^2 + 2 \nor{\tau v \nabla \Phi}{L^2}^2  \leq C\nor{v}{\Ht{1}}^2 ,
$$
and thus~\eqref{Carlthm} implies~\eqref{Carlthm-bis}. Conversely, we have $\nabla v = \nabla (e^{\tau\Phi}u) = e^{\tau\Phi} (\nabla u + \tau u \nabla \Phi)$ so that 
\begin{align*}
\nor{v}{\Ht{1}}^2 & = \nor{\nabla v}{L^2}^2+\tau^2 \nor{v}{L^2}^2 \leq 2 \nor{e^{\tau\Phi} \nabla u }{L^2}^2 + 2 \nor{e^{\tau\Phi} \tau u \nabla \Phi}{L^2}^2+\tau^2 \nor{e^{\tau \Phi }u}{L^2}^2 \\ 
&  \leq C  \left( \tau^2 \nor{e^{\tau \Phi} u}{L^2}^2 + \nor{e^{\tau \Phi }\nabla u}{L^2}^2 \right) ,
\end{align*}
and~\eqref{Carlthm-bis} implies~\eqref{Carlthm}.

 \medskip
We now want to prove that \eqref{ineginverse} implies \eqref{Carlthm}.
Before going further, let us notice that Lemma \ref{lmsympphi} only depends on the leading order of the operator $P$. More precisely, if $\tilde{P} \in \difclas{2}(\Omega)$ has the same principal symbol as $P$, then $P-\tilde{P}\in \difclas{1}(\Omega)$ and $P_\Phi-\tilde{P}_\Phi =R \in \dif{1}(\Omega)$. Henceforth, assuming the Carleman inequality~\eqref{Carlthm} for $\tilde{P}_\Phi$, 
\begin{equation}
\label{e:carleman-tilde}
\tau \nor{v}{\Ht{1}}^2\leq C \nor{\tilde{P}_\Phi v}{L^2}^2
\end{equation}
yields 
$$
\tau \nor{v}{\Ht{1}}^2\leq C \nor{(P_\Phi -R) v}{L^2}^2 \leq C \nor{P_\Phi v}{L^2}^2+  C \nor{R v}{L^2}^2 \leq C \nor{P_\Phi v}{L^2}^2+  D \nor{v}{H^1_\tau}^2
$$
with Item~\ref{i:propSob} in Proposition~\ref{p:calcul}. 
Then, for $\tau$ large enough, we have $\tau -D\geq \tau/2$, and the last term can then be absorbed in the left hand-side, yielding the sought Carleman inequality~\eqref{Carlthm} for $P_\Phi$, with different constants $C$ and $\tau_0$.

Since the operator $P$ has real principal symbol $p_2$, we choose $\tilde{P} = \frac{P+P^*}{2}$, which is selfadjoint and has the same principal symbol $p_2$. Note that we may write 
$P= \sum_{i,j=1}^n a^{ij}(x) D_i D_j$ with $a^{ij} = a^{ji}$ real-valued (see Example~\ref{e:symmetric-aij}), and we have $\tilde{P} = \sum_{i,j=1}^n D_i a^{ij}(x) D_j$ modulo a (selfadjoint) first order operator.
We may thus focus on $\tilde{P}_\Phi  =e^{\tau \Phi} \tilde{P} e^{-\tau \Phi}$ and prove~\eqref{e:carleman-tilde}.
To this aim, we decompose the operator $\tilde{P}_\Phi$ as
\bnan
\label{decomp}
\tilde{P}_\Phi=Q_R+iQ_I , \quad \text{with} \quad Q_R=\frac{\tilde{P}_\Phi+\tilde{P}_\Phi^*}{2};\quad Q_I=\frac{\tilde{P}_\Phi-\tilde{P}_\Phi^*}{2i}
\enan
Note that both $Q_R$ and $Q_I$ are formally selfadjoint ($Q_R^*=Q_R$ and $Q_I^*=Q_I$), and, according to Item~\ref{i:propadjoint} in Proposition~\ref{p:calcul}, we have $Q_R, Q_I \in \dif{2}$ with principal symbols (see Example~\ref{e:general-operator-2})
\begin{align*}
&q_{R}(x,\xi,\tau)=\frac{p_{\Phi}+\overline{p_{\Phi}}}{2}(x,\xi,\tau)=\Re p_{\Phi}(x,\xi,\tau)  =  p_2(x, \xi ) - \tau^2 p_2(x, d\Phi(x)) ,\\
&q_{I}(x,\xi,\tau)=\frac{p_{\Phi}-\overline{p_{\Phi}}}{2i}(x,\xi,\tau)=\Im p_{\Phi} (x,\xi,\tau)=2 \tau \widetilde{p_2}(x, \xi, d\Phi(x)) .
\end{align*}
Moreover (this is a key point), $\tilde{P}_\Phi = \sum_{i,j=1}^n (D_i +i \tau \d_i\Phi)  a^{ij}(x) (D_j+i \tau \d_j\Phi)$ and we may hence write $\tilde{P}_\Phi  = \tilde{P} + \tau M$ for some $M\in \dif{1}(\Omega)$. Since $\tilde{P}$ is selfadjoint, this implies 
\bnan
\label{e:factor-qi}
Q_I=\frac{\tilde{P}_\Phi-\tilde{P}_\Phi^*}{2i} = \frac{\tau M-\tau M^*}{2i}  = \tau \tilde{Q}_I , \quad \text{ with } \quad  \tilde{Q}_I  = \frac{M-M^*}{2i}  \in \dif{1}(\Omega) , 
\enan
i.e. $\tau$ may be factorized in the skewadjoint part of the operator.

Using~\eqref{decomp}, the central computation is now as follows, for $v\in C^{\infty}_c(\Omega)$, 
\begin{align}
\label{e:e:carleman-key}
\nor{\tilde{P}_\Phi v}{L^2}^2&=\left(\tilde{P}_\Phi v, \tilde{P}_\Phi v\right)_{L^2}=\left((Q_R+iQ_I) v,(Q_R+iQ_I)v\right) \nonumber\\
&=\left(Q_R v,Q_Rv\right)+\left(iQ_I v,iQ_Iv\right)+\left(Q_Rv,iQ_Iv\right)+\left(iQ_I v,Q_Rv\right)\nonumber\\
&=\nor{Q_R v}{L^2}^2+\nor{Q_I v}{L^2}^2+\left(i[ Q_R,Q_I]v,v\right). 
\end{align}
Now, we have $2$ kinds of terms
\begin{itemize}
\item  the one with $\nor{Q_R v}{L^2}^2$ (and resp. $\nor{Q_I v}{L^2}^2$) that corresponds to $\left(Q_R^2 v,v\right)$ where $Q_R^2$ is of order $4$ with principal symbol $(\Re p_{\Phi})^2$ (resp. $(\Im p_{\Phi})^2$);
\item 
the one with $i[Q_R,Q_I]$ which is of order $2+2-1=3$ and principal symbol $\{\Re p_{\Phi},\Im p_{\Phi}\}$ by Item~\ref{i:commutator} of Proposition~\ref{p:calcul}.
\end{itemize}
The first two operators have stronger order ($4$) but they can cancel and are therefore not sufficient to obtain the ``coercivity'' estimate. The idea is thus to use the commutator where both $q_{R}$ and $q_{I}$ cancel. However, to compare these terms, we need to bring them to the same order and ``sacrifice'' this main order $4$. More precisely, let $C_1>0$ be as in Assumption~\eqref{ineginverse} (that this is the right constant will appear in~\eqref{inegperte} below). For $\tau \geq C_1$, we have 
$$\frac{1}{\tau^{1/2}}\geq \frac{C_1^{1/2}}{\tau} \geq \frac{C_1^{1/2}}{(|\xi|^2+\tau^2)^{1/2}} \quad \text{for all } \xi\in \R^n .
$$
This implies (using again the Plancherel Theorem)
\begin{align}
\label{e:carleman-sacrifice}
\frac{1}{\tau}\nor{Q_R v}{L^2}^2=\nor{\frac{Q_R v}{\tau^{1/2}}}{L^2}^2&\geq \nor{C_1^{1/2}(-\Delta+\tau^2)^{-1/2}Q_R v}{L^2}^2 =  C_1 \left(Q_R (-\Delta+\tau^2)^{-1}Q_Rv,v\right).
\end{align}
The same estimate applies to $Q_I$.
Combining~\eqref{e:e:carleman-key} with~\eqref{e:carleman-sacrifice}, we have now proved
\bnan
\label{inegperte}
\frac{1}{\tau}\nor{\tilde{P}_{\Phi}v}{L^2}^2\geq \left(Lv,v\right)_{L^2}, 
\enan
with 
\begin{equation*}
L=C_1 \Big( Q_R (-\Delta+\tau^2)^{-1} Q_R+Q_I (-\Delta+\tau^2)^{-1} Q_I \Big)+i \left[Q_R,\frac{Q_I}{\tau}\right] .
\end{equation*}
But we have proved in~\eqref{e:factor-qi} that $Q_I = \tau \tilde{Q}_I$ with $\tilde{Q}_I \in \dif{1}(\Omega)$.
This implies that $\left[Q_R,\frac{Q_I}{\tau}\right] \in \dif{2}$ as well. The operator $L$ is thus precisely of the form of that in Proposition~\ref{lmweakGard}, is moreover selfadjoint, and has principal symbol (in the sense of Proposition~\ref{lmweakGard})
$$
 \frac{C_1}{|\xi|^2+\tau^2}\left( (\Re p_{\Phi})^2+(\Im p_{\Phi})^2\right)+ \left\{\Re p_{\Phi}, \frac{\Im p_{\Phi}}{\tau}\right\}, 
$$
which satisfies~\eqref{ineginverse}. Hence, the G{\aa}rding inequality of Proposition~\ref{lmweakGard} applies and yields the existence of $C, \tau_0 ,r>0$ such that 
$$
 \left(Lv,v\right)_{L^2} \geq C \nor{v}{H^1_\tau}^2 , \quad \text{ for all }v \in C^\infty_c(B(x_0, r)), \quad \tau \geq \tau_0 , 
$$
which, in view of~\eqref{inegperte}, yields~\eqref{e:carleman-tilde} and concludes the proof of the Carleman estimate~\eqref{Carlthm}.
\enp
Note that in~\eqref{e:carleman-sacrifice}, since $Q_R$ is only defined on $\Omega$, and since $(-\Delta+\tau^2)^{-1}Q_Rv \notin C^\infty_c(\Omega)$, the expression $Q_R (-\Delta+\tau^2)^{-1}Q_Rv$ is not well-defined. However, its pairing with the function $v \in C^\infty_c(B(x_0, r))$ is well defined (e.g. as $\left( \chi Q_R\tilde{\chi} (-\Delta+\tau^2)^{-1}Q_Rv,v\right)$ with $\chi \in C^\infty(\Omega)$ such that $\chi=1$ on a neighborhood of $B(x_0,r)$, and $\tilde{\chi} \in C^\infty(\Omega)$ with $\tilde{\chi}=1$ on a neighborhood of $\supp(\chi)$).   

\begin{remark}[Lower order terms]
\label{r:lot}
As seen in the proof, an important feature of the Carleman estimates~\eqref{Carlthm-bis} is its insensitivity with respect to lower order terms. More precisely, if~\eqref{Carlthm-bis} is satisfied for an operator $P$, then it also holds for $P' := P+ \sum_{k=1}^n b_k(x)D_k + c(x)$ as soon as $b_k, c \in L^\infty(\Omega)$. Indeed, applying~\eqref{Carlthm-bis}  for $P$ yields 
\begin{align*}
 \tau^3 \nor{e^{\tau \Phi} u}{L^2}^2 +\tau \nor{e^{\tau \Phi }\nabla u}{L^2}^2 
 & \leq C \Big\|e^{\tau \Phi } \Big(P' -\sum_{k=1}^n b_k(x)D_k + c(x) \Big)u \Big\|_{L^2}^2 \\
 & \leq C \nor{e^{\tau \Phi }P'}{L^2}^2 +C \nor{e^{\tau \Phi }\nabla u}{L^2}^2 +C \nor{e^{\tau \Phi }u}{L^2}^2,
\end{align*}
and the last two terms can be absorbed in the left handside for $\tau$ large enough. Note in particular that no regularity is required on the lower order terms when proceeding that way.
\end{remark}
 
\subsubsection{Carleman estimate for pseudoconvex functions}
\label{s:pseudo-functions}
We now reduce the quantitative symbolic Assumption~\eqref{ineginverse} of the Carleman estimate to a qualitative convexity condition on the weight function $\Phi$ (with respect to the symbol $p_2$).
\begin{definition}[Pseudoconvexity for functions]
\label{defconvexfunctions}
Let $\Omega \ni x_0$ be an open set, $P\in \difclas{2}(\Omega)$ be a (classical) differential operator with real-valued principal symbol $p_2$ and $\Phi\in C^{\infty}(\Omega)$ real-valued. 

We say that the function $\Phi$ is pseudoconvex with respect to $P$ at $x_0$ if it satisfies
\bnan
\label{e:pseudo-function-1}
\left\{ p_2, \{p_2 ,\Phi \} \right\} (x_0 ,\xi) >0  , 
&\text{ if } p_2(x_0,\xi) = 0 \text{ and } \xi \neq 0 ; \\
\label{e:pseudo-function-2}
\frac{1}{i\tau}\{ \overline{p_\Phi},p_\Phi \} (x_0 ,\xi, \tau) >0  , 
&\text{ if } p_\Phi(x_0,\xi,\tau) = 0 \text{ and } \tau >0 ,
\enan
where $p_\Phi(x,\xi,\tau) = p_2(x, \xi + i \tau d\Phi(x))$.
\end{definition}
 Lemma \ref{lmlien} shows that~\eqref{e:pseudo-function-1} is the limit of~\eqref{e:pseudo-function-2} as $\tau \to 0^+$.
\begin{lemma}
\label{lmlien}
Let $p$ be a real-valued smooth function on $\Omega \times \R^n$. Then, we have $\underset{\tau \rightarrow 0^+}{\lim}\frac{1}{i\tau}\{ \overline{p_\Phi},p_\Phi \}(x,\xi,\tau)=2\left\{p\left\{p,\Phi\right\}\right\}(x,\xi)$ for all $(x,\xi)\in\Omega \times \R^n$.
\end{lemma}

We now state the equivalence between Definition~\ref{defconvexfunctions} and the H\"ormander subellipticity condition~\eqref{ineginverse}.
\begin{proposition}
\label{p:pseudo-hormander}
Let $\Omega \ni x_0$ be an open set, $P\in \difclas{2}(\Omega)$ with real-valued principal symbol  $p_2$ and $\Phi\in C^{\infty}$ real-valued. If $\Phi$ is pseudoconvex with respect to $P$ at $x_0$, then the subellipticity condition~\eqref{ineginverse} is satisfied at $x_0$.
\end{proposition}
And hence, if $\Phi$ is a pseudoconvex function in the sense of Definition~\ref{defconvexfunctions}, the Carleman estimate of Theorem~\ref{thmCarlclass} holds with weight $\Phi$.
The proof uses the following (elementary but very useful) lemma.
\begin{lemma}
\label{l:fgh-three-fcts}
Let $K$ be a compact set and $f,g$ two continuous real-valued functions on $K$. Assume that $f \geq 0$ on $K$, and $g>0$ on $\{f=0\}$. Then, there exists $A_0, C>0$ such that for all $A\geq A_0$, we have $g + A f \geq C$ on $K$. 
\end{lemma}
A proof of this elementary lemma can be found e.g. in ~\cite{LL:book}.
We now prove Proposition~\ref{p:pseudo-hormander} from Lemmata~\ref{lmlien} and~\ref{l:fgh-three-fcts}, and finally come back to the proof of Lemma~\ref{lmlien}.

\bnp[Proof of Proposition~\ref{p:pseudo-hormander}]
Note first that since $\{f,f\}=0$ and $\{f,g\}=-\{g,f\}$ for any $f$ and $g$, we have
\begin{align*}
\frac{1}{i\tau}\{ \overline{p_\Phi},p_\Phi \}&=\frac{1}{i\tau} \{\Re p_{\Phi}-i\Im p_{\Phi}, \Re p_{\Phi}+i\Im p_{\Phi} \}\\
&=\frac{1}{\tau}\{\Re p_{\Phi},\Im p_{\Phi} \}-\frac{1}{\tau}\{\Im p_{\Phi},\Re p_{\Phi} \} = \frac{2}{\tau}\{\Re p_{\Phi},\Im p_{\Phi} \}.
\end{align*}
Moreover, we recall that $\frac{\Im p_{\Phi}}{\tau} = 2 \tilde{p}_2(x_0, \xi, d\Phi(x_0))$ is smooth.
We notice that all terms in \eqref{ineginverse} are homogeneous in $(\xi,\tau)$ of order $2$ and continuous thanks to the previous remark. Therefore, it is enough to prove \eqref{ineginverse} on the set $K=\left\{(\xi,\tau), |\xi|^2+\tau^2=1;\tau\geq 0 \right\} $. On this compact set, the result is a consequence of Lemma~\ref{l:fgh-three-fcts} with $f=(\Re p_{\Phi})^2+(\Im p_{\Phi})^2$ and  $g=2\{\Re p_{\Phi},\frac{\Im p_{\Phi}}{\tau}\}$.
Lemma \ref{lmlien} then proves that  the first assumption in Definition \ref{defconvexfunctions} is the limit of the second one on the set $\{\tau=0\}$. Hence, we have $g>0$ on $\{f=0\}$  on the whole $K$, up to the set $\{\tau=0\} \cap \{|\xi|^2+\tau^2=1\}$.
 Lemma~\ref{l:fgh-three-fcts} then concludes the proof of the subellipticity condition~\eqref{ineginverse}.
\enp

\bnp[Proof of Lemma \ref{lmlien}]
We first notice that for $\tau=0$, $\{ \overline{p_\Phi},p_\Phi \}=\{\overline{p},p\}$ so  since $p$ is real, $\{ \overline{p_\Phi},p_\Phi \}=0$ for $\tau=0$. The definition of the derivative in $\tau=0$ then yields 
\bnan
\label{e:lim-tau0}
\underset{\tau \rightarrow 0}{\lim}\frac{1}{\tau}\{ \overline{p_\Phi},p_\Phi \}=\left.\frac{\partial}{\partial \tau} \{ \overline{p_\Phi},p_\Phi \}\right|_{\tau=0} .
\enan
Also, we have $\partial_{\tau}(\{ \overline{p_\Phi},p_\Phi \})=\{ \partial_{\tau}\overline{p_\Phi},p_\Phi \}+\{ \overline{p_\Phi},\partial_{\tau}p_\Phi \}$.
But since $p$ is real, $\overline{p_\Phi}=p(x,\xi-i\tau d \Phi(x))$, so that
\begin{align*}
\partial_{\tau}p_\Phi (x, \xi, \tau) &= i d \Phi \cdot \d_{\xi}p(x,\xi+i\tau d\Phi)=i\{p_{\Phi},\Phi\} (x, \xi, \tau),\\
\partial_{\tau}\overline{p_\Phi}(x, \xi, \tau)&=-i d \Phi \cdot \d_{\xi}p(x,\xi-i\tau d\Phi)=-i\{\overline{p_{\Phi}},\Phi\}(x, \xi, \tau).
\end{align*}
We obtain $\partial_{\tau}(\{ \overline{p_\Phi},p_\Phi \})=-i\{ \{\overline{p_{\Phi}},\Phi\},p_\Phi \}+i\{ \overline{p_\Phi},\{p_{\Phi},\Phi\} \}$, which,  specified to $\tau=0$, writes
\bna
\left.\frac{\partial}{\partial \tau} \{ \overline{p_\Phi},p_\Phi \}\right|_{\tau=0}=-i\{ \{p,\Phi\},p \}+i\{ p,\{p,\Phi\} \}=2i\{ p,\{p,\Phi\} \}.
\ena
Together with~\eqref{e:lim-tau0}, this concludes the proof of the lemma.
\enp
  A very important drawback to Definition \ref{defconvexfunctions} is that, it is not only dependent on the level set of the functions, but also on the ``convexity with respect to the level sets''. This is not a geometric assumption (in general, $g''(x_0)$ is a geometric quantity only if $g'(x_0)=0$). We now need to link this definition to geometric quantities, so that to be able to formulate a result with, at least, a geometric assumption (that is invariant by diffeomorphisms).
  Before that, let us stress an important stability feature of the pseudoconvexity assumption of Definition \ref{defconvexfunctions}.

\subsubsection{Stability of the pseudoconvexity assumption}
We prove that the pseudoconvexity condition of Definition~\ref{defconvexfunctions} is stable by small $C^2$ perturbations of the weight function $\Phi$. This will be very useful for perturbing the hypersurface across which to prove unique continuation.

\begin{proposition}[Stability and Geometric convexification]
\label{propgeometricconvex}
Let $\Omega \ni x_0$ such that $\ovl{\Omega}$ is compact. Assume $P\in \difclas{2}(\Omega)$ has real-valued principal symbol, and $\Phi\in C^{\infty}$ is pseudoconvex with respect to $P$ at $x_0$ (in the sense of Definition \ref{defconvexfunctions}).
Then there exists $\e_0>0$ so that any $\Phi_{\e} \in C^2(\ovl{\Omega})$ with $\|\Phi - \Phi_\e \|_{C^2(\ovl{\Omega})} <\e_0$ is pseudoconvex with respect to $P$ at $x_0$.
\end{proposition}

Note that modifying $\Phi$ allows to slightly change its level sets. For instance, taking $\Phi_\e(x) = \Phi(x)-\e |x-x_0|^2$ (which shall be very useful for applications to unique continuation), the level set $\{\Phi_{\e}=0\}$ is slightly bent (except at $x_0$) into the set $\{\Phi>0\}$ (where $u$ will be assumed to be zero). This slight change will be crucial for the proof of the unique continuation theorem. 
\bnp
First, we notice that we can prove as in the proof of Proposition~\ref{p:pseudo-hormander} (still using Lemma \ref{l:fgh-three-fcts} combined with Lemma \ref{lmlien} for the limit $\tau\to0^+$) that Definition~\ref{defconvexfunctions} implies (and is actually equivalent to) the existence of an inequality of the form
\bna
c_{\Phi}(\xi,\tau)+C_1 \frac{\left|p_{\Phi}(x_0,\xi,\tau)\right|^2}{|\xi|^2+\tau^2}\geq C_2(|\xi|^2+\tau^2), 
\ena
uniformly for $(\xi,\tau)$ with $|\xi|^2+\tau^2=1$, $\tau\geq 0$ (see Lemma~\ref{lmlien}), where
$$c_{\Phi}(\xi,\tau)=\frac{1}{i\tau}\{\overline{p_{\Phi}},p_{\Phi}\}(x_0,\xi,\tau), \text{ for } \tau >0 \quad \text{ and }  \quad c_{\Phi}(\xi,0)= 2\{p_2,\{p_2 ,\Phi\}\}(x_0, \xi) . $$

We then remark that all quantities in the above estimate only involve derivatives of $\Phi$ of order at most $2$ (as a consequence of Lemma \ref{lmlien}) at the point $x_0$. It is therefore stable by the addition of a function small for the $C^2$ norm around $x_0$.  
\enp

\subsection{Strongly pseudoconvex hypersurfaces} 
\label{s:strong-pseudo-surfaces}

Until this point, we have proved a Carleman estimate with weight $\Phi$ provided $\Phi$ satisfies a (weird?) pseudoconvexity condition (Definition \ref{defconvexfunctions}). The main purpose of this section is to provide a geometric characterization of hypersurfaces $S$ for which we can find a function $\Phi$ having $S$ as a level set and being appropriate for the Carleman estimate (i.e. satisfying Definition \ref{defconvexfunctions}).
We shall eventually prove that the sought geometric condition on the hypersurface $S$ is~\eqref{pseudoconvex2def}. We first introduce a seemingly stronger condition.

\begin{definition}[Usual pseudoconvexity for hypersurfaces]
\label{def-pseudosurface}
Let $\Omega \ni x_0$ be an open set, $P\in \difclas{2}(\Omega)$ with real-valued principal symbol $p_2$ and $\Psi\in C^{\infty}(\Omega)$ real-valued. 
We say that the {\em oriented hypersurface} $S = \{\Psi = \Psi(x_0)\} \ni x_0$ is strongly pseudoconvex with respect to $P$ at $x_0$ if
 \bnan
\label{e:pseudo-surface-1}
\left\{ p_2, \{p_2 ,\Psi \} \right\} (x_0 ,\xi) >0  , 
&\text{ if } p_2(x_0,\xi) =\{p_2,\Psi\}(x_0,\xi)= 0 \text{ and } \xi \neq 0 ; \\
\label{e:pseudo-surface-2}
\frac{1}{i\tau}\{ \overline{p_\Psi},p_\Psi \} (x_0 ,\xi, \tau) >0  , 
&\text{ if } p_\Psi(x_0,\xi, \tau) =\{p_\Psi,\Psi\}(x_0,\xi, \tau)= 0 \text{ and } \tau >0 ,
\enan
where $p_\Psi(x,\xi, \tau) = p_2(x, \xi + i \tau d\Psi(x))$.
\end{definition}
Note that the first condition~\eqref{e:pseudo-surface-1} is precisely~\eqref{pseudoconvex2def}.
We shall eventually prove that for differential operators of order two with real principal symbols,~\eqref{e:pseudo-surface-1} implies~\eqref{e:pseudo-surface-2}.
Note that the definition seems to depend on the defining function $\Psi$ for the hypersurface $S$, and not only on the oriented hypersurface $S$ itself. Lemma~\ref{l:defining-function} (se e.g.~\cite{LL:book} for a proof) shows this is not the case, and hence justifies the definition.
\begin{lemma}
\label{l:defining-function}
Assume $S =  \{\Psi_1 = \Psi_1(x_0)\} = \{\Psi_2 = \Psi_2 (x_0)\}$ with $d\Psi_j (x_0)\neq 0$, $j=1,2$ and $d\Psi_2 (x_0) =  \lambda d\Psi_1 (x_0)$ for some $\lambda>0$ (same orientation). Then $\Psi_1$ satisfies~\eqref{e:pseudo-surface-1} if and only if $\Psi_2$ satisfies~\eqref{e:pseudo-surface-1}, and $\Psi_1$ satisfies~\eqref{e:pseudo-surface-2} if and only if $\Psi_2$ satisfies~\eqref{e:pseudo-surface-2}.
\end{lemma}

Remark that Definition~\ref{def-pseudosurface} looks very similar to Definition~\ref{defconvexfunctions}. It is just slightly weaker because the positivity condition is assumed only under the additional conditions $\{p_2,\Phi\}=0$ and $\{p_\Phi,\Phi\}=0$.
In particular, the level sets of a pseudoconvex functions are pseudoconvex oriented hypersurfaces. This is however not useful since Definition~\ref{defconvexfunctions} is not geometric (but rather linked to Carleman estimates).

The importance of Definition~\ref{def-pseudosurface} is twofold:
\begin{itemize}
\item It is a purely geometric definition: this comes from Lemma~\ref{l:defining-function} and the fact that Conditions~\eqref{e:pseudo-surface-1}-\eqref{e:pseudo-surface-2} are invariant by diffeomorphisms.
\item Once $\Psi$ satisfies this geometric condition, one can produce a function $\Phi$ having the same levelsets (hence keeping the geometry unchanged), and that satisfies the stronger pseudoconvexity condition of Definition~\ref{defconvexfunctions}.
This is the goal of the next section.
\end{itemize}

Note that, once again, Condition~\eqref{e:pseudo-surface-1} (on the real domain) is the limit as $\tau \to 0^+$ of Condition~\ref{e:pseudo-surface-2} (on the complex domain). This follows both from Lemma~\ref{lmlien} and the fact that
\begin{align}
\label{e:bracket-limit}
\{p_\Psi,\Psi\}(x,\xi, \tau) & = \d_\xi \big(p(x, \xi + i \tau d\Psi(x) \big) \cdot \d_x \Psi (x) = ( \d_\xi p )(x, \xi + i \tau d\Psi(x)) \cdot \d_x \Psi (x) \nonumber \\
& = \{p_2 , \Psi\} (x, \xi + i \tau d\Psi(x)) \to \{p_2 , \Psi\} (x, \xi ),  \quad \text{ as } \tau\to 0^+ .
\end{align}

\subsubsection{(Analytic) convexification}
\label{Convexification}

\begin{proposition}[Analytic convexification]
\label{propanalyticconvex}
Let $\Omega \ni x_0$ be an open set, $P\in \difclas{2}(\Omega)$ with real-valued principal symbol $p_2$ and $\Psi\in C^{\infty}(\Omega)$ real-valued. 
Assume the oriented hypersurface $S = \{\Psi = \Psi(x_0)\}$ is strongly pseudoconvex with respect to $P$ at $x_0$ (Definition~\ref{def-pseudosurface}).
Then there exists $\lambda_0>0$ such that for all $\lambda \geq \lambda_0$, the function $\Phi=e^{\lambda \Psi}$ is pseudoconvex with respect to $P$ at $x_0$ (Definition~\ref{defconvexfunctions}).
\end{proposition}
Hence, the Carleman estimate of Theorem~\ref{thmCarlclass} holds with weight $\Phi$.
Note that the geometry of the level-sets of $\Phi$ and $\Psi$ are actually the same: only the values of the level sets of $\Phi$ are stretched. 
Here, for any strongly pseudoconvex oriented hypersurface $S=\{\Psi=\Psi(x_0)\}$, this proposition produces an admissible Carleman weight (that is, a pseudoconvex function) $\Phi$ having exactly the same level sets.
In order to simplify the notation for the proof, we recall that $x_0$ is fixed and remark that changing the function $\Psi$ by a constant does not change the assumption. We may thus assume that 
\bnan
\label{e:values-psiphix0} 
\Psi(x_0) = 0 , \quad \text{ and hence } \quad \Phi(x_0) = 1  \quad \text{ and } \quad  d \Phi(x_0) = \lambda d\Psi(x_0) .
\enan
We also denote 
\begin{align}
\label{e:def-c-psi}
c_{\Psi}(\xi,\tau)=\frac{1}{i\tau}\{\overline{p_{\Psi}},p_{\Psi}\}(x_0,\xi,\tau), \text{ for } \tau >0 \quad \text{ and }  \quad c_{\Psi}(\xi,0)= 2\{p_2,\{p_2 ,\Psi\}\}(x_0, \xi) , 
\end{align}
 with a similar definition for $c_{\Phi}(\xi,\tau)$. According to Lemma~\ref{lmlien}, $c_{\Psi}(\xi,\tau)$ and $c_{\Phi}(\xi,\tau)$ are continuous on the whole $\R^n \times \R^+$.
The proof of Proposition~\ref{propanalyticconvex} is then based on the following computation.

\begin{lemma}
\label{l:calcul-cphi-cpsi}
Assume $\Phi=e^{\lambda \Psi}$. For all $(\xi,\tau) \in \R^n \times \R^+$ and all $\lambda>0$, we have 
\bna
c_{\Phi}(\xi,\tau)= \lambda c_{\Psi}(\xi,\lambda \tau)+2\lambda^2 \left|\{p_{\Psi},\Psi\}(x_0,\xi,\lambda \tau)\right|^2.
\ena
\end{lemma}
We first prove the proposition from the lemma and then prove the lemma.
\bnp[Proof of Proposition~\ref{propanalyticconvex} from Lemma~\ref{l:calcul-cphi-cpsi}]
Using Lemma~\ref{l:fgh-three-fcts} (combined with Lemma \ref{lmlien} and \eqref{e:bracket-limit} in the limit $\tau \to 0^+$), Properties \eqref{e:pseudo-surface-1}-\eqref{e:pseudo-surface-2} imply the existence of $C_1, C_2>0$ so that
\bna
c_{\Psi}(\xi,\tau)+C_1\left|\{p_{\Psi},\Psi\}(x_0,\xi,\tau)\right|^2  +C_1 \frac{\left|p_{\Psi}(x_0,\xi,\tau)\right|^2}{|\xi|^2+\tau^2}\geq C_2(|\xi|^2+\tau^2).
\ena
for any $\tau\geq 0$, $|\xi|^2+\tau^2=1$ (note that this takes into account the limit $\tau \to 0^+$).
Replacing $\tau$ by $\lambda\tau$ for $\lambda\geq 1$ and using homogeneity, this can be reformulated as
\bnan
\label{e:interm-estim-cphi}
c_{\Psi}(\xi,\lambda\tau)+C_1\left|\{p_{\Psi},\Psi\}(x_0,\xi,\lambda\tau)\right|^2  +C_1 \frac{\left|p_{\Psi}(x_0,\xi,\lambda\tau)\right|^2}{|\xi|^2+\lambda^2\tau^2}\geq C_2(|\xi|^2+\lambda^2\tau^2).
\enan
for any $(\xi,\tau)\neq (0,0)$ with $\tau\geq 0$. 
Moreover, using Lemma~\ref{l:calcul-cphi-cpsi} and noticing (see~\eqref{e:values-psiphix0}) that 
$$
p_{\Psi}(x_0,\xi,\lambda\tau) = p_2 (x_0 , \xi + i \lambda \tau d \Psi(x_0))= p_2 (x_0 , \xi + i \tau d \Phi(x_0)) = p_{\Phi}(x_0,\xi,\tau),
$$
 we obtain
\begin{align*}
c_{\Phi}(\xi,\tau)+C_1 \lambda \frac{\left|p_{\Phi}(x_0,\xi,\tau)\right|^2}{|\xi|^2+\tau^2}
&= \lambda \left( c_{\Psi}(\xi,\lambda \tau)+2\lambda \left|\{p_{\Psi},\Psi\}(x_0,\xi,\lambda \tau)\right|^2 +C_1 \frac{\left|p_{\Psi}(x_0,\xi,\lambda \tau)\right|^2}{|\xi|^2+\tau^2}\right) .
\end{align*}
Now taking $\lambda \geq \max \{ C_1/2 ,  1\}$ and using~\eqref{e:interm-estim-cphi} yields
\begin{align*}
c_{\Phi}(\xi,\tau)+C_1 \lambda \frac{\left|p_{\Phi}(x_0,\xi,\tau)\right|^2}{|\xi|^2+\tau^2} 
&\geq \lambda \left(c_{\Psi}(\xi,\lambda \tau)+C_1 \left|\{p_{\Psi},\Psi\}(x_0,\xi,\lambda \tau)\right|^2+C_1 \frac{\left|p_{\Psi}(x_0,\xi,\lambda \tau)\right|^2}{|\xi|^2+\lambda^2\tau^2}\right)\\
&\geq C_2\lambda(|\xi|^2+\lambda^2\tau^2)\geq C_2\lambda(|\xi|^2+\tau^2).
\end{align*}
When recalling the definition of $c_{\Phi}$, this readily implies~\eqref{e:pseudo-function-2}, and also~\eqref{e:pseudo-function-1} in the limit $\tau \to 0^+$ (with Lemma~\ref{lmlien}). This concludes the proof that $\Phi$ is pseudoconvex for $P$ at $x_0$ in the sense of Definition~\ref{defconvexfunctions}.
\enp
\bnp[Proof of Lemma~\ref{l:calcul-cphi-cpsi}]
We compute
\begin{align*}
\partial_j \Phi&=\lambda \partial_j \Psi e^{\lambda\Psi} , \quad 
\partial_{j,k} \Phi=\lambda \partial_{j,k} \Psi e^{\lambda\Psi}+\lambda^2 (\partial_j \Psi)(\partial_k \Psi) e^{\lambda\Psi} , 
\end{align*}
which we write in a shorter way as
\begin{align*}
d\Phi&=\lambda e^{\lambda\Psi} d\Psi  , \qquad 
\Hess(\Phi)(\xi,\eta)=\lambda \Hess(\Psi)(\xi,\eta)e^{\lambda\Psi} +\lambda^2 (\xi\cdot \d_x\Psi)(\eta\cdot \d_x\Psi)e^{\lambda\Psi} .
\end{align*}
Taken at the point $x_0$, and  recalling~\eqref{e:values-psiphix0}, this implies
\begin{align}
\label{e:tototo}
d\Phi(x_0)&=\lambda d\Psi(x_0) , \qquad 
\Hess(\Phi)(x_0)(\xi,\eta)=\lambda \Hess(\Psi)(x_0)(\xi,\eta) +\lambda^2 (\xi\cdot \d_x\Psi(x_0))(\eta\cdot \d_x\Psi(x_0)) .
\end{align}
We now compute
\begin{align*}
c_{\Phi}(\xi,\tau)  = \frac{1}{i\tau}\{\overline{p_{\Phi}},p_{\Phi}\}(x_0,\xi,\tau) 
&=\frac{1}{i\tau} \d_{\xi}p(x_0,\xi-i\tau d\Phi(x_0))\cdot \d_x p(x_0,\xi+i\tau d\Phi(x_0)) \\
&\quad+ \Hess(\Phi)(x_0)\left[\d_{\xi}p(x_0,\xi-i\tau d\Phi(x_0)),\d_{\xi}p(x_0,\xi+i\tau d\Phi(x_0))\right]\\
&\quad -\frac{1}{i\tau}\d_{x}p(x_0,\xi-i\tau d\Phi(x_0))\cdot \d_{\xi} p(x_0,\xi+i\tau d\Phi(x_0))\\
&\quad + \Hess(\Phi)(x_0)\left[\d_{\xi}p(x_0,\xi-i\tau d\Phi(x_0)),\d_{\xi}p(x_0,\xi+i\tau d\Phi(x_0))\right]\\
&=\frac{2}{\tau}\Im \left[\d_{\xi}p(x_0,\xi-i\tau d\Phi(x_0))\cdot\d_x p(x_0,\xi+i\tau d\Phi(x_0))\right]\\
&\quad +2 \Hess(\Phi)(x_0)\left[\d_{\xi}p(x_0,\xi-i\tau d\Phi(x_0)),\d_{\xi}p(x_0,\xi+i\tau d\Phi(x_0))\right] .
\end{align*}
Using now~\eqref{e:tototo}, this rewrites (we drop from the notation the fact that $\Psi$ and its derivatives are taken at $x_0$)
\begin{align*}
c_{\Phi}(\xi,\tau)  
&=\frac{2}{\tau}\Im\left[\d_{\xi}p(x_0,\xi-i\tau\lambda d\Psi)\cdot \left(\d_x p(x_0,\xi+i\tau\lambda d\Psi)\right)\right]\\
&+2\lambda \Hess(\Psi)\left[\d_{\xi}p(x_0,\xi-i\tau\lambda d\Psi),\d_{\xi}p(x_0,\xi+i\tau\lambda d\Psi)\right]\\
&+2 \lambda^2 \left(\d_{\xi}p(x_0,\xi-i\tau\lambda d\Psi)\cdot \d_x \Psi\right) \left(\d_{\xi}p(x_0,\xi+i\tau\lambda d\Psi)\cdot \d_x\Psi \right)\\
&=\lambda c_{\Psi}(\xi,\lambda \tau)+2\lambda^2 \left|\{p,\Psi\}(x_0,\xi+i\tau\lambda d\Psi)\right|^2\\
&=\lambda c_{\Psi}(\xi,\lambda \tau)+2\lambda^2 \left|\{p_{\Psi},\Psi\}(x_0,\xi,\lambda \tau)\right|^2,
\end{align*}
proving the lemma.
\enp
 
\subsubsection{Reducing the strong pseudoconvexity assumption to the condition on the real space}
In the particular case of differential operators of order two, with real principal symbol, Condition~\eqref{e:pseudo-surface-1} on the real space implies Condition~\eqref{e:pseudo-surface-2} in the complex space (this is no longer the case if one of these two conditions is not satisfied, see~\cite{Hoermander:V4,Lerner:book-carleman}). That is to say, Definitions~\ref{defconvexsurface} and~\ref{def-pseudosurface} are equivalent (differential operators of order two with real principal symbol).
\begin{proposition}
\label{proppseudosurface}
Let $\Omega \ni x_0$ be an open set, $P\in \difclas{2}(\Omega)$ with real-valued principal symbol $p_2$ and $\Psi\in C^{\infty}(\Omega)$ real-valued. 
Assume that the oriented hypersurface $S=\{\Psi = \Psi(x_0)\}$ satisfies Condition~\eqref{e:pseudo-surface-1} at $x_0$. 
Then $S=\{\Psi = \Psi(x_0)\}$ is strongly pseudoconvex with respect to $P$ at $x_0$ (i.e. both conditions~\eqref{e:pseudo-surface-1} and~\eqref{e:pseudo-surface-2} are satisfied).
\end{proposition}
We split the proof of Proposition~\ref{proppseudosurface} into two lemmata, concerned with the non-characteristic case ($p_2(x_0, d\Psi(x_0))\neq 0$) and the characteristic case ($p_2(x_0, d\Psi(x_0))= 0$), respectively.

\begin{lemma}
\label{l:non-caract-case}
Assume $p_2$ is a real symbol of order two near $x_0$, and $\Psi$ is such that $p_2(x_0, d\Psi(x_0))\neq 0$. Then, for any $\xi \in \R^n$ we have 
\begin{align}
\label{e:asspt-non-caract}
p_\Psi (x_0, \xi , \tau) = \{p_\Psi , \Psi \}(x_0,\xi)  = 0 \quad  \implies \tau =0 .
\end{align}
\end{lemma}
In this case, Assumption~\eqref{e:pseudo-surface-2}  is thus empty.
\begin{lemma}
\label{l:caract-case}
Assume $p_2$ is a real symbol of order two near $x_0$, and $\Psi$ is such that $p_2(x_0, d\Psi(x_0)) = 0$. Assume also~\eqref{e:pseudo-surface-1} for all $\xi \in \R^n\setminus \{0\}$. Then we also have~\eqref{e:pseudo-surface-2}.
\end{lemma}
Both proofs of Lemmata~\ref{l:non-caract-case} and \ref{l:caract-case} rely on the fact that for fixed $\xi\in \R^n$,  
$$
f(z)=p_2(x_0,\xi+ z d\Psi(x_0)) =  p_2(x_0,\xi) + z^2 p_2(x_0, d\Psi(x_0))+ 2 z \tilde{p}_2 (x_0,\xi , d\Psi(x_0)) ,
$$
is a {\em second order polynomial} in the variable $z \in \C$, with {\em real coefficients}. 
Moreover, the assumption of~\eqref{e:asspt-non-caract} (resp. of \eqref{e:pseudo-surface-2}) implies that
\begin{align*}
f(i\tau) & =p_2(x_0,\xi+i\tau d\Psi(x_0)) = p_\Psi (x_0, \xi , \tau) = 0 \quad \text{ and }\\
f'(i\tau)& = \d_\xi p_2(x_0,\xi+i\tau d\Psi(x_0)) \cdot \d_x \Psi (x_0) = \{p_2 , \Psi \} (x_0,\xi+i\tau d\Psi(x_0)) = \{p_\Psi , \Psi \}(x_0,\xi) = 0 ,
\end{align*}
that is to say, $z=i \tau$ ($\tau \in \R^+$) is a double root of the polynomial $f$.

\bnp[Proof of Lemma~\ref{l:non-caract-case}]
Since the coefficient in front of $z^2$, namely $p_2(x_0, d\Psi(x_0))$ is non-zero, the polynomial $f$ has two complex roots which are either both in $\R$, or complex conjugate.
That $z=i \tau$ ($\tau \in \R^+$) is a double root of the polynomial $f$ implies $\tau=0$.
\enp

The proof of Lemma~\ref{l:caract-case} relies on tedious computations, and we refer the reader to~\cite{LL:book}. Note that so far, we have given a complete proof of Theorem~\ref{thmreal2} under the additional non-characteristicity condition $p_2(x_0, d\Psi(x_0)) \neq 0$.

\subsubsection{Unique continuation: end of proof of Theorem~\ref{thmreal2}}
\label{sectUCP}
In this section, we conclude the proof of Theorem~\ref{thmreal2}. 
After a geometric convexification procedure, it consists essentially in using Lemma~\ref{l:uniqueness-from-carleman}. 
\bnp[Proof of Theorem \ref{thmreal2}] 
We first remark that we may assume that $\Psi(x_0)=0$ (up to changing $\Psi$ into $\Psi -\Psi(x_0)$, which does not change the assumption), so that $S= \{\Psi=0\}$. 
Let $u$ be a $C^\infty$ solution of $Pu=0$ in $\Omega$ so that $u=0$ on $\Omega\cap \left\{\Psi> 0\right\}$. The hypersurface $S=\{\Psi= 0\}$ being strongly pseudoconvex at $x_0$, Proposition \ref{propanalyticconvex} shows that for $\lambda$ large enough (but fixed) $\Phi := e^{\lambda \Psi}-1$ is a pseudoconvex function with $\{\Phi =  0\}=\{\Psi = 0\}$, $\{\Phi >  0\}=\{\Psi > 0\}$ and $\{\Phi < 0\}=\{\Psi < 0\}$. 
Proposition \ref{propgeometricconvex} yields the existence of $\e>0$, such that $\Phi_\eps=\Phi-\e|x-x_0|^2$ remains a pseudoconvex function (Definition~\ref{defconvexfunctions}). As a consequence of Proposition~\ref{p:pseudo-hormander} and Theorem~\ref{thmCarlclass}, it therefore satisfies the following properties
\begin{enumerate}
\item \label{enumCarl}there exist $R>0$, $C>0$ and $\tau_0>0$ so that we have the following estimate
\bnan
\label{Carlproof}
 \tau^3 \nor{e^{\tau \Phi_\eps} w}{L^2}^2 +\tau \nor{e^{\tau \Phi_\eps}\nabla w}{L^2}^2 \leq C \nor{e^{\tau \Phi_\eps}P w}{L^2}^2,  
 \enan 
for any $w\in C^{\infty}(B(x_0,R))$ and $\tau\geq \tau_0$.
\item \label{ecart} there exists $\eta>0$ so that $\Phi_\eps(x)\leq -\eta$ for $x\in \{\Phi\leq 0\}\cap \{|x-x_0|\geq R/2\}$,
\item \label{pres}there exists a neighborhood $V\subset B(x_0,R/2)$ of $x_0$ so that $\Phi_\eps(x)\geq -\eta/2$ for $x\in V$. 
\end{enumerate}
Property \ref{enumCarl} is a consequence of Theorem \ref{thmCarlclass}, and $R$ is fixed by that theorem.
Property \ref{ecart} is true thanks to the parameter $\e$ in the geometric convexification. Indeed, for $|x-x_0|\geq R/2$, we have $\Phi_\eps(x)\leq \Phi(x)-\e R^2/4$. If $\Phi(x)\leq 0$, this implies $\Phi_\eps(x)\leq -\e R^2/4$, so that we can take $\eta=-\e R^2/4$.
Property \ref{pres} is only a continuity argument since $\Phi_\eps(x_0)=0$.
From this point forward, it suffices to follow the strategy described in Section~\ref{subsectstratCarle} to conclude the proof of Theorem~\ref{thmreal2}.

Note finally that, in order for the result to hold for $u \in H^1(\Omega)$, we need to remark that a density argument shows that the Carleman estimate is still valid for all $w \in H^1(\Omega)$ such that $\supp(w) \subset B(x_0,R)$ and $Pw \in L^2(\Omega)$. Here, in case $u \in H^1(\Omega)$ with $Pu=0$, we have $w = \chi u \in H^1(\Omega)$ with $\supp(w) \subset \supp(\chi) \subset B(x_0,R)$ and $Pw = 0 + [P,\chi]u \in L^2(\Omega)$ since $[P,\chi] \in \difclas{1}(\Omega)$ and $u \in H^1(\Omega)$. Hence, the Carleman estimate applies and the remainder of the proof remains unchanged.
\enp

\subsection{Necessity of strong pseudoconvexity for stable unique continuation}
\label{s:Alinhac-baouendi}
In this section, we discuss optimality/limitations of the H\"ormander's theorem for the wave operator, via two counterexamples due to Alinhac and Alinhac-Baouendi, respectively. We recall (see Remark~\ref{r:lot}) that H\"ormander's theorem is insensitive to addition of lower order terms to the operator.
The following result is a particular case of~\cite[Th\'eor\`eme~2]{Alinhac:83}.
\begin{theorem}[Alinhac]
\label{t:alinhac}
Let $\Omega$ be an open subset of $\R^n$ let $x_0 \in \Omega$, and let $P \in \difclas{2}(\Omega)$ with real principal symbol $p_2(x,\xi) = a^{ij}(x) \xi_i\xi_j$. Assume $a^{ij}$ is a symmetric real-valued matrix defined in a neighborhood of $x_0$ and such that $(a^{ij}(x_0))$ is non-degenerate.
Let $\Psi \in C^\infty(\Omega;\R)$ be such that $\Psi(x_0)=0$ and $p_2(x_0, d\Psi(x_0)) \neq 0$,  and assume that there exists $\xi_0 \in \R^n\setminus \{0\}$ such that 
\begin{equation}
\label{e:non-strong-pseudo}
p_2(x_0,\xi_0)=\{ p_2 ,\Psi\}(x_0, \xi_0)=0 \quad \text{ and } \quad  \{p_2,\{p_2,\Psi\}\}(x_0, \xi_0) < 0  .
\end{equation}
Then, there exist $U \subset \Omega$ a neighborhood of $x_0$ and $q, u \in C^\infty(U;\C)$ such that 
\begin{align*}
Pu+qu = 0 \text{ in } U , \quad \text{ and } \quad \supp(u) =\{\Psi \geq 0\} \cap U.
\end{align*}
\end{theorem}
In particular, under Assumption~\eqref{e:non-strong-pseudo} on the oriented hypersurface $\{\Psi=0\}$, unique continuation from $\{\Psi > 0\}$ to a neighborhood of $x_0$ does not hold for the operator $P+q$.
This applies to the wave operator. Condition~\eqref{e:non-strong-pseudo} is a strong negation of strong pseudoconvexity (Definition~\ref{defconvexsurface})
and thus Theorem~\ref{t:alinhac} is an almost converse to Theorem~\ref{thmreal2} if one consider ``stable unique continuation'' for $P$, that is to say unique continuation for all zero-order perturbations of $P$.

Note that Assumption~\eqref{e:non-strong-pseudo} can be reformulated on the tangent space as in Remark~\ref{r:waves-tangent}.
If we denote by $a(x) = a_{ij}(x) = (a^{ij}(x))^{-1}$ the pre-dual (pseudo-Riemannian) metric as in Remark~\ref{r:waves-tangent}, 
then Assumption~\eqref{e:non-strong-pseudo} is equivalent to: there exists $X_0 \in T_{x_0}\Omega \setminus \{0\}$ such that 
$$
(X_0, X_0)_a =0 , \quad  d\Psi(x_0)(X_0)=0 \quad \text{ and } \quad  \Hess_a \Psi(x_0) (X_0 , X_0) <0  .
$$

The following result is another counterexample to stable unique continuation in the limit case where  $\{p_2,\{p_2,\Psi\}\}(x_0, \xi_0) =0$.
It is a particular case of~\cite[Th\'eor\`eme~2]{AB:79},~\cite[Theorem]{AB:95}.
\begin{theorem}[Alinhac-Baouendi]
\label{t:alinhac-baouendi}
Assume $d\geq 2$ and consider $P := D_t^2- D_{x_1}^2 -D_{x_2}^2-\cdots - D_{x_d}^2$ near the point $0\in \R^{1+d}$. There is an open set $U\subset \R^{1+d}$ with $0 \in U$, there exist $q, u \in C^\infty(U;\C)$ such that 
\begin{align*}
Pu+qu = 0 \text{ in } U , \quad  \text{ and } \quad  \supp(u) =\{x_1\geq 0\} \cap U.
\end{align*}
\end{theorem}
This implies that unique continuation fails for the operator $P+q$ across the hypersurface $\{x_1=0\}$, even though $P$ has constant coefficients (hence the Holmgren-John theorem~\ref{thmholmgren} applies to $P$ across $\{x_1=0\}$) and the perturbation $q$ is of lower order and smooth.
Note that the principal symbol of $P+q$ is $p_2(t,x,\xi_t,\xi_x) = \xi_t^2 - \xi_{x_1}^2-\xi_{x_2}^2-\cdots - \xi_{x_d}^2$ and the hypersurface $\{x_1=0\}$ barely fails to be strongly pseudoconvex (see Definition~\ref{defconvexsurface}). Indeed, we have $\{p_2,x_1\} = -2\xi_{x_1}$ and $\{p_2 , \{p_2,x_1\} \}= 0$ so that, if one chooses $\xi_0 := (1,0, 1, 0, \cdots , 0)$ (that is to say $\xi_t=\xi_{x_2}=1$, $\xi_{x_1}=0$ and $\xi_{x_j}=0$ for $j \in \{3,\cdots ,d\}$), we have 
\begin{equation*}
p_2(0,\xi_0)=0, \quad \{ p_2 ,x_1 \}(0, \xi_0)=0 \quad \text{ and } \quad  \{p_2,\{p_2,x_1\}\}(x_0, \xi_0) = 0  .
\end{equation*}
In some sense, this is a weaker form of violation of strong pseudoconvexity (see Definition~\ref{defconvexsurface}) compared to~\eqref{e:non-strong-pseudo}.

\section{Unique continuation for time-independent wave operators}
\label{chapterwave}
To summarize the discussion so far, if one considers a general wave operator $P = D_i a^{ij}D_j+$ lower order terms, we have on the one hand the Holmgren-John Theorem~\ref{thmholmgren}: we assume analyticity of all coefficients and unique continuation holds across any noncharacteristic hypersurface. The latter geometric condition appears to be the appropriate one in applications (see Sections~\ref{s:semiglobal} and~\ref{s:GLUC} below) and is essentially optimal, whereas analyticity is a very demanding condition.
On the other hand the H\"ormander Theorem~\ref{thmreal2} assumes only $C^\infty$ regularity of the principal part of the operator and $L^\infty$ regularity of the lower order terms (these regularity assumptions can even be relaxed, see~\cite{Hoermander:63}, but we do not discuss this issue here) and yields unique continuation across any strongly pseudoconvex hypersurface. The regularity assumption is suitable for applications, but the geometric pseudoconvexity condition is extremely demanding (see e.g. the geometric discussion in Section~\ref{rkconvexflow} and the examples in Section~\ref{s:examples-hor}). As explained in Section~\ref{s:Alinhac-baouendi}, it is however optimal if one considers stable unique continuation.

In the present chapter, we explore further the case where all coefficients of the wave operator are {\em time-independent} with only finite regularity in space. The result we present has a long history with several ancestors and descendents.  A historical account is provided in Section~\ref{s:history}. One important fact noticed along the way is the role of time analyticity (which obviously holds in the time-independent case). Here, we
 focus our attention to the operator $\partial_t^2 + Q$ where $Q= q(x, D_x)$ is a positive elliptic operator. 
The result presented in this chapter has been proved in~\cite{Tataru:95} and our presentation is inspired by~\cite{Hor:97}.


\subsection{Setting and statement of the unique continuation result}
In the following, we denote the generic variable by $\x=(t,x)\in\R^{1+d}=\R^n$ (note the slight difference with the notation in Section~\ref{chapterclassical}, where $x\in \R^n$ includes the time variable) with dual variable $\xi=(\xi_t,\xi_x)\in\R^{1+d}$.
The main theorem of this chapter is as follows.
\begin{theorem}[Wave type operator with time-independent coefficients]
\label{thmwave}
Let $T>0$ and $\Omega_x$ an open set of $\R^d$. Denote $\Omega=]-T,T[\times \Omega_x$.
Let $$Q=\sum_{i,j=1}^d g^{ij}(x)D_i D_j +\sum_{k=1}^d b_k(x)D_k +c(x)$$ be a differential operator of order $2$ with $g^{ij}\in C^{\infty}(\Omega_x)$ real-valued, $b_k$, $c\in L^{\infty}(\Omega_x)$. Assume also that $Q$ is positive elliptic, i.e. there exists $C>0$ so that 
\begin{align}
\label{ellipticity}
 q(x,\xi_x) := \sum_{i,j=1}^d g^{ij}(x)\xi_i \xi_j\geq C |\xi_x|^2, \textnormal{ for all } (x,\xi_x)\in \Omega_x \times \R^d.
\end{align}
Let $P=\partial_t^2 + Q$ on $\Omega$ and set $p_2(t,x,\xi_t,\xi_x): = -\xi_t^2 + q(x,\xi_x)$.
Let $\x_0=(t_0,x_0)\in \Omega$ and $\Psi\in C^{2}(\Omega)$ with $d\Psi(\x_0)\neq 0$ so that $p_2(\x_0,d\Psi(\x_0)) \neq 0$, i.e.
\bna
(\partial_t \Psi(\x_0))^2\neq \sum_{i,j} g^{ij}(x_0)(\partial_i\Psi(\x_0))(\partial_j\Psi(\x_0)).
\ena
 Then, there exists a neighborhood $V$ of $\x_0$ so that for any $u\in H^1(\Omega)$,
\bnan
\left\{\begin{array}{rcl}
Pu &=& 0  \textnormal{ in  } \Omega,\\ u&=&0 \textnormal{ in } \Omega\cap \left\{\Psi> \Psi(\x_0)\right\}\end{array}\right.\Longrightarrow u= 0 \textnormal{ on } V.
\enan
\end{theorem}

Note that there is no link between the strong pseudoconvexity condition in Definition~\ref{defconvexsurface} and the non-characteristicity condition in Theorem~\ref{thmwave}. Many useful hypersurfaces are non-characteristic but not strongly pseudoconvex (see Section~\ref{s:semiglobal} below), but one may also construct characteristic hypersurfaces that are strongly pseudoconvex (Indeed, the former condition is a first order condition whereas the latter is a second order condition).
 The result of Theorem~\ref{thmwave} actually holds under the weaker condition that the hypersurface is $\{\Psi=0\}$ is strongly pseudoconvex in $\xi_t=0$ (see Definition~\ref{def-pseudosurface-xit=0} below).
However, in all applications we have in mind, only the non-characteristicity condition is useful to deduce optimal results. See see Sections~\ref{s:semiglobal} and~\ref{s:GLUC} below.

\begin{remark}
As in Remark~\ref{r:waves-tangent} concerning the pseudoconvexity condition, the non-characteristicity condition $p_2(\x_0,d\Psi(\x_0)) \neq 0$, formulated here as a cotangent condition, may be equivalently rephrased on the tangent space. Denoting again by $a=(a_{ij})$ the metric on the tangent space $T((-T,T)\times \Omega_x)$, the hypersurface $S=\{\Psi=0\}$ is non-characteristic at $\x_0 \in S$ iff $(\nabla_a \Psi, \nabla_a \Psi)_a (\x_0) \neq 0$. Notice that under the assumptions of Theorem~\ref{thmwave}, the metric $a$ is assumed of the (particular) form $a(X,Y)= (X,Y)_a = -X_tY_t + g(X_x,Y_x)$ where $g$ is a time-independent Riemannian metric and $X=(X_t,X_x), Y=(Y_t,Y_x)$.
\end{remark}

The proof of Theorem~\ref{thmwave} relies on an inequality of Carleman type, but with an additional weight in the Fourier variable.
Namely, we let $e^{-\e\frac{|D_t|^2}{2\tau}}$ be the Fourier multiplier defined naturally by
$$
\F  \left(e^{-\e\frac{|D_t|^2}{2\tau}}u\right)(\xi)=e^{-\e\frac{|\xi_t|^2}{2\tau}}\widehat{u}(\xi) , \quad u\in\mathcal{S}(\R^{1+d}) ,
$$
where $\xi_t$ is the Fourier variable corresponding to the variable $t$ and $\xi= (\xi_t,\xi_x)$. Note that this amounts to solving the heat equation with $t$ as a ``spatial'' variable, during a ``time'' $\frac{\e}{2\tau}$. 
Using the explicit expression of the Fourier transform of the Gaussian $e^{-\e\frac{|\xi_t|^2}{2\tau}}$, this may be rewritten as a convolution with a heat kernel:
\begin{equation}
\label{e:Fourier-gaussienne}
\left(e^{-\e\frac{|D_t|^2}{2\tau}}u \right) (t,x) = \left(\frac{\tau}{2\pi \eps} \right)^{1/2} \int_\R e^{-\frac{\tau}{2\eps} (t-s)^2} u(s,x) ds .
\end{equation}
This operator has several interesting features: it localizes close to $D_t=0$ (i.e. in low frequencies w.r.t. the time variable $t$), in an analytic way (the function $e^{-\e\frac{|D_t|^2}{2\tau}}u$ produced is an entire function in the $t$-variable). However (and consequently), note that $e^{-\frac{\e}{2\tau}|D_t|^2}$ is not local; in particular, $e^{-\e\frac{|D_t|^2}{2\tau}}u$ is not compactly supported, even if $u$ is.

For a smooth real-valued weight function $\Phi$ (later on, we will assume that it is polynomial of order $2$), the Carleman estimate below will make use of the operator
\bna
Q_{\e,\tau}^\Phi u:=e^{-\e\frac{|D_t|^2}{2\tau}}e^{\tau \Phi} u.
\ena

The following is an analogue of the Definition~\ref{defconvexfunctions}, under which the Carleman estimate of Theorem~\ref{thmCarlclass} holds. Here, the condition is weaker for it is only restricted to $\xi_t=0$.

\begin{definition}[Pseudoconvex function in $\xi_t=0$]
\label{defpseudofunctionwave}
With the above assumptions for $P$, let $\Phi$ be smooth and real-valued. We say that $\Phi$ is a pseudoconvex function with respect to $P$ in $\xi_t=0$ at $\x_0$ if 
\bnan
\label{hypopseudoconvecCarletau0wave}
 \left\{ p_2, \{p_2 ,\Phi \} \right\} (\x_0,\xi) >0  , 
&\text{ if } p_2(\x_0,\xi) = 0, \quad  \xi_t=0, \quad \xi \neq 0 ; \\
\label{hypopseudoconvecCarlewave}
\frac{1}{i\tau}\{ \overline{p}_\Phi,p_\Phi \} (\x_0 ,\xi,\tau) >0  , 
&\text{ if } p_\Phi(\x_0,\xi, \tau) = 0,  \quad  \xi_t = 0 , \quad \tau >0 ,
\enan
where $ p_\Phi(z,\xi,\tau) = p_2(z, \xi + i \tau d\Phi(z))$.
\end{definition}
\begin{theorem}[Carleman estimate for wave type operators with coefficients constant in time]
\label{th:carlemanwave}
Let $P$ satisfy the assumptions of Theorem~\ref{thmwave}. Let $\Phi$ be a {\em quadratic real-valued polynomial} such that $\Phi$ is a pseudoconvex function with respect to $P$ in $\xi_t=0$ at $\x_0$, in the sense of Definition \ref{defpseudofunctionwave}.
Then, there exist $r, \eps , \mathsf{d}, C, \tau_0>0$ such that for all $\tau \geq \tau_0$ and $w \in H^1_{\comp}(B(\x_0,r))$, we have 
\bnan
\label{Carlemanwave}
 \tau \|Q_{\e,\tau}^{\Phi}w\|_{\Ht{1}}^2 \leq C
 \nor{Q_{\e,\tau}^{\Phi}P w}{L^2}^2+ Ce^{-\mathsf{d}\tau}\nor{e^{\tau\Phi}w}{\Ht{1}}^2 .
\enan
\end{theorem}
Note that for $\e=0$, this would be a classical Carleman estimate. 
Yet, the role of the Fourier multiplier $e^{-\e\frac{|D_t|^2}{2\tau}}$ is to truncate the ``high frequencies'' (with respect to $\tau$) in the variable $t$. So, we just need to look at small frequencies in $\xi_t$ (compared to $\tau$, namely $|\xi_t| \leq \varsigma \tau$). Note that the set $|\xi_t| \geq \varsigma \tau$ only contributes to $e^{-\eps \varsigma^2 \tau/2}$ to the estimate, which is an admissible remainder in view of~\eqref{Carlemanwave}.
 This is why the pseudoconvexity assumption is only made in $\xi_t=0$. 

\subsection{Proving unique continuation using the Carleman estimate}
In this section, we assume that Theorem \ref{th:carlemanwave} is proved and we prove Theorem \ref{thmwave}. Part of the argument is similar to the classical case: constructing an appropriate pseudoconvex function $\Phi$  in $\xi_t=0$ from the function $\Psi$ defining the hypersurface $S=\{\Psi=  \Psi(\x_0)\}$. 
The main differences are the following:
\begin{itemize}
\item the pseudoconvexity is only on $\xi_t=0$, so it requires a small adaptation of the convexification procedure. Moreover, we want $\Phi$ to be quadratic.
\item the Carleman estimate implies an additional Fourier multiplier ($e^{-\e\frac{|D_t|^2}{2\tau}}$) that changes the proof of unique continuation. The additional difficulty comes from the fact that the Carleman estimate~\eqref{Carlemanwave} only dominates the low frequencies in $\xi_t$ of the function $u$.
\end{itemize}

\subsubsection{Convexification}
Similarly to the classical case studied in Section~\ref{chapterclassical}, the natural assumption for the unique continuation Theorem~\ref{thmwave} is a strong pseudoconvexity condition similar to that of Definition~\ref{def-pseudosurface}, but restricted to the set $\{\xi_t=0\}$. We define this notion, and then check that any noncharacteristic hypersurface is strongly pseudoconvex  in $\xi_t=0$.
\begin{definition}[Pseudoconvex hypersurface in $\xi_t=0$]
\label{def-pseudosurface-xit=0}
Let $\Omega \ni \x_0$ be an open set, $P\in \difclas{2}(\Omega)$ with real-valued principal symbol $p_2$ and $\Psi\in C^{\infty}(\Omega)$ real-valued. 
We say that the {\em oriented hypersurface} $S = \{\Psi = \Psi(\x_0)\} \ni \x_0$ is strongly pseudoconvex with respect to $P$ at $\x_0$ in $\xi_t=0$ if
 \bnan
\label{e:pseudo-surfacewave-1}
\left\{ p_2, \{p_2 ,\Psi \} \right\} (\x_0 ,\xi) >0  , 
&\text{ if } p_2(\x_0,\xi) =\{p_2,\Psi\}(\x_0,\xi)= \xi_t = 0 \text{ and } \xi \neq 0 ; \\
\label{e:pseudo-surfacewave-2}
\frac{1}{i\tau}\{ \overline{p_\Psi},p_\Psi \} (\x_0 ,\xi, \tau) >0  , 
&\text{ if } p_\Psi(\x_0,\xi, \tau) =\{p_\Psi,\Psi\}(\x_0,\xi, \tau) = \xi_t = 0 \text{ and } \tau >0 ,
\enan
where $p_\Psi(z,\xi, \tau) = p_2(z, \xi + i \tau d\Psi(z))$.
\end{definition}

The next lemma explains that the noncharacteristicity condition assumed in Theorem~\ref{thmwave} is a particular case of Definition~\ref{def-pseudosurface-xit=0}.

\begin{lemma}[Noncharacteristicity implies strong pseudoconvexity in $\xi_t=0$]
\label{proppseudosurfacewave}
Let $\Omega, P$ as in Theorem~\ref{thmwave}. If the hypersurface $S = \{\Psi = \Psi(\x_0)\} \ni \x_0$ is noncharacteristic for $P$ at $\x_0$ ($p_2(\x_0,d\Psi(\x_0)) \neq 0$), then it is strongly pseudoconvex with respect to $P$ at $\x_0$ in $\xi_t=0$.
\end{lemma}

\bnp
The principal symbol of $P$ is $p_2(t,x,\xi_t,\xi_x)=-\xi_t^2+q(x, \xi_x)$ where $q(x,\xi_x)=\sum_{i,j}a^{ij}(x)\xi_{i}\xi_{j}$. 
We first notice that for $\xi_t=0$, we have $p_2(t,x,0,\xi_x)=q(x,\xi_x)$. Since $q$ is assumed to be elliptic, the assumption $p(\x_0,\xi) =\xi_t= 0$ implies $\xi=0$ and therefore Condition~\eqref{e:pseudo-surfacewave-1} is empty.
Second, we have proved in Lemma~\ref{l:non-caract-case} that \eqref{e:pseudo-surfacewave-2} is empty if $p_2(\x_0,d\Psi(\x_0)) \neq 0$, which concludes the proof.
\enp

Next, we follow the same convexification procedure as Section \ref{Convexification}.
\begin{proposition}[Analytic convexification]
\label{propanalyticconvexwave}
Let $\Omega, P$ satisfy the assumptions of Theorem \ref{thmwave}. 
Assume that the hypersurface $S= \{\Psi= \Psi(\x_0)\}$ is strongly pseudoconvex with respect to $P$ at $\x_0$ in $\xi_t=0$, in the sense of Definition~\ref{def-pseudosurface-xit=0}.
Then there exists $\lambda_0>0$ such that for all $\lambda \geq \lambda_0$, the function $\Phi=e^{\lambda \Psi}$ is a pseudoconvex function with respect to $P$ at $\x_0$ in $\xi_t=0$, in the sense of Definition \ref{defpseudofunctionwave}.
\end{proposition}

Note however that, as opposed to the classical case, the Carleman estimate of Theorem~\ref{th:carlemanwave} does not apply to the weight function $\Phi$ since it is not (yet) quadratic.

\bnp
The proof is very similar to that of Proposition \ref{propanalyticconvex}. Again, we assume that $\Psi(\x_0)=0$ for simplicity, and use the notation  $c_{\Psi}(\xi,\tau),c_{\Phi}(\xi,\tau)$ in~\eqref{e:def-c-psi}.
  Lemma~\ref{lmlien} still applies and $c_{\Psi}(\xi,\tau)$ and $c_{\Phi}(\xi,\tau)$ are both continuous on the whole $\R^n \times \R^+$.
Then, using Lemma \ref{l:fgh-three-fcts}, Definition \ref{def-pseudosurface-xit=0} may be equivalently reformulated as the existence of constants $C_1$, $C_2>0$ so that
\bna
c_{\Psi}(\xi,\tau)+C_1\left[\left|\{p_{\Psi},\Psi\}(\x_0,\xi,\tau)\right|^2  + \frac{\left|p_{\Psi}(\x_0,\xi,\tau)\right|^2}{|\xi|^2+\tau^2}+|\xi_t|^2\right]\geq C_2(|\xi|^2+\tau^2).
\ena
Lemma~\ref{l:calcul-cphi-cpsi} still applies, and the argument in the Proof of Proposition \ref{propanalyticconvex} then implies 
\bna
c_{\Phi}(\xi,\tau)+\tilde{C}_1\left[\frac{\left|p_{\Phi}(\x_0,\xi,\tau)\right|^2}{|\xi|^2+\tau^2}+|\xi_t|^2\right]\geq \tilde{C}_2(|\xi|^2+\tau^2),
\ena
for $\lambda$ large enough. This implies the sought result.
\enp
It remains to perform the Geometric convexification and to ensure that we can take the weight function $\Phi$ quadratic.
 \begin{proposition}[Geometric convexification]
\label{propgeometricconvexwave}
Let $\Phi$ be a pseudoconvex function for $P$ at $\x_0$ in $\xi_t=0$, in the sense of Definition \ref{defpseudofunctionwave} with $\Phi(\x_0)=0$.
Then there exists a function $\varphi$ such that 
\begin{enumerate}
\item $\varphi$ pseudoconvex function for $P$ at $\x_0$ in $\xi_t=0$,
\item $\varphi$ is a quadratic polynomial,
\item \label{ecartwavelm} $\varphi(\x_0)=0$ and there exists $R_0>0$ such that for any $0<R<R_0$, there exists $\eta>0$ so that $\varphi(\x)\leq -\eta$ for $\x\in \{\Phi\leq 0\}\cap \{R/2\leq |\x-\x_0|\leq R \}$.
\end{enumerate}
\end{proposition}
\bnp
For $\delta>0$, we take $\varphi(\x)=\Phi_T(\x)-\delta |\x-\x_0|^2$, where 
\bna
\Phi_T(\x)=\sum_{|\alpha|\leq 2 }\frac{1}{\alpha !}(\partial^{\alpha}\Phi)(\x_0)(\x-\x_0)^{\alpha}
\ena
is the Taylor expansion of $\Phi$ at order $2$.
First, we notice that the pseudoconvexity condition only involves derivative up to order $2$ at $\x_0$. Hence, $\Phi_T$ is also a strongly pseudoconvex function in $\xi_t=0$ at $\x_0$. Moreover, the same stability argument as in Proposition~\ref{propgeometricconvex} applies. So, for $\delta$ small enough, $\varphi$ is as well a strongly pseudoconvex function in $\xi_t=0$ at $\x_0$. We fix $\delta>0$ sufficiently small. It remains to prove the geometric properties of Item~\ref{ecartwavelm}.
Since $\Phi_T$ is the Taylor expansion of $\Phi$ at order $2$, there exists $R_0$ small enough so that $|\Phi_T-\Phi|\leq |\x-\x_0|^2\delta/2$ for $|\x-\x_0|\leq R_0$. Now, take $R\leq R_0$.
Let $\x\in \{\Phi\leq 0\}\cap \{R/2\leq |\x-\x_0|\leq R \}$. Since $\Phi(\x)\leq 0$, we have $\Phi_T(\x)\leq |\x-\x_0|^2\delta/2$. Therefore,
$\varphi(\x)\leq -\delta |\x-\x_0|^2/2$.
So, in particular since $|\x-\x_0|^2\geq R^2/4$, we get $\varphi(\x)\leq -\delta R^2/8$ and we can take $\eta=\delta R^2/8$.
\enp

\subsubsection{Unique continuation}
In this section, we conclude the proof of the unique continuation Theorem \ref{thmwave} assuming the Carleman estimate of Theorem~\ref{th:carlemanwave}.
\bnp[Proof of Theorem \ref{thmwave}] Let $u$ solution of $Pu=0$ in $\Omega$ so that $u=0$ on $\Omega\cap \left\{\Psi> 0\right\}$. The hypersurface $S= \{\Psi= \Psi(\x_0)\}$ is strongly pseudoconvex at $\x_0$ in $\xi_t=0$. Propositions \ref{propanalyticconvexwave} and \ref{propgeometricconvexwave} allow to produce a quadratic function $\Phi$ (it is the function called $\varphi$ in Proposition~\ref{propgeometricconvexwave}, which we now rename as $\Phi$) that satisfies the pseudoconvexity for functions at $\x_0$ in $\xi_t=0$. In particular, Theorem~\ref{th:carlemanwave} applies. We therefore obtain:
\begin{enumerate}
\item \label{enumCarlwave}
 there exist $R, \eps , \mathsf{d}, C, \tau_0>0$ such that for all $\tau \geq \tau_0$ and $w \in H^1_{\comp}(B(\x_0,r))$, the Carleman estimate~\eqref{Carlemanwave} holds,
\item \label{ecartwave} $\Phi(\x_0)=0$ and there exists $\eta>0$ so that $\Phi(\x)\leq -\eta$ for $\x\in \{\Psi\leq 0\}\cap \{|\x-\x_0|\geq R/2\}$,
\item \label{majore} $\Phi(\x)\leq \mathsf{d}/4$ in $B(\x_0,R)$.
\end{enumerate}
We only added Item~\ref{majore}, which follows from a continuity statement (holding up to reducing $R$) using $\Phi(\x_0)=0$.
Now we pick $\chi\in C^{\infty}_c(B(\x_0,R))$ so that $\chi=1$ on $B(\x_0,R/2)$. We apply the Carleman estimate~\eqref{Carlemanwave} to $w=\chi u \in H^1_{\comp}(B(\x_0,R))$, solution of $Pw=\chi Pu+[P,\chi]u=[P,\chi]u$. 
Again, $[P,\chi]$ is a classical differential operator of order $1$ with coefficients supported in the set $\{\frac{R}{2}\leq|\x-\x_0|\leq R\}$. Moreover, we have $\supp (u) \subset \{\Psi\leq 0\}$, and thus $[P,\chi]u$ is supported in $\{\Psi\leq 0\} \cap \{\frac{R}{2}\leq|\x-\x_0|\}$, where $\Phi \leq -\eta$. In particular, we have $\nor{Q_{\e,\tau}^{\Phi}P w}{L^2}\leq \nor{e^{\tau\Phi} Pw}{L^2}\leq C e^{-\tau\eta}\nor{u}{H^1}$.
As for the second term in the right hand side of~\eqref{Carlemanwave}, we use Item~\ref{majore} to deduce
\bna
e^{-\mathsf{d}\tau}\nor{e^{\tau\Phi}w}{\Ht{1}}^2\leq e^{-\mathsf{d}\tau}e^{\mathsf{d}\tau/2}\nor{w}{\Ht{1}}^2\leq e^{-\mathsf{d}\tau/2}\tau^2\nor{w}{H^1}^2\leq e^{-\mathsf{d}\tau/4}\nor{u}{H^1}^2,
\ena
for $\tau$ large enough.
From~\eqref{Carlemanwave}, we have obtained that there exist $C, \delta, \tau_0>0$ so that
\bnan
\label{e:Qwleqexpo}
\|Q_{\e,\tau}^{\Phi}w\|_{L^2}\leq C e^{-\delta \tau} , \quad \text{ for all }\tau\geq \tau_0. 
\enan
We now use the following lemma, which is an analogue in the present setting to Lemma~\ref{l:uniqueness-from-carleman}, used in the classical case.
\begin{lemma}
\label{lmsupportQe}
Let $\Phi\in C^{\infty}(\Omega)$ be a real-valued function such that $d\Phi \neq 0$ on $\Omega$. Let $v\in L^2_{\comp}(\Omega)$ and assume there exists $C_0,\tau_0,\eps>0$ such that 
\bnan
\label{e:asspt-lemma-Qphi}
\|Q_{\e,\tau}^{\Phi}v\|_{L^2}\leq C_0 \quad \text{ for all } \tau \geq \tau_0.
\enan
Then, $v$ is supported in $\left\{\Phi\leq 0\right\}$.
\end{lemma}
To apply the lemma, we rewrite~\eqref{e:Qwleqexpo} as $\|Q_{\e,\tau}^{\Phi}e^{\delta \tau}w\|_{L^2}\leq C$ i.e. $\|Q_{\e,\tau}^{\Phi+\delta}w\|_{L^2}\leq C$.
Lemma~\ref{lmsupportQe} applied to the function $\Phi+\delta$ implies that $w$ is supported in the set $\left\{\Phi+\delta\leq 0\right\} = \left\{\Phi\leq -\delta\right\}$. Since we have $\Phi(\x_0)=0$ and $\chi=1$ on $B(\x_0, R/2)$, the set $V = B(\x_0, R/2) \cap \left\{\Phi > -\delta\right\}$ is a neighborhood of $\x_0$ on which  $u= w=0$, concluding the proof of the theorem.
\enp
For the proof to be complete, it remains to prove Lemma \ref{lmsupportQe}. Note that if we had $\e=0$, this is precisely Lemma~\ref{l:uniqueness-from-carleman} (and the proof is straighforward). Before describing the details of the proof, we first give a sketch of it to present the main new ideas, recalling that $n=1+d$:
\begin{enumerate}
\item Proving that $\supp(v) \subset \{\Phi \leq 0\}$ is equivalent to proving that $\x\mapsto \chi\circ \Phi(\x)v(\x)$ vanishes a.e. on $\R^{n}$ for all test function $\chi \in C^\infty(\R)$, such that $\supp(\chi)\subset [0,+\infty)$. Again, this may be reformulated equivalently in a weak form (still for all $\chi \in C^\infty(\R)$ such that $\supp(\chi)\subset [0,+\infty)$) as
$$
 \int_{\R^{n}} f(\x)v(\x) \chi(\Phi(\x)) d\x = 0 , \quad \text{ for all } f \in \calS(\R^{n}) .
$$ 
\item We change slightly the point of view and, considering $f$ fixed, see this quantity as a distribution on $\R$, with $\chi$ as test function:
\begin{equation}
\label{e:def-hf}
\langle h_f , \chi \rangle_{\mathcal{E}'(\R), C^\infty(\R)}=\langle fv , \chi (\Phi) \rangle_{\mathcal{E}'(\R^{n}), C^\infty(\R^{n})} = \int_{\R^{n}} f(\x)v(\x) \chi(\Phi(\x)) d\x .
\end{equation}
This corresponds to defining the distribution $h_f=\Phi_*(fv)$. Heuristically, $h_f(s)$ is the integral of $fv$ on the level set $\left\{\Phi(x)=s\right\}$. According to the first point, $\supp(v) \subset \{\Phi \leq 0\}$ is now equivalent to $\supp (h_f) \subset (-\infty, 0]$.

\item We shall see that the Fourier transform of $h_f$ is 
$$
\widehat{h_f}(\zeta) = \langle h_f , s \mapsto e^{-i \zeta s} \rangle_{\mathcal{E}'(\R), C^\infty(\R)}  
= \int_{\R^{n}} f(\x)v(\x) e^{-i\zeta \Phi(\x)} d\x ,
$$ and can be extended to the complex domain if $v$ is compactly supported (which is assumed here). In particular, for $\zeta\in i\R^+$, $\zeta=i\tau$, we have $\widehat{h_f}(i \tau)= \langle f , ve^{\tau \Phi} \rangle$. The assumption~\eqref{e:asspt-lemma-Qphi} gives an information on the norm of $e^{\tau \Phi}v$ for $\tau$ large which can be translated in a uniform bound on $|\widehat{h_f}|$ on the upper imaginary axis. A Phragm\'en-Lindel\"of type argument allows to transfer this uniform bound on $|\widehat{h_f}|$ to the {\em whole upper half plan}. 
\item From the bound $|\widehat{h_f}|\leq C$ on the whole upper half plan, a  Paley-Wiener theorem (roughly saying $\supp(g) \subset (-\infty,0] \Longleftrightarrow |\hat{g}|\leq C$ uniformly on the upper half complex plane) allows to conclude that $\supp(h_f) \subset (-\infty,0]$ for all $f$, which is the sought result according to the first two points.
\end{enumerate}

Let us now proceed to the details of the proof. 
\bnp[Proof of Lemma \ref{lmsupportQe}]
Let $f\in\mathcal{S}(\R^{n})$ with Fourier transform $\widehat{f}$ compactly supported in $B(0,R)$ for $R$ large. We define the distribution $h_f\in \mathcal{E}'(\R)$ by~\eqref{e:def-hf}.
Note that $h_f$ is a distribution of order zero since 
$$
|\langle h_f , \chi\rangle_{\mathcal{E}'(\R), C^\infty(\R)}| \leq  \int_{\R^{n}} | f(\x)v(\x)| | \chi(\Phi(\x))| d\x  \leq \nor{f}{L^2}\nor{v}{L^2} \sup_{\Phi(\supp(v))}|\chi| , 
$$
and is indeed compactly supported because $\supp(h_f)\subset \Phi(\supp(v))= \left\{\Phi(\x);\x\in \supp(v)\right\}$ which is compact. 
Since $h_f\in \mathcal{E}'(\R)$, the Fourier transform of $h_f$ can be computed for $\zeta \in \R$ by
\bna
\widehat{h_f}(\zeta)=\left\langle h_f,e^{-is\zeta }\right\rangle_{\mathcal{E}'(\R_s), C^\infty(\R_s)}
=\langle fv , e^{-i\zeta \Phi} \rangle_{\mathcal{E}'(\R^{n}), C^\infty(\R^{n})}
= \int_{\R^{n}} f(\x)v(\x) e^{-i\zeta \Phi(\x)} d\x .
\ena
We notice that this formula still defines a function for $\zeta \in \C$ satisfying the bound
\bnan
\label{e:upper-apriori-expo}
|\widehat{h_f}(\zeta)| \leq \int_{\supp(v)} |f(\x)v(\x)| e^{\Im(\zeta) \Phi(\x)} d\x \leq e^{C_1 |\Im(\zeta)| } \nor{f}{L^2}\nor{v}{L^2} , \quad C_1 = \max_{\supp(v)}|\Phi| .
\enan
Holomorphy of the integrand with respect to $\zeta$ implies that $\widehat{h_f}(\zeta)$ is holomorphic on the whole $\C$.
For $\zeta \in\R$, the Cauchy-Schwarz inequality~\eqref{e:upper-apriori-expo} yields the general bound
\bna
|\widehat{h_f}(\zeta)|\leq \nor{f}{L^2}\nor{v}{L^2}= C_{f,v}.
\ena
Now, we use the assumption of the lemma, namely~\eqref{e:asspt-lemma-Qphi}, to obtain a bound on the upper imaginary axis. Indeed, for $\zeta=i\tau$, and $\tau\geq \tau_0$,~\eqref{e:asspt-lemma-Qphi} implies 
\begin{align*}
|\widehat{h_f}(i\tau)|&=\left|\langle fv , e^{\tau\Phi} \rangle_{\mathcal{E}'(\R^{n}), C^\infty(\R^{n})}\right| 
= \left|\langle f , ve^{\tau\Phi} \rangle_{\mathcal{S}'(\R^{n}), \mathcal{S}(\R^{n})}\right|\\
&= \left|\langle e^{\e\frac{|D_t|^2}{2\tau}}f , e^{-\e\frac{|D_t|^2}{2\tau}}ve^{\tau\Phi} \rangle_{\mathcal{S}'(\R^{n}), \mathcal{S}(\R^{n})}\right|\\
&\leq \nor{e^{\e\frac{|D_t|^2}{2\tau}}f}{L^2(\R^{n})}\nor{e^{-\e\frac{|D_t|^2}{2\tau}}ve^{\tau\Phi}}{L^2(\R^{n})}\\
&\leq \nor{e^{\e\frac{|\xi_t|^2}{2\tau}}}{L^{\infty}(\supp(\widehat{f}))}\nor{f}{L^2(\R^{n})}\nor{Q_{\e,\tau}^\Phi v}{L^2(\R^{n})} \\
& \leq C e^{\frac{\e R^2}{2\tau }}\nor{f}{L^2(\R^{n})} \leq  C e^{\frac{\e R^2}{2\tau_0 }}\nor{f}{L^2(\R^{n})}
 =C_{\eps, f,\tau_0}C_0.
\end{align*}
Note that the operator $e^{\e\frac{|D_t|^2}{2\tau}}$ is harmless because the Fourier transform of $f$ is assumed compactly supported in $B(0,R)$ (in general, $e^{\e\frac{|D_t|^2}{2\tau}}f$ does not have any meaning, even for $f\in\mathcal{S}(\R^{n})$; this is the reason why working by duality).
Moreover, for $\tau\in[0,\tau_0]$, the estimate
$|\widehat{h_f}(i\tau)|\leq C$
follows by compactness and continuity, with some appropriate constant $C$ independent on $\tau$. 

Now, $|\widehat{h_f}|$ has a uniform bound on $\R \cup i \R_+$, as well as an a priori subexponential growth~\eqref{e:upper-apriori-expo}. We are thus in position to transfer the uniform bounds to the whole upper half plane by the Phragm\'en-Lindel\"of Theorem. 
\begin{lemma}[Phragm\'en-Lindel\"of Theorem]
\label{l:phragmen}
Let $g$ be a holomorphic function in $Q_1=\left\{x+iy;x > 0,y> 0 \right\}$, continuous in $\bar{Q}_1$.
Assume that there exist $c>0$ and $C>0$ such that
\begin{align*}
|g(z)| &\leq Ce^{c|z|}  , \quad \text{for all } z \in Q_1 , \\
|g(z)| &\leq 1 , \quad \text{for all } z \in \partial Q_1 = \R_+ \cup i \R_+ .
\end{align*}
Then, we have $|g(z)| \leq 1$ for all $z \in Q_1$.
\end{lemma}
We refer e.g. to~\cite[Theorem~3.4]{SS:03} for a proof.
Applying this result to the function $g= \widehat{h_f}$ on both $Q_1$ and the quarter plane $\left\{x+iy;x < 0,y> 0 \right\}$, we obtain that 
\bna
|\widehat{h_f}(\zeta)| \leq C \quad \text{ for all } \zeta\in\C, \Im(\zeta)\geq 0 .
\ena
We may now apply the following version of the Paley-Wiener theorem to $h_f$.
\begin{theorem}[Paley-Wiener-Schwartz]
\label{t:paley-wiener}
Suppose that $g \in \mathcal{E}'(\R)$ is of order zero. Then the following two statements are equivalent:
\begin{itemize}
\item $\supp(g) \subset (-\infty, 0]$,
\item $\widehat{g}$ can be extended continuously as an entire function which is uniformly bounded in the closed upper half-plane $$\C^+=\left\{ x + iy ; x \in \R , y \geq 0\right\} . $$ 
\end{itemize}
\end{theorem}
This is a particular case of the general Paley-Wiener-Schwartz theorem, see e.g.~\cite[Theorem~7.3.1]{Hoermander:V1}.
Applying this result to the function $h_f$  gives $\supp(h_f)\subset ]-\infty,0]$. Therefore, we have proved that for $\chi\in C^{\infty}(\R)$, 
\bna
\supp (\chi) \subset [0, + \infty) \quad  \implies \quad  0= \left\langle h_f, \chi \right\rangle = \left\langle f v,\chi(\Phi)\right\rangle=\left\langle f ,\chi(\Phi)v\right\rangle .
\ena
Since this is true for a subset of function $f$ dense in $\mathcal{S}$ (those having compactly supported Fourier transform), this means that the function $\chi(\Phi)v$ is identically zero on $\R^{n}$ as soon as $\supp (\chi) \subset [0, + \infty)$. 
That is to say $v=0$ a.e. on $\Phi>0$ or $\supp (v) \subset \{\Phi \leq 0\}$, which concludes the proof of the lemma.
\enp

\subsection{The Carleman estimate}
To complete the proof of Theorem~\ref{thmwave}, we are now left to proving the Carleman estimate of Theorem~\ref{th:carlemanwave}.
\subsubsection{The ``conjugated operator''}
As in the classical case, we first compute the ``conjugated operator''. Yet, we have to be a little careful, since $e^{\e\frac{|D_t|^2}{2\tau}}$ is not well defined on any Sobolev space and even not on $\mathcal{S}$. 
As before, we make the change of variable $v=e^{\tau \Phi}w$ and~\eqref{Carlemanwave}, rewrites 
\bna
 \tau \|e^{-\e\frac{|D_t|^2}{2\tau}} v\|_{\Ht{1}}^2  \leq  C
 \nor{e^{-\e\frac{|D_t|^2}{2\tau}}P_{\Phi} v}{L^2}^2+ Ce^{-\mathsf{d}\tau}\nor{v}{\Ht{1}}^2
\ena
The operator $P$ commutes with $e^{-\e\frac{|D_t|^2}{2\tau}}$ since its coefficients are independent on $t$. Yet, the operator $P_{\Phi}=e^{\tau\Phi}Pe^{-\tau\Phi}$ now depends on $t$ through $\Phi$. We take advantage of the fact that $\Phi$ is quadratic, hence the coefficients of the operator $P_{\Phi}$ only involves derivative of $\Phi$ of order at most $1$ and therefore are linear in $t$.
We first prove the following simple lemma.
\begin{lemma}
\label{lmcommuteps} 
Let $u\in\mathcal{S}(\R^{1+d})$, then
\bna
e^{-\e\frac{|D_t|^2}{2\tau}}(tu)=\left(t+i\e\frac{D_t}{\tau} \right)e^{-\e\frac{|D_t|^2}{2\tau}}u.
\ena
\end{lemma}
\bnp
We first recall that $\F(tv)(\xi)=i\partial_{\xi_t}\widehat{v}(\xi)$ and hence 
\begin{align*}
\F\left(e^{-\e\frac{|D_t|^2}{2\tau}}(tu)\right)(\xi)
&=e^{-\e\frac{|\xi_t|^2}{2\tau}}\widehat{(tu)}(\xi)=e^{-\e\frac{|\xi_t|^2}{2\tau}}i\partial_{\xi_t}\widehat{u}(\xi)=i\partial_{\xi_t}\left[e^{-\e\frac{|\xi_t|^2}{2\tau}}\widehat{u}(\xi)\right]+i\frac{\e\xi_t}{\tau}e^{-\e\frac{|\xi_t|^2}{2\tau}}\widehat{u}(\xi)\\
& = \F\left(te^{-\e\frac{|D_t|^2}{2\tau}}u+i \frac{\e D_t}{\tau}e^{-\e\frac{|D_t|^2}{2\tau}}u\right)(\xi) ,
\end{align*}
which proves the lemma.
\enp
\begin{remark}
Lemma \ref{lmcommuteps} can be iterated to deduce that
\bna
e^{-\e\frac{|D_t|^2}{2\tau}}(t^ku)=\left(t+i\e\frac{D_t}{\tau} \right)^k e^{-\e\frac{|D_t|^2}{2\tau}}u ,
\ena
where the exponent $k$ is meant in the sense of composition. For $f$ polynomial in $t$, we obtain
\bna
e^{-\e\frac{|D_t|^2}{2\tau}}(f(t)u)=f\left(t+i\e\frac{D_t}{\tau}\right) e^{-\e\frac{|D_t|^2}{2\tau}}u .
\ena
This means that the ``formal'' conjugated operator of $f(t)$ by $e^{-\e\frac{|D_t|^2}{2\tau}}$ is a differential operator, whose order is given by the degree of the polynomial $f$. 
For a general real-analytic function $f(t)$, giving a meaning to $f\left(t+i\e\frac{D_t}{\tau}\right)$ is one of the difficulties in~\cite{Tataru:95,RZ:98,Hor:97,Tataru:99}. 
 \end{remark} 
We now want to understand how $Q_{\e,\tau}^{\Phi}$ ``commutes'' with an operator $P$. To this aim, let us first consider the simplest case in which $P=D_j$. 
We have the following key lemma.
\begin{lemma}
\label{l:conjugation-D}
Assume $\Phi$ is a real polynomial of degree two in the variable $t$. For all $k \in \{0,\cdots ,d\}$ (with the convention $t=\x_0$, $D_0=D_t$, $\x_k=x_k$ and $D_k=D_{x_k}$ for $1\leq k \leq d$)
$$
Q_{\e,\tau}^{\Phi} D_k =(D_k)_{\Phi,\e}Q_{\e,\tau}^{\Phi} ,
$$
where (denoting $\Phi''_{t,\x_k} = \d_t \d_{\x_k}\Phi$) 
$$
(D_k)_{\Phi,\e}= D_k +i\tau \d_k\Phi(\x)-\e \Phi''_{t,\x_k}D_t .
$$
\end{lemma}
Note that since $\Phi$ is quadratic in the variable $t$, the quantity $\Phi''_{t,\x_j}$ is actually constant in $t$! In particular, the principal symbol of $(D_k)_{\Phi,\e}$ is $\xi_k +i\tau \d_k\Phi-\e \Phi''_{t,\x_k}\xi_t$.
\bnp
Since $\Phi$ is quadratic in the variable $t$, $\partial_k\Phi$ is a polynomial of degree $1$ in $t$ and can be written as $\partial_k \Phi=f_1(x)+t f_0$,
where $f_1(x)$ (resp. $f_0$) is polynomial in $x$ of order $1$ (resp. a constant). In particular, Lemma \ref{lmcommuteps} gives 
\begin{align*}
 e^{-\e\frac{|D_t|^2}{2\tau}}\left[(D_k +i\tau\partial_k \Phi)u\right]& =e^{-\e\frac{|D_t|^2}{2\tau}}\left[(D_k +i\tau(f_1(x)+tf_0))u\right]\\
&=\left[D_k +i\tau\left(f_1(x)+\left(t+i\e\frac{D_t}{\tau} \right)f_0\right)\right]e^{-\e\frac{|D_t|^2}{2\tau}}u\\
&= (D_k +i\tau\partial_k \Phi-\e  f_0D_t)e^{-\e\frac{|D_t|^2}{2\tau}}u.
\end{align*}
To get an intrinsic expression, we notice that $f_0=\partial_t \partial_k\Phi$, so $f_0D_t$ can be written $\partial_t\partial_k\Phi D_t $. This concludes the proof of the lemma.
\enp

This lemma allows to compute the principal symbol of the  ``conjugated operator'' of general differential operators (with coefficients independent of $t$).
\begin{corollary}[The ``conjugated operator'']
\label{lmsympphiwave}
Let $\Omega \subset \R^{1+d}=\R_t\times \R_x$ and $P \in \difclas{m}(\Omega)$ be a (classical) differential operator with principal symbol $p_m$. Assume also that all its coefficients are independent on $t$ (that is $p_{\alpha}(\x)=p_{\alpha}(x)$ for all $|\alpha|\leq m$).
Let $\Phi$ be a real-valued quadratic function.
Then, for any $\e>0$, there exists a unique $P_{\Phi,\e}\in \dif{m}(\Omega)$ so that  
\bna
Q_{\e,\tau}^{\Phi}P =P_{\Phi,\e}Q_{\e,\tau}^{\Phi} .
\ena
Moreover, the principal symbol of $P_{\Phi,\e}$ is
\bna
p_{\Phi,\e}(\x,\xi,\tau)=p_m\big(\x,\xi+i\tau d\Phi(\x)-\e \Phi''_{t,\x}\xi_t\big), \ena
where we use the notation $\Phi''_{t,\x}\xi_t=\Hess(\Phi)((\xi_t,0,\cdots,0) , \cdot)=\xi_t V$ with $V$ the constant vector with coefficients $V_k=(\partial_t \partial_{\x_k} \Phi)$ (using the convention of Lemma~\ref{l:conjugation-D}).
\end{corollary}
We stress the fact that all coefficients of $P$ should be independent of $t$: this is not an assumption on the principal part of the operator only. The proof is similar to that of Lemma \ref{lmsympphi}, using Lemma~\ref{l:conjugation-D}, Item~\ref{i:composition} of Proposition~\ref{p:calcul} together with the fact that the coefficients of $P$ commute with $Q_{\e,\tau}^{\Phi}$.

\begin{remark}
\label{rem:tau-divise-im}
In the case of a second order operator $P$ (with coefficients independent of $t$), with real symbol $p_2$, we have (denoting by $\tilde{p}_2$ the polar symmetric bilinear form of $p_2$), 
\begin{align*}
p_{\Phi,\e}(\x,\xi, \tau)&=p_2(\x,\xi+i\tau d\Phi(\x)-\e \Phi''_{t,\x}\xi_t)\\ 
&=p_2(\x,\xi-\e \Phi''_{t,\x}\xi_t) - \tau^2  p_2(\x,d\Phi(\x)) +2i\tau \tilde{p}_2 (\x,\xi-\e \Phi''_{t,\x}\xi_t ,d\Phi(\x)) .
\end{align*}
As in the classical case, an important point here is that $\Im(p_{\Phi,\e}(\x,\xi, \tau)) = \tau 2\tilde{p}_2 (\x,\xi-\e \Phi''_{t,\x}\xi_t ,d\Phi(\x))$ and may be divided by $\tau$. 
\end{remark}
An important feature of the Corollary~\ref{lmsympphiwave} is that the principal symbol of $P_{\Phi,\e}$ is actually close to the principal symbol of $P_{\Phi}$ if $\e$ is small. So, we can expect that it satisfies the same subelliptic estimates. 

\subsubsection{A first subelliptic estimate}

We first write the following Lemma on $p_{\Phi}$, that we have actually already used and proved in Proposition~\ref{propanalyticconvexwave}, using Lemma~\ref{l:fgh-three-fcts} and homogeneity, so we skip the proof.
\begin{lemma}
\label{lminegpsiwave}
Let $\Omega, P$ satisfy the assumptions of Theorem \ref{thmwave}. 
Assume that the function $\Phi$ is pseudoconvex with respect to $P$ at $\x_0$ in $\xi_t=0$, in the sense of Definition~\ref{defpseudofunctionwave}.
Then there exist $C_1,C_2>0$ such that for any $(\xi, \tau) \in \R^n \times \R^+$, we have
\bna
\frac{1}{i\tau}\{\overline{p_{\Phi}},p_{\Phi}\}(\x_0,\xi,\tau) +C_1\left[\frac{\left|p_{\Phi}(\x_0,\xi,\tau)\right|^2}{|\xi|^2+\tau^2}+|\xi_t|^2\right]\geq C_2(|\xi|^2+\tau^2).
\ena
where we have extended $\frac{1}{i\tau}\{\overline{p_{\Phi}},p_{\Phi}\}(\x_0,\xi,\tau)$ by continuity at $\tau=0$ with the value $2\{p,\{p,\Phi\}\}(\x_0,\xi)$.
\end{lemma}
By perturbation, we obtain a similar conclusion for the perturbated operator.
\begin{lemma}
\label{lmpsudoconvexinege}
Let $\Omega, P$ satisfy the assumptions of Theorem \ref{thmwave}. 
Assume that the function $\Phi$ is pseudoconvex with respect to $P$ at $\x_0$ in $\xi_t=0$, in the sense of Definition~\ref{defpseudofunctionwave}.
Then there exists $\e_0>0$ so that for any $0\leq \e<\e_0$, there exist $C_1, C_2>0$ such that for any $(\xi, \tau) \in \R^n \times \R^+$, we have
\bna
\frac{1}{i\tau}\{\overline{p_{\Phi,\e}},p_{\Phi,\e}\}(\x_0,\xi,\tau)+C_1\left[\frac{\left|p_{\Phi,\e}(\x_0,\xi,\tau)\right|^2}{|\xi|^2+\tau^2}+|\xi_t|^2\right]\geq C_2(|\xi|^2+\tau^2).
\ena
where the quantity $\frac{1}{i\tau}\{\overline{p_{\Phi,\e}},p_{\Phi,\e}\}(\x_0,\xi,\tau)$ is extended by continuity at $\tau=0$.
\end{lemma}
\bnp
The lemma mainly follows from the fact that $p_{\Phi,\e}$ is a perturbation of $p_{\Phi}$ and using Lemma \ref{lminegpsiwave}. Yet, we have to be a little careful because of the factor $\frac{1}{\tau}$. Noticing again that $\frac{1}{i\tau}\{\overline{p_{\Phi,\e}},p_{\Phi,\e}\}=\frac{2}{\tau}\{\Re p_{\Phi,\e},\Im p_{\Phi,\e}\}$ and using Remark~\ref{rem:tau-divise-im}  
we can write $\Im p_{\Phi,\e}=\tau \widetilde{p_{\Phi,\e}}^i$. Moreover, $\widetilde{p_{\Phi,\e}}^i$ and all its derivatives are all continuous in $\e$. Hence, we can write $\frac{1}{i\tau}\{\overline{p_{\Phi,\e}},p_{\Phi,\e}\}=2\{\Re p_{\Phi,\e},\widetilde{p_{\Phi,\e}}^i\}$, which may therefore be extended by continuity to $\tau=0$. The result then follows by a perturbation of Lemma~\ref{lminegpsiwave}.
\enp

We are now ready to prove a first subelliptic estimate that will be crucial for the final proof of Theorem~\ref{th:carlemanwave}.
\begin{proposition}
\label{p:subellipti-xit=0}
Let $\Omega, P$ satisfy the assumptions of Theorem \ref{thmwave}. 
Assume that the function $\Phi$ is real-valued, quadratic and strongly pseudoconvex with respect to $P$ at $\x_0$ in $\xi_t=0$, in the sense of Definition~\ref{defpseudofunctionwave}.
Then, there exist $\e,r,C,\tau_0>0$ so that we have the estimate
\bnan
\label{Carlpropwaveinterm}
\tau \nor{v}{\Ht{1}}^2\leq C \nor{P_{\Phi,\e}v}{L^2}^2+C\tau \nor{D_t v}{L^2}^2,
\enan
for any $v\in C^{\infty}_c(B(\x_0,r))$ and $\tau\geq \tau_0$.
\end{proposition}
Note that the parameter $\eps>0$ is fixed by this proposition (in fact, by Lemma~\ref{lmpsudoconvexinege}). This estimate is very close to the usual Carleman estimate~\eqref{Carlthm} of Theorem~\ref{thmCarlclass}. The only difference is the last term $\tau \nor{D_t v}{L^2}^2$ in the right hand-side. This term comes from the fact that the pseudoconvexity assumption (and hence the symbolic estimate of Lemma~\ref{lmpsudoconvexinege}) is made on $\xi_t=0$ only, i.e. on $D_t=0$ only.
Also, remark that this additional term has precisely the same strength as the term $\tau \nor{v}{\Ht{1}}^2$ on the left handside of the estimate.

\bnp
The proof is as well very similar to that of Theorem \ref{thmCarlclass}. A little care is needed to factorize the skew-adjoint part of the operator.
Note first that the form of Estimate~\eqref{Carlpropwaveinterm} remains unchanged under addition to $P$ of a (classical) differential operator in $\difclas{1}(\Omega)$, the coefficients of which do not depend on the variable $t$. Indeed, after conjugation, the latter perturbation will yield a perturbation of $P_{\Phi,\e}$ being in $\dif{1}(\Omega)$, which, applied to $v$, is bounded by $\nor{v}{\Ht{1}}^2$ and thus can be absorbed in the left handside for $\tau\geq \tau_0$ with $\tau_0$ large enough.

We then notice that, with $P$ satisfying the assumptions of Theorem \ref{thmwave}, we have
$$
P=\tilde{P} + R_1, \quad \text{with}\quad  \tilde{P} = -D_t^2 + \sum_{1\leq i,j\leq d}D_i g^{ij}(x)D_j, \quad \text{and}\quad  R_1\in \dif{1}(\Omega) ,
$$
where $g^{ij}(x)=g^{ji}(x)$. See also Example~\ref{e:general-operator-2}.  
According to the previous discussion, it is sufficient to prove Estimate~\eqref{Carlpropwaveinterm} for $P$ replaced by $\tilde{P}$.
Applying Lemma~\ref{l:conjugation-D} and using that $Q_{\e,\tau}^{\Phi}$ exactly commutes with $g^{ij}(x)$, we have the exact formula:
$$
\tilde{P}_{\Phi,\e} =-(D_t +i\tau\partial_t \Phi-\e  \Phi''_{t,t}D_t)^2 + \sum_{1\leq i,j\leq d}(D_i +i\tau\partial_i \Phi-\e  \Phi''_{t,x_i}D_t) g^{ij}(x)(D_j +i\tau\partial_j \Phi-\e  \Phi''_{t,x_j}D_t) .
$$
We now collect all terms being factorized by $\tau$ to obtain, for some $M \in \dif{1}(\Omega)$,
$$
\tilde{P}_{\Phi,\e} =-(D_t -\e  \Phi''_{t,t}D_t)^2 + \sum_{1\leq i,j\leq d}(D_i -\e  \Phi''_{t,x_i}D_t) g^{ij}(x)(D_j -\e  \Phi''_{t,x_j}D_t)  + \tau M, 
$$
and remark that $\tilde{P}_{\Phi,\e}- \tau M$ is a formally selfadjoint operator. As a consequence, when defining
$$
P_{R,\e} = \frac{\tilde{P}_{\Phi,\e}+\tilde{P}_{\Phi,\e}^*}{2} , \qquad P_{I,\e} = \frac{\tilde{P}_{\Phi,\e}-\tilde{P}_{\Phi,\e}^*}{2i} , 
$$
we notice that we have, as in the proof of Theorem \ref{thmCarlclass}, $P_{I,\e} = \frac{\tau M - \tau M^*}{2i} =: \tau \widetilde{P_{I,\e}}$ (that is, $\tau$ can be factorized in the skew-adjoint part of $P_{I,\e}$).
With this decomposition, we have $P_{\Phi,\e}=P_{R,\e}+i P_{I,\e}=P_{R,\e}+i\tau \widetilde{P_{I,\e}}$, and may now proceed to the key computation, following the proof of Theorem~\ref{thmCarlclass}. We obtain
$$
\nor{P_{\Phi,\e}v}{L^2}^2 
=\nor{P_{R,\e} v}{L^2}^2+\nor{P_{I,\e} v}{L^2}^2+\tau \left(i[ P_{R,\e},\widetilde{P_{I,\e}}]v,v\right) , 
$$
and the same computations as in the proof of Theorem~\ref{thmCarlclass} lead to
\begin{align*}
&\frac{1}{\tau}\nor{P_{\Phi,\e}v}{L^2}^2\geq \left(Lv,v\right) , \quad \text{ with } \\
&L= C_1 P_{R,\e}(-\Delta+\tau^2)^{-1}P_{R,\e}+C_1 P_{I,\e}(-\Delta+\tau^2)^{-1}P_{I,\e}+i[P_{R,\e},\widetilde{P_{I,\e}}] ,
\end{align*}
for $\tau \geq \tau_0$, $\tau_0$ large enough, and $C_1$ being taken as in the conclusion of Lemma~\ref{lmpsudoconvexinege}.
We thus obtain 
$$
\frac{1}{\tau}\nor{P_{\Phi,\e}v}{L^2}^2+C_1 \nor{D_t v}{L^2}^2\geq \left((L+C_1 D_t^2)v,v\right).
$$
The principal symbol of $L+C_1 D_t^2$ is 
\bna
\frac{1}{i\tau}\{\overline{p_{\Phi,\e}},p_{\Phi,\e}\}(\x,\xi,\tau)+C_1 \left[\frac{\left|p_{\Phi,\e}(\x,\xi,\tau)\right|^2}{|\xi|^2+\tau^2}+|\xi_t|^2\right].
\ena
We conclude as in the proof of Theorem~\ref{thmCarlclass} using the symbolic estimate of Lemma~\ref{lmpsudoconvexinege} at the point $\x_0$, together with the G{\aa}rding inequality of Proposition~\ref{lmweakGard}.
\enp

\subsubsection{End of the proof of the Carleman estimate}
Equipped with the subelliptic estimate~\eqref{Carlpropwaveinterm}, we may now proceed to the proof of Theorem~\ref{th:carlemanwave}.
 Setting $v = Q_{\eps,\tau}^{\Phi} w = e^{-\frac{\e}{2\tau}|D_t|^2}(e^{\tau\Phi}w)$, we need to prove the estimate
\bna
\tau\nor{v}{\Ht{1}}^2\leq C\nor{P_{\Phi,\e}v}{L^2}^2+Ce^{-\mathsf{d}\tau}\nor{e^{\tau\Phi}w}{\Ht{1}}^2.
\ena
The latter is very close to~\eqref{Carlpropwaveinterm}, except for the last term, and it is very tempting to apply~\eqref{Carlpropwaveinterm} to $v = Q_{\eps,\tau}^{\Phi} w$. The hope is then that the term $\tau \nor{D_t v}{L^2}^2$ is estimated by using that the multiplier $e^{-\e\frac{|D_t|^2}{2\tau}}$ ``localizes'' where $D_t$ is small. This will indeed be done at the end of the proof. However, the first problem we have to face is that, even if $w$ is compactly supported, the function $v = Q_{\eps,\tau}^{\Phi} w$ is not compactly supported in the variable $t$. Indeed, the operator $e^{-\frac{\e}{2\tau}|D_t|^2}$ is not local. We thus need to introduce an additional cutoff in time, and estimate the remainder it produces. 

\bnp[Proof of Theorem~\ref{th:carlemanwave}]
We assume for simplicity that the point $\x_0$ involved is $\x_0 = (0, x_0)$, i.e. $t_0=0$.
We let $\e>0$ and $r>0$ be fixed by Proposition~\ref{p:subellipti-xit=0}. We choose $r_0>0$ with $2r_0= r$, and, all along the proof, we consider functions $u\in C^\infty_c (B(\x_0,r_0/4))$.
Let $\chi \in C^\infty_c (]-r_0, r_0[)$ such that $\chi=1$ on $]-r_0/2, r_0/2[$.
Since $v$ is not compactly supported in the variable $t$, we set $f(\x) = \chi(t)v(\x)$ and we have $\supp (f) \subset [-r_0, r_0] \times B(x_0,r_0/4) \subset B(\x_0,2r_0)=B(\x_0,r)$, so that Proposition~\ref{p:subellipti-xit=0} will apply to the function $f$. To estimate $v$, we write 
\bna
\|v\|_{\Ht{1}} \leq \|f\|_{\Ht{1}} + \|v-f\|_{\Ht{1}} ,
\ena
where  
$$
v-f = (1-\chi)Q_{\e,\tau}^{\Phi} u = (1-\chi) e^{-\frac{\e}{2\tau}|D_t|^2}(\check{\chi}e^{\tau\Phi}u),
$$
for $\check{\chi} \in C^\infty_c (]-r_0/3, r_0/3[)$ with $\check{\chi} =1$ in a neighborhood of $[-r_0/4, r_0/4]$ so that $\check{\chi}u=u$.
We are in position to apply the following lemma to estimate the remainder $v-f$.

\begin{lemma}
\label{l:noyau-chaleur-hormander}
Let $\chi_1\in C^{\infty}(\R^{n})$, $\chi_2\in C^{\infty}(\R^{n})$ with all derivatives bounded such that $\dist(\supp(\chi_1),\supp(\chi_2))>0$. Then there exist $C,c>0$ such that for all $u\in\mathcal{S}(\R^n)$ and all $\lambda\geq 0$, we have 
\begin{align*}
\nor{\chi_1 e^{-\frac{|D_t|^2}{\lambda}}(\chi_2 u) }{L^2} \leq C e^{-c \lambda}\nor{u}{L^2} , \qquad 
\nor{\chi_1 e^{-\frac{|D_t|^2}{\lambda}}(\chi_2 u) }{\Ht{1}} \leq C e^{-c \lambda}\nor{u}{\Ht{1}} .
\end{align*}
\end{lemma}
As a consequence of Lemma~\ref{l:noyau-chaleur-hormander}, we obtain, for $\tau \geq \tau_0$
 \bnan
 \label{e:Carl-xit-1}
\|v\|_{\Ht{1}} \leq \|f\|_{\Ht{1}} + C e^{-c\frac{\tau}{\e}} \|e^{\tau\Phi}u\|_{\Ht{1}}.
\enan
The subelliptic estimate \eqref{Carlpropwaveinterm} applied to $f$ gives
 \bnan
 \label{e:Carl-xit-2}
\tau\nor{f}{\Ht{1}}^2\leq C \nor{P_{\Phi,\e}f}{L^2}^2+C\tau\nor{D_t f}{L^2}^2,
\enan
and we need to estimate the two terms on the right handside in terms of $v$. 
First, we estimate the term $\nor{P_{\Phi,\e}f}{L^2}=\nor{P_{\Phi,\e}\chi v}{L^2}\leq \nor{\chi P_{\Phi,\e}v}{L^2}+\nor{[P_{\Phi,\e},\chi ]v}{L^2}$. For the commutator, we write $[P_{\Phi,\e},\chi ]v=[P_{\Phi,\e},\chi ]e^{-\frac{\e}{2\tau}|D_t|^2}\check{\chi}e^{\tau\Phi}u$. We notice that $[P_{\Phi,\e},\chi ] \in \dif{1}$ with coefficients supported in $\supp(\chi'_{t})$ that is, away from $\supp (\check{\chi})$. In particular, Lemma \ref{l:noyau-chaleur-hormander} implies $\nor{[P_{\Phi,\e},\chi ]v}{L^2}\leq C e^{-c\frac{\tau}{\eps} }\nor{e^{\tau\Phi}u}{\Ht{1}}$. This yields
 \bnan
 \label{e:Carl-xit-3}
\nor{P_{\Phi,\e}f}{L^2}\leq \nor{P_{\Phi,\e}v}{L^2}+C e^{-c\frac{\tau}{\eps} }\nor{e^{\tau\Phi}u}{\Ht{1}}.
\enan
Second, we estimate the term $\| D_t f\|_{L^2}$. We obtain in a similar way
 \bnan
 \label{e:Carl-xit-4}
\| D_t f\|_{L^2}=\| D_t (\chi v)\|_{L^2} \leq \| \chi D_t v\|_{L^2}+\| \chi'(t) e^{-\frac{\e}{2\tau}|D_t|^2}\check{\chi}e^{\tau\Phi}u \|_{L^2}\leq \|  D_t v\|_{L^2} +C e^{-c\frac{\tau}{\eps} }\nor{e^{\tau\Phi}u}{L^2}
\enan
where we have used again Lemma \ref{l:noyau-chaleur-hormander} in the last inequality.

Let $\varsigma$ a small constant to be fixed later on. We distinguish between frequencies of size smaller and bigger than $\varsigma \tau$. We obtain
\begin{align*}
\|  D_t v\|_{L^2} =\|  D_t  e^{-\frac{\e}{2\tau}|D_t|^2}e^{\tau\Phi}u\|_{L^2}
& \leq \|  D_t \mathds{1}_{|D_t|\leq \varsigma \tau}v\|_{L^2}+\|  D_t \mathds{1}_{|D_t|\geq \varsigma \tau} e^{-\frac{\e}{2\tau}|D_t|^2}e^{\tau\Phi}u\|_{L^2} \\
& \leq \|  D_t \mathds{1}_{|D_t|\leq \varsigma \tau}v\|_{L^2}+\max_{\xi_t\in [\varsigma \tau, +\infty)} (\xi_t e^{-\frac{\e}{2\tau}|\xi_t|^2})\| e^{\tau\Phi}u\|_{L^2} .
\end{align*}
Now, on $\R^+$, the function $s\mapsto s e^{-\frac{\e}{2\tau}s^2}$ reaches its maximum at $s =\sqrt{\frac{\tau}{\eps}}$, and is decreasing on $[\sqrt{\frac{\tau}{\eps}}, + \infty)$.
Hence, if $\tau \geq \frac{1}{\varsigma^2 \eps}$, then $\sqrt{\frac{\tau}{\eps}} \leq \varsigma \tau$, the function $s\mapsto s e^{-\frac{\e}{2\tau}s^2}$ is decreasing on the interval $[\varsigma \tau, +\infty)$, and thus bounded by its value at $\varsigma \tau$. This yields, for all $\tau \geq \max (\tau_0, \frac{1}{\varsigma^2 \eps})$, the estimate
 \begin{align}
 \label{e:Carl-xit-5}
\|  D_t v\|_{L^2} \leq  \varsigma \tau \|  v\|_{L^2}+\varsigma \tau e^{-\frac{\tau\varsigma^2 \e}{2}}\| e^{\tau\Phi}u\|_{L^2}.
\end{align}
Combining all estimates so far, namely~\eqref{e:Carl-xit-1}-\eqref{e:Carl-xit-2}-\eqref{e:Carl-xit-3}-\eqref{e:Carl-xit-4}-\eqref{e:Carl-xit-5}, we have proved that there are some constants $c>0$ (depending on $\eps$) and $C>0$ so that for any $\varsigma>0$, we have for $\tau\geq \max (\tau_0, \frac{1}{\varsigma^2 \eps})$, 
\bna
\tau \nor{v}{\Ht{1}}^2\leq C\nor{P_{\Phi,\e}v}{L^2}^2 +C \varsigma^2 \tau^3 \|  v\|_{L^2}^2+C\left(e^{-c\tau }+\varsigma^2 \tau^3 e^{-\tau\varsigma^2 \e}\right)\nor{e^{\tau \Phi}u}{\Ht{1}}^2.
\ena
We now fix the constant $\varsigma$ small enough so that the term $C \varsigma^2 \tau^3 \|  v\|_{L^2}^2\leq C \varsigma^2 \tau \|  v\|_{\Ht{1}}^2$ can be absorbed in the left handside of the estimate. This yields the sought estimate for $\tau\geq \max (\tau_0, \frac{1}{\varsigma^2 \eps})$, and concludes the proof of the theorem.
\enp
\bnp[Proof of Lemma \ref{l:noyau-chaleur-hormander}]
Recalling the explicit expression of the Fourier transform of the Gaussian in~\eqref{e:Fourier-gaussienne}, we have 
\begin{align*}
\chi_1 e^{-\frac{|D_t|^2}{\lambda}}(\chi_2 u)(t,x)& = \left(\frac{\lambda}{4\pi}\right)^{\frac{1}{2}}\int_{\R_s}\chi_1(t,x)e^{-\frac{\lambda}{4} |s- t|^2}(\chi_2 u)(s,x)~ds\\
&= \left(\frac{\lambda}{4\pi}\right)^{\frac{1}{2}}\chi_1(t,x)\int_{s,|t-s|\geq d}e^{-\frac{\lambda}{4} |s- t|^2}(\chi_2 u)(s,x)~ds
\end{align*}
where we have used the support properties in the second equality. This yields 
\begin{align*}
|\chi_1 e^{-\frac{|D_t|^2}{\lambda}}(\chi_2 u)|(t,x)
& \leq  \nor{\chi_1}{L^{\infty}}
\left(\frac{\lambda}{4\pi}\right)^{\frac{1}{2}}\int_{s,|t-s|\geq d}e^{-\frac{\lambda}{4} |s- t|^2}|\chi_2 u|(s,x)~ds\\
& \leq   \nor{\chi_1}{L^{\infty}}
\left(\frac{\lambda}{4\pi}\right)^{\frac{1}{2}}\Big( \mathds{1}_{|\cdot|\geq d}e^{-\frac{\lambda}{4}|\cdot|^2} *_{\R_s} |\chi_2 u|(\cdot , x) \Big)(t) .
\end{align*}
As a consequence, using the Young inequality, we have 
\begin{equation}
\label{e:estim-gauss-con-young}
\|\chi_1 e^{-\frac{|D_t|^2}{\lambda}}(\chi_2 u)\|_{L^2} \leq  \nor{\chi_1}{L^{\infty}}\left(\frac{\lambda}{4\pi}\right)^{\frac{1}{2}}\nor{ \mathds{1}_{|\cdot|\geq d}e^{-\frac{\lambda}{4}|\cdot|^2}}{L^1(\R)} \|\chi_2 u\|_{L^2(\R^{n+1})} .
\end{equation}
Next, we notice that
 \begin{align*}
 \frac12 \nor{ \mathds{1}_{|\cdot|\geq d}e^{-\frac{\lambda}{4}|\cdot|^2}}{L^1(\R)} & =  \int_d^{+\infty}e^{-\frac{\lambda}{4}s^2}ds
 = \frac{2}{\sqrt{\lambda}} \int_{d\sqrt{\lambda}/2}^{+\infty}e^{-y^2}dy
 = \frac{2}{\sqrt{\lambda}} \int_{d\sqrt{\lambda}/2}^{+\infty}e^{-\frac{y^2}{2}} e^{-\frac{y^2}{2}}dy \\
& \leq  \frac{2}{\sqrt{\lambda}} e^{-\frac{d^2}{8}\lambda}\int_{d\sqrt{\lambda}/2}^{+\infty} e^{-\frac{y^2}{2}}dy \leq  \frac{2}{\sqrt{\lambda}} e^{-\frac{d^2}{8}\lambda}\int_{0}^{+\infty} e^{-\frac{y^2}{2}}dy = \frac{2}{\sqrt{\lambda}} e^{-\frac{d^2}{8}\lambda} \sqrt{\frac{\pi}{2}}.
\end{align*}
Coming back to~\eqref{e:estim-gauss-con-young}, we have obtained the existence of a constant $C>0$ such that for all $\lambda>0$, 
\bna
\|\chi_1 e^{-\frac{|D_t|^2}{\lambda}}(\chi_2 u)\|_{L^2} \leq C \nor{\chi_1}{L^{\infty}} \nor{\chi_2}{L^{\infty}} e^{-\frac{d^2}{8} \lambda} \| u\|_{L^2(\R^{n+1})} ,
\ena
 which implies the result in $L^2$. The proof  in $\Ht{1}$ is a consequence of that in $L^2$.
\enp

\subsection{Semiglobal statements  and non characteristic hypersurfaces}
\label{s:semiglobal}
In this section, we first describe a geometric setting that encompasses those described in Section~\ref{s:motivation}. We then prove a semiglobal unique continuation statement relying on the construction of a family of noncharacteristic hypersurface.
This semiglobal statement will be the cornerstone in  all applications discussed in Section~\ref{s:motivation}.

\subsubsection{Distance, metric, Laplace-Beltrami operator}
\label{s:geom-def}
We consider a connected $d$-dimensional Riemannian manifold $(\M,g)$ with or without boundary $\d \M$.  In case $\d \M \neq \emptyset$, we denote by $\Int(\M)$ the interior of $\M$, so that $\M = \d \M \sqcup \Int(\M)$ (see e.g.~\cite[Chapter~1]{Lee:book}).
The metric $g$ (a bilinear from on $T\M$) is in local charts a smooth family of symmetric elliptic matrices $(g_{ij}(x))$, i.e. $(g_{ij}(x))X^iX^j \geq c |X|^2$ for all $X \in T_x\M$.
We consider the  Riemannian volume density $d\Vol_g$, given in local charts by $d\Vol_g(x) =\sqrt{\det g(x)}dx$. 
We also define the co-metric $g^*$ to $g$, defined on $T^*\M$, given in local charts by the smooth family of symmetric elliptic (i.e. satisfying~\eqref{ellipticity}) matrices $g^*(x)=(g_{ij}(x))=(g^{ij}(x))^{-1}$.
Associated to the metric $g$ and the volume density $d\Vol_g$, the natural (negative) elliptic operator is the Laplace-Beltrami operator on $\M$, given in local charts by
$$
\Delta_g f = \sum_{i,j=1}^n \frac{1}{\sqrt{\det g}}\partial_i \left(\sqrt{\det g} \, g^{ij}\partial_j f\right) , \quad f \in C^\infty(\M)
$$
which is formally selfadjoint on $L^2(\M) := L^2(\M,d\Vol_g)$.

Given a path $\gamma \in C^1( [0,1] ;  \Int(\M))$ (or even $\gamma \in W^{1,1}( [0,1] ;  \Int(\M))$), its length (according to the Riemannian metric $g$)  is given by
\bna
\length(\gamma)=\int_0^1 |\dot{\gamma}(t)|_{\gamma(t)}dt ,
\ena
where we have written $|X|_{x}^2:=g_x(X,X)$ for $x\in\M$ and $X\in T_x\M$ (i.e. $|X|_{x}^2=g_{ij}(x)X^iX^j$ in local charts).
This allows to define the Riemannian distance associated to $g$ as
\begin{align*}
\dist(x_1,x_2)& = \inf \left\{ \length(\gamma) , \gamma \in W^{1,1}( [0,1] ;  \Int(\M)) , \gamma(0)=x_1; \gamma(1)=x_2 \right\} \\
& = \inf \left\{ \length(\gamma) , \gamma \in C^\infty( [0,1] ;  \Int(\M)) , \gamma(0)=x_1; \gamma(1)=x_2 \right\} ,
\end{align*}
where equality in the second line follows from a classical regularization argument.

Note that in case $\M \subset \R^d$,  any (locally uniformly) elliptic operator on $\M$ can be written under the form $\Delta_g$  modulo smooth lower order terms.

\subsubsection{The semiglobal theorem}
\label{s:semiglo}
The key semiglobal result in all applications discussed in Section~\ref{s:motivation}  is the following.
\begin{theorem}[Semi-global unique continuation for waves]
\label{thmwaveUCPglobal}
Let $q \in L^\infty_{\loc}(\Int(\M))$ and  consider on $\R \times \M$ the operator $P := \d_t^2 - \Delta_g + q$, with $\M,g$ satisfying the assumptions in Section~\ref{s:geom-def}. Let $x_0$, $x_1\in \Int(\M)$, let $\omega_0$ be neighborhood of $x_0$ in $\Int(\M)$.  Then, for any $T> \dist(x_0,x_1)$, there exist $\e>0$ and $V_{x_1}$ a neighborhood of $x_1$ such  that 
\bnan
\left\{\begin{array}{l}
u \in H^1_{\loc} \big((-T,T)\times \Int(\M)\big) \\ 
Pu  =  0  \textnormal{ in  } (-T,T)\times \Int(\M) \\
 u = 0 \textnormal{ in } (-T,T)\times\omega_0\end{array}\right.\Longrightarrow u= 0 \textnormal{ in } (-\e,\e)\times V_{x_1}.
\enan
\end{theorem}
The proof relies first on fixing suitable coordinates along a path joining $x_0$ and $x_1$ and having length $<T$, and second on constructing in these coordinates appropriate noncharacteristic hypersurfaces in which to apply Theorem~\ref{thmwave}. Note that the result remains true if one adds to $P$ first order differential operators with time-independent $L^\infty_{\loc}$ coefficients. 
\bnp[Proof of Theorem~\ref{thmwaveUCPglobal}]
According to the definition of $\dist$, there is a smooth injective path $\gamma : [0,1] \to \Int(\M)$ such that $\gamma(0)=x_0$, $\gamma(1)=x_1$ and $\length(\gamma) = \ell_0$ with $\dist(x_0,x_1)<\ell_0< T$.
According to Lemma~\ref{lmcoord} below, we can find local coordinates $(w, l)$ near $\gamma$ in which the path $\gamma$ by $\gamma(s)=(0,s\ell_0)$ for $s\in [0,1]$ and the cometric $g^*$ (defined on $T^*\M$) is given by the matrix $m(w,l) \in M_d(\R)$ with
\bnan
\label{fromdiagm}
m(w,l)  = \left( \begin{array}{cc}
 m' (l) & 0 \\
 0 & 1  
 \end{array}
\right) + \mathcal{O}_{M_n(\R)}(|w|) , \quad \text{ for } w\in B_{\R^{d-1}}(0,\delta), \delta >0 ,
\enan
with $m' (l)\in M_{d-1}(\R)$ (uniformly) definite symmetric. With these coordinates in the space variable, and still using the straight time variable, the symbol of the wave operator is given by
\bnan
\label{e:ptauwxn}
p_2(t,w,l, \xi_t, \xi_w, \xi_l)=p_2(w,l, \xi_t, \xi_w, \xi_l) = -\xi_t^2  +  m(w,l) \xi  \cdot \xi   ,  \quad\xi = (\xi_w, \xi_l), 
\enan
where we have used $\xi_t$ for the cotangent variable to the time variable and $\xi_w$, $\xi_l$ for the dual to $w\in B_{\R^{d-1}}(0,\delta)$ and $l\in [0,\ell_0]$ respectively.
We now aim to apply Theorem \ref{thmwave} and we need to construct appropriate non characteristic hypersurfaces.
To this aim, we let $t_0$ with $\ell_0<t_0<T$. For $b <\delta$ small, to be fixed later on, we define
 \begin{align*}
  D& =\left\{(t,w) \in [-t_0,t_0]\times B_{\R^{d-1}}(0,\delta)  \left| \Big(\frac{w}{b}\Big)^2+\Big(\frac{t}{t_0}\Big)^2\leq 1\right.\right\} , \\
 G(t,w,\e)& =  \e \ell_0\zeta \left(\sqrt{\Big(\frac{w}{b}\Big)^2+\Big(\frac{t}{t_0}\Big)^2} \right), 
 \quad  \Psi_\eps(t,w,l): =  G(t,w , \eps) - l , \quad  \e \in [0,1] ,
\end{align*}
where $\zeta$ is a fixed function such that 
 \bna
 \zeta :[-1,1]\to \R^+\text{ even,}\quad  \zeta(\pm1) =0, \quad \zeta(0)=1, \quad 
\zeta(s)\geq 0, \quad |\zeta'(s)| \leq \alpha,  \textnormal{ for }s\in [-1,1], 
 \ena
with $1<\alpha <\frac{t_0}{\ell_0}$. This is possible since $\frac{t_0}{\ell_0}>1$.
Note that the fact that $\zeta$ is even implies that $G(t,w,\e)$ is actually smooth.
Note also that the point $(t=0,w=0,l=\ell_0)$ corresponding in the local coordinates to $x_1$ belongs to the hypersurface $\left\{\Psi_{1}=0\right\}$.
 We have
  $$
d \Psi_\eps(t,w,x_n) 
=   \e \ell_0  \left(\Big(\frac{w}{b}\Big)^2+\Big(\frac{t}{t_0}\Big)^2\right)^{-1/2} \zeta' \left(\sqrt{\Big(\frac{w}{b}\Big)^2+\Big(\frac{t}{t_0}\Big)^2} \right) \left( \frac{tdt}{t_0^2}  + \frac{w dw}{b^2}  \right) - d x_n .
 $$
Given the form of the principal symbol of the wave operator in these coordinates (see \eqref{fromdiagm}-\eqref{e:ptauwxn}), we obtain
\begin{align*}
p( w,l, d \Psi_\eps(t,w,l) ) &=  -\e^2 \ell_0^2\frac{t^2}{t_0^4}    \left(\Big(\frac{w}{b}\Big)^2+\Big(\frac{t}{t_0}\Big)^2 \right)^{-1}|\zeta'|^2   \\
& \quad + \ell_0^2\frac{\e^2}{ b^4}  \langle m' (l) w, w \rangle \left(\Big(\frac{w}{b}\Big)^2+\Big(\frac{t}{t_0}\Big)^2 \right)^{-1} |\zeta'|^2  + 1 \\
& \quad +O(|w|^2) \left( 1 + \frac{\e^2 \ell_0^2}{b^4}  |w|^2 \left(\Big(\frac{w}{b}\Big)^2+\Big(\frac{t}{t_0}\Big)^2 \right)^{-1}|\zeta'|^2 \right),
\end{align*}
where $|\zeta'|^2 $ is taken at the point $\sqrt{\Big(\frac{w}{b}\Big)^2+\Big(\frac{t}{t_0}\Big)^2}$. Now, since $\alpha <\frac{t_0}{\ell_0}$ and $m' (l)$ is uniformly (for $l\in [0,\ell_0]$) definite positive, there are $\eta>0$ and $b>0$ small enough so that for $|w|\leq b$, we have
\begin{align*}
1+O(|w|^2)&\geq  \alpha^2 \frac{\ell_0^2}{t_0^2} + \eta , \\
 \langle m' (l) w, w \rangle +O(|w|^2)|w|^2&\geq  \frac{1}{2} \langle m' (l) w, w \rangle\geq 0.
\end{align*}
Hence, there is a sufficiently small neighborhood (taking again $b$ small enough) of the path (i.e. of $w = 0$), in which we have (for any $\e\in [0,1]$), and any $(t,w,l)\in {D}\times [0,\ell_0]$,
\begin{align*}
p( w,l, d \Psi_\eps(t,w,l) )& \geq   - \frac{\e^2}{t_0^2}\ell_0^2   \Big(\frac{t}{t_0}\Big)^2\left(\Big(\frac{w}{b}\Big)^2+\Big(\frac{t}{t_0}\Big)^2 \right)^{-1}  |\zeta'|^2  
+  \alpha^2 \frac{\ell_0^2}{t_0^2}+\eta  \\
& \geq   -\frac{\ell_0^2}{t_0^2}  |\zeta'|^2 +\alpha^2 \frac{\ell_0^2}{t_0^2}+\eta \geq \eta.
 \end{align*}
As a consequence, for any $\eps \in [0,1]$, the hypersurface $\{\Psi_\eps =0\}$ is noncharacteristic for $P$ near any of its points. Theorem~\ref{thmwave} thus applies 
Now, define $K_{\e} :=\left\{(t,w,l)\in {D}\times [0,\ell_0], l\leq G(t,w,\e)\right\}\cap \left\{l\geq 0\right\}$ for $\eps\in [0,1]$ and consider $\e_0:=\sup \left\{\e\in [0,1], u=0 \textnormal{ in }K_{\e}\right\}$. A continuity argument yields that that $u=0$ on $K_{\e_0}$. A compactness argument on the compact set $K_\eps$ (taking into account the ``corners'') and successive applications of Theorem~\ref{thmwave} proves that $\eps_0=1$.  Theorem~\ref{thmwave} applied once again across $\{\Psi_1=0\}$ then implies that $u=0$ in a neighborhood of $\{\Psi_1=0\}$ which contains the point $(t,w,l)=(0,0,\ell_0)=x_1$ (in these coordinates), and concludes the proof of the theorem.
\enp

The following result is proved e.g. in~\cite{LL:book}.
\begin{lemma}
\label{lmcoord}
Let $\gamma: [0,1]\rightarrow \Int(\M)$ be a smooth path without self intersection (i.e. $\gamma$ is injective) of length $\ell_0$ so that $\gamma(0)=x_0$ and $\gamma(1)=x_1$.
Then, there are coordinates $(w,l)\in B_{\R^{d-1}}(0,\e)\times [0,\ell_0]$ in an open neighborhood $U$ of $\gamma([0,1])$ such that
\begin{itemize}
\item  $\gamma([0,1])=\{w=0 \}\times [0,\ell_0]$,
\item the cometric $g^*$ (defined on $T^*\M$) is of the form $
m(w,l)  = \left( \begin{array}{cc}
 m' (l) & 0 \\
 0 &  1 
 \end{array}
\right) + O_{M_d(\R)}(|w|) ,
$ where $m'$ is a smooth family on $[0,1]$ of positive definite matrices in $M_{d-1}(\R)$.
\end{itemize}
\end{lemma}

\begin{remark}
Note that all hypersurfaces constructed in the proof of Theorem~\ref{thmwaveUCPglobal} have points where they are {\em not } stronglypseudoconvex, e.g. for $t=0$ (actually, one could prove that they are strongly pseudoconvex near none of their points).
As a consequence of Theorems~\ref{t:alinhac} and~\ref{t:alinhac-baouendi}, for each of these these hypersurfaces, one can modify the operator by a time dependent zero order term so that unique continuation does not hold.
 \end{remark}

\subsection{Global unique continuation statements: back to applications}
\label{s:GLUC}
We now come back to the motivating applications presented in Section~\ref{s:motivation}.

\subsubsection{Region of dependence}
\label{s:UC-problem}
The following result might be seen as a counterpart to finite speed of propagation for waves. recall that we set
$$
 \dist(x, E) = \inf_{y \in E} \dist(x , y) , \quad \text{ for } E \subset \M, x \in \M .
$$
\begin{theorem}
\label{t:region-dependence}
Let $(\M,g)$ be a connected Riemannian manifold with or without boundary $\d \M$, $q \in L^\infty_{\loc}(\Int(\M))$, and $\omega \subset \M$ be a nonempty open set. 
Let $u \in H^1_{\loc}((-T,T)\times \M)$ be such that $\d_t^2u-\Delta_gu + qu = 0$ in $\mathcal{D}' \big((-T,T)\times \Int(\M)\big)$. 
Assume that $u|_{(-T,T)\times \omega}=0$, then  $u|_{\mathcal{U}}=0$ where 
$$
\mathcal{U}=\mathcal{U}(\omega,T) := \big\{ (t,x) \in (-T,T)\times \M , \dist (x,\omega) < T-|t| \big\} .
$$
Similarly, if  $\Gamma \subset \d \M$ is a nonempty open set of $\d\M$, if in addition $q \in L^\infty, u \in H^1$ in a neighborhood of $\Gamma$, and if $u|_{(-T,T)\times \Gamma}=\d_n u|_{(-T,T)\times \Gamma}=0$, then  $u|_{\mathcal{U}'}=0$ where $\mathcal{U}' =\mathcal{U}(\Gamma,T)$.
\end{theorem}
The first statement in this result is a direct consequence of Theorem~\ref{thmwaveUCPglobal} together with translation invariance in time.
Concerning the second statement, it suffices to notice that for any extension $\tilde\M$ of $\M$ such that $\d\M\setminus \Gamma \subset \d\tilde{\M}\cap \d\M$, one may extend $q$ as (any) $\tilde q \in L^\infty_{\loc}(\tilde \M)$ and extend $u$ as $\tilde u$ equal to $0$ in $\tilde\M \setminus \M$. The first part of the statement then applies to $\tilde u$ and yields the second statement.

Theorem~\ref{t:region-dependence} provides the largest region $\mathcal{U}(\omega,T)$ of time-space where solutions to wave equations vanishing on the cylinder $(-T,T)\times \omega$ has to vanish. That one cannot improve the size of the region $\mathcal{U}(\omega,T)$  is a consequence of finite speed of propagation, see e.g. the constructions in~\cite{Russell:71-1,Russell:71-2}.

Theorem~\ref{t:region-dependence} is also the unique continuation statement used in the Boundary Control method to solve the hyperbolic inverse problem presented in Section~\ref{s:BC-method}. We refer the reader to~\cite{NO:22} for a presentation of the method of~\cite{Belishev:87}.

\subsubsection{Penetration into shadow and approximate controllability}
We are now prepared to discuss the motivation to ``Penetration into shadow'' and to Approximate controllability introduced in Sections~\ref{s:shadow-intro} and~\ref{s:approx-wave-intro}.
Given two subsets $E_0,E_1 \subset \M$, we introduce the largest distance of $E_0$ to a point of $E_1$ (note that this quantity is not symmetric with respect to $E_0,E_1$):
\bnan
\label{e:def-L}
\mathcal{L}(E_1 ,E_0):= \sup_{x \in E_1} \dist(x, E_0)  , \quad   \dist(x, E) = \inf_{y \in E} \dist(x , y) .
\enan
The following result is a direct corollary of Theorem~\ref{thmwaveUCPglobal} (or Theorem~\ref{t:region-dependence}) together with a compactness argument.
\begin{theorem}[Global unique continuation for waves]
\label{c:obserwaveintro}
Let $(\M,g)$ be a connected Riemannian manifold with or without boundary $\d \M$, $q \in L^\infty_{\loc}(\Int(\M))$. Let $\omega \subset \M$ be a nonempty open set and $K \subset \M$ be a compact set.
If $T> \mathcal{L}(K ,\omega)$, and if $u$ satisfies
\begin{align}
\label{e:UC-presque-final}
 u \in H^1_{\loc} ((-T,T) \times \Int(\M)), \quad (\partial_t^2 -\Delta_g +q) u= 0 \text{ in } (-T,T) \times \Int(\M) ,\quad  u = 0 \text{ on }(-T,T)\times \omega, 
\end{align}
then, there is $\eps >0$ such that $u = 0$ identically in $(-\eps, \eps)\times K$.
\end{theorem}
This result answers Question~\ref{e:UCshadow} by the affirmative if $T> \mathcal{L}(K ,\omega)$. Under this condition, all waves supported at time $0$ in $K$ are visible from $(-T,T)\times \omega$.
As a corollary (together with well-posedness for the Boundary-value problem), we deduce the following global unique continuation result.
\begin{theorem}
\label{c:global-UC}
Let $(\M,g)$ be a {\em compact} connected Riemannian manifold with or without boundary $\d \M$, $q \in L^\infty(\Int(\M))$ and  fix $T>2\mathcal{L}(\M ,\omega)$. Then~\eqref{e:global-UC} holds.
\end{theorem}
\bnp
Indeed, under the assumption of~\eqref{e:global-UC}, the solution $w$ to~\eqref{e:ondes-obs} belongs to $$C^0([0,T]; H^1_0(\M))\cap C^1([0,T];L^2(\M)) \subset H^1 ((0,T) \times \Int(\M))$$ and we can thus apply Theorem~\ref{c:obserwaveintro} to $K=\M$ and $u(t) = w(t+T/2)$, defined on the time interval $(-T/2,T/2)$. According to time invariance of the wave equation, $u$ satisfies~\eqref{e:UC-presque-final} with $T/2$ in place of $T$ and we deduce that $u = 0$ in $\mathcal{D}'((-\eps, \eps)\times \M)$. That is to say $w = 0$ in $\mathcal{D}'((T/2-\eps, T/2+\eps)\times \M)$ but since $C^0([0,T]; H^1_0(\M))\cap C^1([0,T];L^2(\M))$ this implies $(w,\d_tw)(T/2)=0$ and well-posedness of the Boundary-value problem~\eqref{e:ondes-obs} implies $(w_0,w_1)=(0,0)$.
\enp
Recalling the discussion of Section~\ref{s:approx-wave-intro}, we now have, as a last corollary, the following.
\begin{corollary}
\label{c:approx-control}
Assume $(\M,g)$ is a compact connected Riemannian manifold with or without boundary $\d \M$, $q \in L^\infty(\Int(\M))$, let $\omega \subset \M$ be a nonempty open set,  and fix $T>2\mathcal{L}(\M ,\omega)$.
Then, Equation~\eqref{e:controlled-onde} is {\em approximately controllable} from $(\chi_\omega, T)$.
\end{corollary}
Again, on account to finite speed of propagation, the minimal time $\mathcal{L}(K ,\omega),\mathcal{L}(\M ,\omega),2\mathcal{L}(\M ,\omega)$ needed in Theorems~\ref{c:obserwaveintro},~\ref{c:global-UC} and Corollary~\ref{c:approx-control} are all optimal, see e.g.~\cite{Russell:71-1,Russell:71-2}.

\section{Notes}
\label{notes}
To conclude, we briefly discuss in this section the related question of obtaining quantitative estimates associated to unique continuation, historical notes and pointers to the literature.
\subsection{Quantitative unique continuation and the cost of approximate controls}
\label{sec:intro-wave-schrod}
All the results we presented in this course have their {\em quantitative } counterpart. We discuss a few of them in this section and only provide statements without any proofs (as opposed to the previous sections). The latter are rather technical and would not fit in these introductory lecture notes.

The local quantitative estimates associated to the H\"ormander Theorem~\ref{thmreal2} are of H\"older type. This was first noticed by Bahouri~\cite{Bahouri:87} and leads to the following type of statement. 

\begin{theorem}
\label{thmreal2quantif}
Under the assumptions of Theorem \ref{thmreal2}, there exist a neighborhood $V$ of $x_{0}$, $\alpha\in (0,1)$ and $C>0$ so that we have 
\bna
\nor{u}{H^{1}(V)}\leq C \nor{u}{H^{1}(\Omega \cap  \{\Psi \geq  \Psi(x_0)\})}^{\alpha}\nor{u}{H^{1}(\Omega)}^{1-\alpha}
\ena
 for all $u\in H^1(\Omega)$ solution of $Pu=0$.
\end{theorem}
It was also understood in~\cite{Robbiano:95,LR:95}, for elliptic operators, that this kind of interpolation estimates can be iterated to deduce global results. We also refer to \cite{LL:15} for a global geometric statement for general operators, across a global foliation of strongly pseudoconvex hypersurfaces. 

Concerning the unique continuation result described in Section \ref{chapterwave}, the quantitative stability estimates are of logarithmic type.
The following (global) result is the quantitative version of the global unique continuation result of Theorem~\ref{c:global-UC}, proved in~\cite{LL:15} (see also~\cite{LL:17Hypo} for obtaining the observation term in $L^2$ norm). It followed several earlier results on the subject and we refer to Section \ref{s:history} for some summary of the previous literature.
\begin{theorem}[Quantitative unique continuation for waves~\cite{LL:15}]
\label{thmobserwaveintro}
Let $\M$ be a compact Riemannian manifold with (or without) boundary. For any nonempty open subset $\omega$ of $\M$ and any $T> 2\mathcal{L}(\M ,\omega)$, there exist $C, \kappa ,\mu_0>0$ such that for any $(u_0,u_1)\in H^1_0(\M)\times L^2(\M)\setminus \{(0,0)\}$ and associated solution $u$ of 
\bneqn
\label{e:free-wave}
\partial_t^2u-\Delta_g u= 0& & \text{ in } (0,T) \times \Int(\M) , \\
u_{\left|\partial \M\right.}=0 & &  \text{ in } (0,T) \times \d \M ,\\
(u,\partial_tu)_{\left|t=0\right.}=(u_0,u_1)& & \text{ in } \Int(\M)  ,
\eneqn
we have
\begin{align}
\label{e:QUCP-mu}
\nor{(u_0,u_1)}{L^2\times H^{-1}} & \leq C e^{\kappa \mu}\nor{u}{L^2((0,T)\times\omega)}+\frac{1}{\mu}\nor{(u_0,u_1)}{H^1\times L^2}, \quad \text{for all } \mu\geq \mu_0, \\
\label{e:QUCP-log-form}
\nor{(u_0,u_1)}{L^2\times H^{-1}} & \leq C \frac{\nor{(u_0,u_1)}{H^1\times L^2}}{\log\left(1+\frac{\nor{(u_0,u_1)}{H^1\times L^2}}{\nor{u}{L^2((0,T)\times \omega)}}\right)} ,\\
\label{e:QUCP-Lambda-form}
\nor{(u_0,u_1)}{H^1\times L^2} & \leq C e^{\kappa \Lambda}\nor{u}{L^2((0,T)\times \omega)} , \quad \text{ with }\Lambda=\frac{\nor{(u_0,u_1)}{H^1\times L^2}}{\nor{(u_0,u_1)}{L^2\times H^{-1}}} .
\end{align}
\end{theorem}
The three inequalities~\eqref{e:QUCP-mu}--\eqref{e:QUCP-log-form}--\eqref{e:QUCP-Lambda-form} are actually equivalent (up to changing the values of the constants).
 Note that the statement~\eqref{e:QUCP-mu} remains valid for all $\mu>0$ (not only $\mu\geq \mu_0$), the estimate for $\mu$ bounded being trivial/useless.
In Estimate~\eqref{e:QUCP-log-form}, the function on the right hand-side is to be understood as being $\left(\log(1+\frac{1}{x})\right)^{-1}$ for $x>0$ and $0$ for $x=0$.
In Estimate~\eqref{e:QUCP-Lambda-form}, $\Lambda$ has to be considered as the typical frequency of the initial data. So, the estimate states a cost of observability of the order of an exponential of the typical frequency. As an illustration, taking for initial data $(u_0,u_1) = (\psi_\lambda , 0)$ with $\psi_\lambda$ a normalized eigenfunction of the Laplace-Dirichlet operator on $\M$, associated to the eigenvalue $\lambda$, one has $\Lambda \sim \sqrt{\lambda}$ and~\eqref{e:QUCP-Lambda-form} recovers the tunneling estimate $\nor{\psi_\lambda}{L^2(\omega)} \geq C^{-1} e^{-C\sqrt{\lambda}}$ (see~\cite{DF:88,LR:95,LeLe:09,LRLR:book1,LL:18} for a discussion on quantitative unique continuation for eigenfunctions of elliptic operators).  A generalization to hypoelliptic operators has been proved in~\cite{LL:17Hypo}.

As stated by Lebeau~\cite[Section~2, pages~5 and~6]{Leb:Analytic} in the analytic context, the  exponential (resp. logarithmic) dependence in Estimates~\eqref{e:QUCP-mu}--\eqref{e:QUCP-Lambda-form} (resp. in~\eqref{e:QUCP-log-form}) is sharp in  general:
 the form of the estimates of Theorem~\ref{thmobserwaveintro} is optimal if the geometric control condition of~\cite{RT:74,BLR:88,BLR:92} is violated.
 
As a consequence of Theorem~\ref{thmobserwaveintro}, we obtain the cost of the approximate controllability of the wave equation (i.e. a quantitative version of Corollary~\ref{c:approx-control}).
\begin{theorem}[Cost of approximate controls~\cite{LL:15}]
\label{t:cost-app-cont-wave}
Assume $\M$ is compact and connected. For any nonempty open set $\omega \subset \M$ and any $T> 2 \mathcal{L}(\M ,\omega)$, there exist $C,c>0$ such that for any $\eps >0$ and any $(u_0 ,u_1) \in H^1_0(\M) \times L^2(\M)$, there exists $f \in L^2((0,T) \times \omega)$ with 
$$
\|f\|_{L^2((0,T) \times \omega)} \leq C e^{\frac{c}{\eps}} \nor{(u_0 ,u_1)}{H^1_0(\M) \times L^2(\M)} ,
$$
such that the solution of~\eqref{e:controlled-onde}
satisfies $\nor{(u ,\d_t u)_{|t=T}}{L^2(\M) \times H^{-1}(\M)} \leq \eps\nor{(u_0 ,u_1)}{H^1_0(\M) \times L^2(\M)}$. 
\end{theorem}
An application of Theorem~\ref{thmobserwaveintro} to the exact controllability problem (and a uniform quantitative version of it, with dependence of the observability constants with respect to different parameters) is developed in~\cite{LL:16}.

\subsection{Historical remarks on unique continuation for waves}
\label{s:history}
The first general unique continuation result of the form~\eqref{e:UCP-general-local} is the Holmgren-John theorem~\ref{thmholmgren} (due to Holmgren~\cite{Holmgren} in a special case, and to John~\cite{John:49} in the general case). This local unique continuation result enjoys a global version proved by John~\cite{John:49}, where uniqueness is propagated through a family of noncharateristic hypersurfaces (see Section~\ref{s:semiglo} above). 

When focusing on operators with (only) smooth ($C^\infty$) coefficients, the first result is due to Carleman \cite{Carleman:39} who first had the idea to conjugate the operator with an exponential weight to obtain unique continuation. He proved the result in the case of elliptic operators of order $2$ in dimension $2$. Calder\'on \cite{Calderon:58} extended the result to operators with simple characteristics. Namely, to situations where $p_{\Psi}(x_0,\xi)=0 \implies \{p_{\Psi},\Psi\}(x_0,\xi)\neq0$. The most general result of Theorem~\ref{thmreal2} was proved by H\"ormander~\cite[Chapter~VIII]{Hoermander:63},~\cite[Chapter~XXVIII]{Hoermander:V4}, still using Carleman estimates. Uniqueness across a hypersurface holds assuming the strong pseudoconvexity condition (Definition~\ref{defconvexsurface}). 
Other works consider the limit case where higher order of cancelation are considered. We refer e.g. to~\cite{Zuily:83,Lerner:book-carleman} for more details on this topic.

Motivation to study the wave operator arised both from geoseismics~\cite {Symes:83}, control theory~\cite{Lio:88, Lio2:88} and inverse problems~\cite{Belishev:87,NO:22}. 
For the wave operator $P = \d_t^2 - \Delta_g$ on $\Omega = (-T,T) \times \M$, the central question remained for long the validity of the unique continuation property~\eqref{e:global-UC}.
If $\M$ is analytic (and connected), the Holmgren-John theorem~\ref{thmholmgren} applies, which together with the argument of John~\cite{John:49}, allows to prove unique continuation from $(-T,T) \times \omega$ for any nonempty open set $\omega$ as soon as $T>  \mathcal{L} (\M ,\omega)$.
Removing the analyticity condition on $\M$ has led to a considerable difficulty, since H\"ormander general uniqueness result does not apply in this setting and the Carleman-H\"ormander strategy fails, as showed by the Alinhac-Baouendi counterexamples of Theorems~\ref{t:alinhac} and~\ref{t:alinhac-baouendi}.

This uniqueness problem in the $C^\infty$ setting was first solved by Rauch-Taylor~\cite{RT:72} and Lerner~\cite{Lerner:88} in the case $T= \infty$, and $\M = \R^d$ (under different assumptions at infinity). Then, Robbiano~\cite{Robbiano:91} managed to prove that unique continuation from $(-T,T) \times \omega$ holds in any domain $\M$ as soon as $\omega \neq \emptyset$ and $T \geq C_0 \mathcal{L}(\M ,\omega)$, with $C_0$ sufficiently large (a result qualified as ``striking and surprising'' by H\"ormander~\cite{Hor:97}). H\"ormander~\cite{Hormander:92} improved this result down to $T> \sqrt{\frac{27}{23}} \mathcal{L}(\M ,\omega)$. That these two results fail to hold in time $\mathcal{L}$ translates the fact that the local uniqueness results of these two authors are not valid across any noncharacteristic hypersurface.
The proof of local uniqueness results across any noncharacteristic hypersurface for $\d_t^2 - \Delta_g$ was reached by Tataru in~\cite{Tataru:95} (Theorem~\ref{thmwave}), leading to the global unique continuation result in optimal time $T> \mathcal{L}(\M ,\omega)$ (Theorem~\ref{t:region-dependence} and Corollary~\ref{c:obserwaveintro}). The result of Tataru is not restricted to the wave operator:
it holds for operators with coefficients that are analytic in part of the variables, interpolating between the Holmgren theorem and the H\"ormander theorem. Technical assumptions of this article were successively removed by Robbiano-Zuily~\cite{RZ:98}, H\"ormander~\cite{Hor:97} and Tataru~\cite{Tataru:99},  leading to a very general local unique continuation result for operators with partially analytic coefficients (containing as particular cases both Holmgren and H\"ormander theorems).

The local, semiglobal and global quantitative versions of all these unique continuation results  are proved in~\cite{LL:15}.
For the particular case of  wave equations with time-independent coefficients, a global quantitative statement takes the form of Theorem~\ref{thmobserwaveintro}.
In the analytic setting, this result is a global quantitative version of the Holmgren-John theorem and can be proved with the theory developed by Lebeau in~\cite{Leb:Analytic}.
In the $C^\infty$ case, Robbiano~\cite{Robbiano:95} first proved a similar inequality for $T$ sufficiently large and $Ce^{\kappa \mu}$ replaced by $Ce^{\kappa \mu^2}$, improved by Phung \cite{Phung:10} to $C_\eps e^{\kappa_\eps \mu^{1+\eps}}$. In~\cite{Tatarunotes}, Tataru suggested a strategy to obtain $C_\eps e^{\kappa_\eps \mu^{1+\eps}}$ in optimal time (in domains without boundaries). At the same time we proved the above Theorem~\ref{thmobserwaveintro}, Bosi, Kurylev and Lassas~\cite{BKL:16} obtained a related result with $C_\eps e^{\kappa_\eps \mu^{1+\eps}}$ dependence, on manifolds without boundaries.

\subsection{Related references}
\paragraph{Treatises on unique continuation.}
There are many references on  unique continuation and the related topics of Carleman estimates.  
We mention here a nonexhaustive list of books and monographs treating these topics.
The classical reference on unique continuation for partial differential operators is the Chapter XXVIII of Lars H\"ormander's treatise \cite{Hoermander:V4}. The latter gives a more general framework for what is described in Chapter~\ref{chapterclassical}. 
The recent book of Nicolas Lerner~\cite{Lerner:book-carleman} also contains several results covered in the present notes.
A short lecture notes version of this book can be found in~\cite{Lernernotes}.
We also refer to the book of Claude Zuily \cite{Zuily:83} and the notes of Daniel Tataru~\cite{Tatarunotes} for related questions. Carleman estimates for hyperbolic equations are also studied in \cite{BYBook:17} with a special emphasis on applications to inverse problems.
 The presentation of Chapter~\ref{chapterwave}, concerning the wave operator, is inspired by the article~\cite{Hor:97}.
 Finally, the survey article by J\'er\^ome Le Rousseau and Gilles Lebeau~\cite{LeLe:09} is a smooth introduction to unique continuation for elliptic operators (which are only alluded here in Remark~\ref{r:elliptic-rem}). See also the recent books of J\'er\^ome Le Rousseau, Gilles Lebeau and Luc Robbiano~\cite{LRLR:book1,LRLR:book2} for an extensive treatment of Carleman estimates and unique continuation for elliptic operators, with plenty of applications.

\paragraph{Some other recent developments for waves.}
We briefly discuss recent developments in the field of unique continuation and Carleman estimates for waves, that are not alluded in Section~\ref{s:history}. The list is of course very far from exhaustive. 
Most of these developments concern elaborations on global Carleman estimates for waves in a (global) setting close to that of Theorem~\ref{thmreal2}, with applications to control theory, inverse problems or general relativity.

Global Carleman estimates for waves were proved in~\cite[Chapter~4]{FI:96} and \cite{TY:02} with applications to controllability.
Various boundary conditions have been considered in~\cite{IY:00}.  Global Carleman estimates with limited regularity coefficients were investigated in~\cite{Imanuvilov:02}, with applications to controllability and inverse problems (determining lower order terms in the equation). 

Admitting lower order terms in unique continuation results is often crucial for treating nonlinear problems. For instance, a nonlinearity of the form $u^p$ can be treated as a term $q u$ with potential $q=u^{p-1}$, in general having limited regularity. A global unique continuation statement  (close to Proposition~\ref{p:global-pseudoconvex} above) was proved in~\cite{Ruiz:92}, with application to energy decay for nonlinear waves. Sometimes having nonlinear problems in mind, it has therefore been a goal to minimize the regularity of the admissible lower order terms. This led to some ``dispersive'' Carleman estimates with Strichartz type spaces.  The literature is vast, and we refer for instance to~\cite{KRS:87,DDSF:05, KT:05} and the references therein.

Global unique continuation results for nonlinear waves have been proved in~\cite{JL:13,JL:20}, taking advantage of the nonlinear term. The global geometric assumptions, like Geometric Control Condition or weak trapping, allows to prove the analyticity in time of the nonlinear solution. Results in the spirit of Section~\ref{chapterwave} can then be used to conclude the unique continuation.

Another problem not covered by these notes is the so called ``strong unique continuation'' property. Instead of assuming that the solution $u$ is zero on an open set $\omega$, one assumes infinite order of vanishing at a point or along a submanifold. We refer for instance to \cite{L:99} for the case of unique continuation for waves from $[0,T]\times \{x_0\}$ and a quantified version in \cite{V:17}. 

Numerical analysis and algorithms were introduced in~\cite{BDBE:13} based on Carleman estimates. 
Carleman estimates in a geometric context (extending Proposition~\ref{p:global-pseudoconvex} above) were proved in~\cite{DZZ:08,Shao:19,JS:21}, with application to observability estimates and null-controllability.
Recent developments involving Carleman estimates for waves to solve inverse problems include~\cite{BDBEO:21} (for developing recovery algorithms) or~\cite{AFO:22,AFO:23} (adapting the Boundary Control method to a geometric setting in which to apply the H\"ormander theorem~\ref{thmreal2}).
In~\cite{LR:15}, the authors keeps track of two large parameters (here $\tau$ and $\lambda$) in Carleman estimates in the context of the (general version of the) H\"ormander Theorem~\ref{thmreal2}.

Unique continuation problems for wave operators also arise from mathematical general relativity. 
In~\cite{IK:09,Lerner:19}, the authors proved a unique continuation statement through tranversally intersecting characteristic hypersurfaces, relying on the H\"ormander theorem~\ref{thmreal2}. See also~\cite{AS:15,CS:22} where unique continuation statements for waves from infinity are proved.
 
 The quantitative Theorem~\ref{thmobserwaveintro} was generalized in~\cite{LL:17Hypo} to hypoelliptic waves, that is to say in which the elliptic operator $\Delta_g$ is replaced by a H\"ormander sum of squares of vector fields. See also the review article~\cite{LL:XEDP}. Finally, Theorem~\ref{thmobserwaveintro} was also generalized recently in~\cite{F:22} to the case in which the metric $g$ has a jump across a hypersurface (a situation arising e.g. in geosismics).

\bigskip
\noindent
{\textbf{\em Acknowledgements.}}
These notes correspond to (and expand slightly) a course given at CIRM by the second author in November 2022 for the SMF School ``Spectral Theory, Control and Inverse Problems''.
The authors would like to warmly thank Lucie Baudouin, Franck Boyer, J\'er\'emi Dard\'e, Sylvain Ervedoza and Julien Royer for having organized this very nice conference.
These notes also originate from a (part of a) M2 course (see~\cite{LL:book} for an expanded version of the notes) the first author taught in Sorbonne Universit\'e (2015--2017) and the second author taught in Universit\'e Paris-Saclay (2017--2020).  
The authors would like to thank all students and colleagues who followed one (or even several) of these courses, for their questions and comments which allowed to improve the presentation of the notes.
The second author is partially supported by the Agence Nationale de la Recherche under grants SALVE (ANR-19-CE40-0004) and ADYCT (ANR-20-CE40-0017).

\small
\bibliographystyle{alpha} 
\bibliography{bibli}
\end{document}

%% file: mymacros.tex
\newtheorem{lemma}{Lemma}[section]
\newtheorem{theorem}[lemma]{Theorem}
\newtheorem{proposition}[lemma]{Proposition}
\newtheorem{corollary}[lemma]{Corollary}

\theoremstyle{definition}
\newtheorem{example}[lemma]{Example}
\newtheorem{definition}[lemma]{Definition}
\newtheorem{remark}[lemma]{Remark}

\makeatletter
\def\keywords{
    \vspace{1ex}
    \noindent
    \if@twocolumn
      \small{\bf  Keywords}\/---$\!$    \else
      \begin{center}\small\ {\bf Keywords}\end{center}\quotation\small
    \fi}
\def\endkeywords{\vspace{0.6em}\par\if@twocolumn\else\endquotation\fi
    \normalsize\rm}
\makeatother

\newcommand{\calS}{\ensuremath{\mathcal S}}

\newcommand{\x}{\ensuremath{\bf x}}

\DeclareMathOperator{\loc}{loc}
\DeclareMathOperator{\comp}{comp}

\DeclareMathOperator{\Vol}{Vol}

\DeclareMathOperator{\Int}{Int}

\DeclareMathOperator{\Hess}{Hess}
\DeclareMathOperator{\range}{range}

\newcommand{\mb}[1]{\ensuremath{\mathbb{#1}}}
\newcommand{\N}{{\mb{N}}}

\newcommand{\R}{{\mb{R}}}
\newcommand{\C}{{\mb{C}}}

\newcommand{\eps}{\varepsilon}

\newcommand{\M}{\ensuremath{\mathcal M}}

\let \Re \relax
\DeclareMathOperator{\Re}{Re}
\let \Im \relax
\DeclareMathOperator{\Im}{Im}

\newcommand{\ovl}[1]{\overline{#1}}

\renewcommand{\S}{\ensuremath{\mathbb S}}

\DeclareMathOperator{\supp}{supp}

\DeclareMathOperator{\diag}{diag}
\DeclareMathOperator{\dist}{dist}

\newcommand{\transp}{\ensuremath{\phantom{}^{t}}}

\renewcommand{\d}{\ensuremath{\partial}}

\newcommand{\dsp}{\displaystyle}

\let \div \relax
\DeclareMathOperator{\div}{div}

\DeclareMathOperator{\length}{length}

\def\x{x}

\newcommand{\F}{\mathscr F}

\newcommand{\todo}[1]{$\clubsuit$ {\tt #1}}